
\def\LaTeX{%
  \let\Begin\begin
  \let\End\end
  \let\salta\relax
  \let\finqui\relax
  \let\futuro\relax}

\def\UK{\def\our{our}\let\sz s}
\def\USA{\def\our{or}\let\sz z}



\LaTeX

\USA


\salta

\documentclass[twoside,12pt]{article}
\setlength{\textheight}{24cm}
\setlength{\textwidth}{16cm}
\setlength{\oddsidemargin}{2mm}
\setlength{\evensidemargin}{2mm}
\setlength{\topmargin}{-15mm}
\parskip2mm


\usepackage{amsmath}
\usepackage{amsthm}
\usepackage{amssymb}
\usepackage{hyperref}
\usepackage[mathcal]{euscript}

\usepackage[usenames,dvipsnames]{color}
%
%


\definecolor{viola}{rgb}{0.3,0,0.7}
\definecolor{ciclamino}{rgb}{0.5,0,0.5}

\def\pcol #1{{#1}}
\def\ele #1{{#1}}
\def\pcgg #1{{#1}}
\def\elena #1{{#1}}
\def\pier #1{{#1}}
\def\gianni #1{{#1}}


\bibliographystyle{plain}


%

\finqui

\def\Beq{\Begin{equation}}
\def\Eeq{\End{equation}}
\def\Bsist{\Begin{eqnarray}}
\def\Esist{\End{eqnarray}}
\def\Bdef{\Begin{definition}}
\def\Edef{\End{definition}}
\def\Bthm{\Begin{theorem}}
\def\Ethm{\End{theorem}}
\def\Blem{\Begin{lemma}}
\def\Elem{\End{lemma}}
\def\Bprop{\Begin{proposition}}
\def\Eprop{\End{proposition}}

\def\Brem{\Begin{remark}\rm}
\def\Erem{\End{remark}}

\def\Bnot{\Begin{notation}\rm}
\def\Enot{\End{notation}}
\def\Bdim{\Begin{proof}}
\def\Edim{\End{proof}}
\let\non\nonumber




\def\step #1 \par{\medskip\noindent{\bf #1.}\quad}


\def\Lip{Lip\-schitz}
\def\holder{H\"older}
\def\aand{\quad\hbox{and}\quad}

\def\lsc{lower semicontinuous}

\def\lhs{left hand side}
\def\rhs{right hand side}
\def\sfw{straightforward}


\def\characteriz{characteri\sz}
\def\discretiz{discreti\sz}
\def\generaliz{generali\sz}

\def\organiz{organi\sz}

\def\regulariz{regulari\sz}

\def\bhv{behavi\our}


\let\fontt\tt  

\def\multibold #1{\def\arg{#1}%
  \ifx\arg\pto \let\next\relax
  \else
  \def\next{\expandafter
    \def\csname #1#1#1\endcsname{{\bf #1}}%
    \multibold}%
  \fi \next}

\def\multislantedbold #1{\def\arg{#1}%
  \ifx\arg\pto \let\next\relax
  \else
  \def\next{\expandafter
    \def\csname #1#1\endcsname{\hbox{\slantedbf #1}}%
    \multislantedbold}%
  \fi \next}

\def\pto{.}

\def\multical #1{\def\arg{#1}%
  \ifx\arg\pto \let\next\relax
  \else
  \def\next{\expandafter
    \def\csname cal#1\endcsname{{\cal #1}}%
    \multical}%
  \fi \next}


\def\multimathop #1 {\def\arg{#1}%
  \ifx\arg\pto \let\next\relax
  \else
  \def\next{\expandafter
    \def\csname #1\endcsname{\mathop{\rm #1}\nolimits}%
    \multimathop}%
  \fi \next}

\multibold
qwertyuiopasdfghjklzxcvbnmQWERTYUIOPASDFGHJKLZXCVBNM.

\multislantedbold
qwertyuiopasdfghjklzxcvbnmQWERTYUIOPASDFGHJKLZXCVBNM.

\multical
QWERTYUIOPASDFGHJKLZXCVBNM.

\multimathop
dist div dom meas sign supp .


\def\accorpa #1#2{\eqref{#1}--\eqref{#2}}
\def\Accorpa #1#2 #3 {\gdef #1{\eqref{#2}--\eqref{#3}}%
  \wlog{}\wlog{\string #1 -> #2 - #3}\wlog{}}


\def\QED{\hfill $\square$}
\def\spacca{\noalign{\allowbreak}}

\def\supess{\mathop{\rm sup\,ess}}
\def\Neto{\mathrel{{\scriptstyle\nearrow}}}
\def\Seto{\mathrel{{\scriptstyle\searrow}}}

\def\seto{\mathrel{{\scriptscriptstyle\searrow}}}

\def\somma #1#2#3{\sum_{#1=#2}^{#3}}
\def\tonde #1{\left(#1\right)}

\def\graffe #1{\mathopen\{#1\mathclose\}}
\def\Graffe #1{\left\{#1\right\}}
\def\<#1>{\mathopen\langle #1\mathclose\rangle}
\def\norma #1{\mathopen \| #1\mathclose \|}

\def\iot {\int_0^t}
\def\ioT {\int_0^T}
\def\iO{\int_\Omega}
\def\iG{\int_\Gamma}
\def\intQt{\int_{Q_t}}
\def\intQ{\int_Q}

\def\dt{\partial_t}
\def\ds{\,ds}

\def\cpto{\,\cdot\,}

\def\checkmmode #1{\relax\ifmmode\hbox{#1}\else{#1}\fi}
\def\aeO{\checkmmode{a.e.\ in~$\Omega$}}
\def\aeQ{\checkmmode{a.e.\ in~$Q$}}

\def\aaO{\checkmmode{for a.a.~$x\in\Omega$}}

\def\aat{\checkmmode{for a.a.~$t\in(0,T)$}}
\def\aeT{\checkmmode{a.e.\ in~$(0,T)$}}


\def\erre{{\mathbb{R}}}




\def\genspazio #1#2#3#4#5{#1^{#2}(#5,#4;#3)}
\def\spazio #1#2#3{\genspazio {#1}{#2}{#3}T0}

\def\L {\spazio L}
\def\H {\spazio H}
\def\W {\spazio W}

\def\C #1#2{C^{#1}([0,T];#2)}

\def\Wz{W_0}
\def\Vp{V^*}


\def\Lx #1{L^{#1}(\Omega)}
\def\Hx #1{H^{#1}(\Omega)}
\def\Wx #1{W^{#1}(\Omega)}
\def\Cx #1{C^{#1}(\overline\Omega)}
\def\LxG #1{L^{#1}(\Gamma)}

\def\Luno{\Lx 1}
\def\Ldue{\Lx 2}
\def\Linfty{\Lx\infty}

\def\Huno{\Hx 1}
\def\Hdue{\Hx 2}


\def\LQ #1{L^{#1}(Q)}
\def\CQ {C^0(\overline Q)}


\let\theta\vartheta
\let\eps\varepsilon

\let\TeXchi\chi                         
\newbox\chibox
\setbox0 \hbox{\mathsurround0pt $\TeXchi$}
\setbox\chibox \hbox{\raise\dp0 \box 0 }
\def\chi{\copy\chibox}


\let\hat\widehat\let\tilde\widetilde

\def\normaV #1{\norma{#1}_V}
\def\normaH #1{\norma{#1}_H}
\def\normaW #1{\norma{#1}_W}

\def\normaZ #1{\norma{#1}_Z}

\let\TeXnu\nu
\font\slantedbf=cmmib10 scaled \magstep 1    
\def\nu{\hbox{\slantedbf\char 23}}
\def\sigma{\hbox{\slantedbf\char 27}}
\def\dn{\partial_\TeXnu}

\def\Funo{F_1}
\def\Fdue{F_2}
\def\thetac{\theta_{\!c}}
\def\RO{R_\Omega}
\def\BO{B_\Omega}
\def\bG{b_\Gamma}
\def\BG{B_\Gamma}
\def\thetaz{\theta_0}
\def\uz{u_0}
\def\uzp{\uz'}
\def\chiz{\chi_0}
\def\wz{w_0}
\def\sigmaz{\sigma_0}

\def\Bz{B_0}

\def\chiuno{\chi_1}

\def\umu{u_{-1}}

\def\ceps{c_\eps}
\def\cdelta{c_\lambda}
\def\cz{c_0}
\def\deltaz{\lambda_0}

\def\deltastar{\lambda_*}

\def\cstar{C_*}

\def\cq{c_q}
\def\gammastar{\gamma_*}
\def\jV{j_V}
\def\jH{j_H}
\def\JV{J_V}

\def\vstar{v^*}

\def\alphaeps{\alpha_{\eps}}
\def\auno{a_1}
\def\adue{a_2}
\def\leps{\lambda_\eps}
\def\Deps{D_\eps}
\def\Zeps{\zeta_\eps}
\def\betaeps{\beta_\eps}
\def\thetaeps{\theta_\eps}
\def\ueps{u_\eps}
\def\chieps{\chi_\eps}
\def\sigmaeps{\sigma_\eps}
\def\weps{w_\eps}

\def\phieps{\phi_\eps}
\def\Feps{F_\eps}
\def\feps{f_\eps}
\def\Funoeps{F_{1,\eps}}
\def\Funoepsp{F_{1,\eps'}}
\def\Geps{G_\eps}
\def\Gpeps{\Geps'}
\def\geps{\gamma_\eps}
\def\gpeps{\geps^{\pier{\,\prime}}}
\def\Reps{R_\eps}        \let\Reps R
\def\BOeps{B_{\Omega,\eps}}
\def\BGeps{B_{\Gamma,\eps}}
\def\Beps{B_\eps}
\def\thetazeps{\theta_{0,\eps}}  
\def\uzeps{u_{0,\eps}}
\def\uzepsp{\uzeps'}
\def\chizeps{\chi_{0,\eps}}
\def\wzeps{w_{0,\eps}}
\def\sigmazeps{\sigma_{0,\eps}}
\def\thetam{\theta_m}
\def\thetan{\theta_n}
\def\thetamp{\theta_{m+1}}
\def\thetanp{\theta_{n+1}}

\def\un{u_n}
\def\um{u_m}
\def\ump{u_{m+1}}
\def\unp{u_{n+1}}
\def\unm{u_{n-1}}
\def\chim{\chi_m}
\def\chimp{\chi_{m+1}}
\def\chin{\chi_n}
\def\chinp{\chi_{n+1}}

\def\sigman{\sigma_n}
\def\sigmanp{\sigma_{n+1}}

\def\sigmamp{\sigma_{m+1}}

\def\an{\alphaeps(\thetan)}

\def\apn{\alphaeps'(\thetan)}

\def\Gn{\Geps(\chin)}

\def\Gpn{\Geps'(\chin)}

\def\gn{\geps(\chin)}

\def\gpn{\geps^{\pier{\,\prime}}(\chin)}
\def\gpnp{\geps^{\pier{\,\prime}}(\chinp)}

\def\Rn{R_n}

\def\Bn{B_n}

\def\zn{z_n}
\def\znp{z_{n+1}}
\def\znpp{z_{n+2}}
\def\znm{z_{n-1}}
\def\zm{z_m}
\def\zmp{z_{m+1}}
\def\zN{z_N}
\def\zNp{z_{N+1}}
\def\zNm{z_{N-1}}
\def\vn{v_n}
\def\vnp{v_{n+1}}
\def\etam{\eta_m}
\def\etan{\eta_n}

\def\etanm{\eta_{n-1}}

\def\zetan{\zeta_n}
\def\zetamp{\zeta_{m+1}}
\def\zetanp{\zeta_{n+1}}

\def\fn{f_n}
\def\fm{f_m}
\def\fz{f_0}
\def\gibbs{{\cal G}}

\def\dtau{\delta_\tau}
\def\dzn{\dtau z_n}

\def\overz{\overline z_\tau}
\def\underz{\underline z_\tau}
\def\hv{\hat v_\tau}
\def\hz{\hat z_\tau}
\def\tz{\tilde z_\tau}
\def\underv{\underline v_\tau}
\def\hv{\hat v_\tau}

\def\underR{\underline R_\tau}

\def\underB{\underline B_\tau}

\def\overtheta{\overline\theta_\tau}
\def\undertheta{\underline\theta_\tau}
\def\htheta{\hat\theta_\tau}
\def\overu{\overline u_\tau}
\def\underu{\underline u_\tau}
\def\hu{\hat u_\tau}
\def\tu{\tilde u_\tau}
\def\overchi{\overline\chi_\tau}
\def\underchi{\underline\chi_\tau}
\def\hchi{\hat\chi_\tau}
\def\tchi{\tilde\chi_\tau}
\def\oversigma{\overline\sigma_\tau}
\def\undersigma{\underline\sigma_\tau}
\def\hsigma{\hat\sigma_\tau}

\def\Tk{{\cal T}_k}
\def\Qk{Q^k}
\def\intQk{\int_{\Qk}}

\def\vdelta{v_\pier{\lambda}}
\def\Deltastar{\Delta^{\!*}}



\Begin{document}


\title{{\bf \pier{Existence of solutions for}\\ a mathematical model related\\ to solid-solid 
phase transitions\\ in shape memory alloys}\footnote{{\bf Acknowledgments.}\quad\rm
The present paper benefits from the support of the MIUR-PRIN Grant 2010A2TFX2 ``Calculus of variations''.}}
\author{}
\date{}
\maketitle
\Begin{center}
\vskip-2cm
{\large\bf Elena Bonetti$^{(1)}$}\\
{\normalsize e-mail: {\fontt elena.bonetti@unipv.it}}\\[.2cm]
{\large\bf Pierluigi Colli$^{(1)}$}\\
{\normalsize e-mail: {\fontt pierluigi.colli@unipv.it}}\\[.2cm]
{\large\bf Mauro Fabrizio$^{(2)}$}\\
{\normalsize e-mail: {\fontt fabrizio@dm.unibo.it}}\\[.2cm]
{\large\bf Gianni Gilardi$^{(1)}$}\\
{\normalsize e-mail: {\fontt gianni.gilardi@unipv.it}}\\[.4cm]
$^{(1)}$
{\small Dipartimento di Matematica ``F. Casorati'', Universit\`a di Pavia}\\
{\small via Ferrata 1, 27100 Pavia, Italy}\\[.2cm]
$^{(2)}$
{\small Dipartimento di Matematica, Universit\`a di Bologna}\\
{\small piazza di Porta San Donato, 5, 40126 Bologna, Italy}\\[.8cm]
\End{center}

\Begin{abstract}
\elena{We consider a strongly nonlinear PDE system describing 
solid-solid phase transitions in shape memory alloys. 
The system accounts for the evolution of an order parameter~$\chi$ 
(related to different symmetries of the crystal lattice in the phase configurations), 
of the stress (and the displacement~$u$), 
and of the absolute temperature~$\theta$. 
The resulting equations present several technical difficulties to be \pier{tackled}\pier{: 
in particular, we emphasize} the presence of nonlinear coupling terms, 
higher order dissipative contributions, possibly multivalued operators. 
\pier{As for the evolution of temperature, a highly nonlinear parabolic equation 
has to be solved for a right hand side that is controlled only in~$L^1$.}
We prove \pier{the} existence of a solution for a \pcgg{regularized}  version, 
by use of a time discretization technique. 
Then, we perform suitable a priori estimates \pier{which 
allow us pass to the limit and find a weak \pier{global-in-time} solution to the system.}}
\\[2mm]
\\[2mm]
\\[5mm]
\End{abstract}

\salta

\pagestyle{myheadings}
\newcommand\testopari{\sc Bonetti \ --- \ Colli \ --- \ Fabrizio \ --- \ Gilardi}
\newcommand\testodispari{\sc Study of a mathematical model}
\markboth{\testodispari}{\testopari}

\finqui


\section{Introduction}
\label{Intro}
\setcounter{equation}{0}

\elena{This paper deals with a strongly nonlinear differential system\ele{,}
which may be related to austenite-martensite phase transitions in shape memory alloys.
These materials are characterized by the fact that they can be permanently deformed by mechanical loads and
then recover their original shape just by heating.
This phenomenon is justified by  a change of symmetry of the mesoscopic
structure, as the transition involves a deformation of the crystalline cells. 
In particular, the austenite phase (which is present at high temperatures) 
is more symmetric with respect to the martensite variants.
The model we are considering (see \cite{bfg} and \cite{dfg} for a detailed derivation) couples a Ginzburg-Landau type equation,
which describes the evolution of a phase (order) parameter~$\chi$, 
with the momentum balance (accounting for \gianni{accelerations}) 
written in the displacement~$u$, and the energy balance governing the evolution of the absolute temperature~$\theta$. 
Note that here, just for the \pier{sakes of simplicity and better readability of the paper, we let} the displacement be a scalar variable. 
As a consequence, deformations are accounted for by $\nabla u$ and the stress is a vector. 
In the more general situation (but in the small strain regime) 
deformations should be described by the linearized symmetric strain tensor.}

\elena{Here is the resulting PDE system:}
\Bsist
  && \bigl(\cz - \theta \alpha''(\theta) G(\chi) \bigr) \dt\theta
  - \theta \alpha'(\theta) G'(\chi) \dt\chi
  -\Delta\theta
  = \RO + |\dt\chi|^2
  \label{Iprima}
  \\
  && \dt^2 u - \div\sigma = \BO
  \quad \hbox{where} \quad
  \sigma \pier{{}={}} 
  \nabla u - \gamma(\chi)\ee
  \label{Iseconda}
  \\
  && \dt\chi - \Delta\chi
  + \thetac F'(\chi) + \alpha(\theta) G'(\chi) - \sigma\cdot\ee \, \gamma^{\pier{\,\prime}}(\chi) = 0
  \label{Iterza}
\Esist
\Accorpa\Ipbl Iprima Iterza
in the unknown fields $\theta$, $u$, and~$\chi$, 
\elena{with $\sigma$ denoting the stress.} 
\gianni{As usual,}  
the partial differential equations are meant to hold in \elena{a} bounded domain $\Omega\subset\erre^3$
and in some time interval~$(0,T)$.
In the above equations, $\cz$, $\kappa$ and $\thetac$ are positive constants,
$\ee$~is a fixed unit vector,
and $\alpha$, $F$, $G$\elena{,} and $\gamma$ are given nonlinear functions
satisfying suitable properties
(which \elena{in particular} ensure a parabolic character for~\eqref{Iprima}).
\pier{One may think to $F$ as a potential with two wells located for 
instance in $-1$ and~$+1$; $G$~is a nonnegative, bounded \ele{function} 
such that $G(0)= G'(0) =0$;} \elena{moreover, the function $\gamma$ in \eqref{Iterza} is related to $G$~by
\Beq
  \pier{\gamma(r) = G(r) \sign (r), \quad r\in \erre}
  \label{defgamma}
\Eeq
with \pier{$\sign(r)$ taking the values: $+1$ if $r>0$\ele{,} \  $0$ if $r=0$\ele{,}} \ $-1$ if $r<0$. 
\pier{We point out that both $G$ and $\gamma$ are sufficiently smooth\ele{:  
 for} their precise regularity we refer to the subsequent 
assumptions}  \eqref{hpG}, \eqref{hpbdd}, and~\eqref{hpgamma}. 
\elena{Actually, with respect to the model introduced in \cite{bfg} and~\cite{dfg}, 
we are \pier{taking} a smoother function $\alpha$ in the energy functional. 
Indeed, \pier{in \cite{bfg} and \cite{dfg} it is postulated that $\alpha$ is simply of the  
type $\alpha(r)=(r -\theta_M)^+$ for $r\in\erre$ ($\theta_M >0$~being a critical transition 
temperature and $(\,\cdot\,)^+ $ denoting the positive part function), 
which entails the embarrassing presence of a Dirac measure in the equation corresponding to~\eqref{Iprima}. 
Here, instead\gianni,} we consider $\alpha$ smooth 
(to~give a meaning to $\alpha''$ entering the definition of the specific heat),
\gianni{and bounded (for technical reasons)}. 
We point out that \pier{the boundedness} assumption is, at the end, not restricting from a modeling point of view, 
as it preserves the required \bhv\ between different phases 
and \pier{corresponds to a} change in the free energy (preserving minima) just for very high temperatures, 
when \pier{only} the austenite phase may be present. 
On the contrary, while in some classical models for (solid-solid) phase transitions (and in~\cite{bfg}),  
$F$~is just a \pier{quartic double-well potential}, \ele{here}
we include the possibility that $F$ accounts for internal (non-smooth) constraints \pier{on} the phase variable. 
In particular, equation \eqref{Iterza} has to be read as a differential inclusion
if the monotone part of $F'$ is replaced by a subdifferential 
(e.g., \pier{of the indicator function of the interval~$[-1,1]$, 
so that $\chi$ is compelled to take} values in a physically consistent range).}
Finally, $\RO$~and $\BO$ are given forcing terms.}

\elena{The system \pier{\accorpa{Iprima}{Iterza}} is then complemented with the proper initial conditions
\Beq
  \theta(0) = \thetaz , \quad u(0) = \uz, \quad \dt u(0) = \uz', \quad \chi(0) = \chiz 
  \label{Ibczero}
\Eeq
as well as} the boundary conditions for the fluxes, namely
\Beq
  \dn\theta = 0, \quad \sigma\cdot\nu = \bG \,,  \quad \dn\chi = 0
  \quad \hbox{on \pier{$\Gamma \times (0,T) $}}
  \label{Ibc}
\Eeq
\pier{where $\Gamma:=\partial\Omega$, $\nu$ denotes the outward normal unit vector on~$\Gamma$,  
$\dn:=\nu\cdot\nabla$ stands for the normal derivative
and $\bG$ is a given datum on the boundary}.

\pcgg{In this paper, we are
mainly interested in the analytical study of the initial-boundary value problem, 
which represents an interesting mathematical issue in itself. Thus, before proceeding, we briefly comment on the main difficulties we are going to deal with.}

\pcgg{First, let us point out the presence of the nonlinear coefficient of $\dt\theta$ in~\eqref{Iprima}, 
as well as of other nonlinear terms. 
In particular, the quadratic dissipative term on the right hand side 
of \eqref{Iprima} has to be handled and
is, a priori, estimated just in some $L^1$ space 
(once the time derivative of $\chi$ is estimated, as expected, in $L^2$ from \eqref{Iterza}). 
Hence, some \emph{ad~hoc} techniques for equations with $L^1$ data have to be applied.}

\pcgg{We also notice the presence of a non-smooth and possibly multivalued 
operator in~\eqref{Iterza}. \pcol{As for \eqref{Iterza}, our approach is very general and gives us the possibility to set} some internal constraint on the phase without using any a posteriori maximum principle type technique. On the other hand, it is clear that the treatment of a possibly singular and multivalued operator leads to additional mathematical difficulties.}

\pcgg{Next, we point out the presence of the inertial term $ \dt^2 u$ in  
equation~\eqref{Iseconda} which \pcol{is evolutionary and} hyperbolic. 
\pcol{The coupling of \eqref{Iseconda} with other equations \eqref{Iprima} and \eqref{Iterza} and with conditions \eqref{Ibczero}--\eqref{Ibc}} provides an absolutely non-trivial problem.}

\pcgg{Furthermore, even though some formal a priori estimates could be shown
with rather standard techniques,
the necessity of dealing with approximating problems makes the whole argument
difficult.
In particular, the precise choice of the \regulariz ation we make
is crucial and its construction is necessarily involved.
Eventually, such an approximating problem 
still couples a hyperbolic equation with two strongly nonlinear 
equations of parabolic type, whence its solvability is not obvious.
This forces us to \pcol{additionally} use a time \discretiz ation technique,
with turns out to be rather heavy.}

\ele{Concerning the physical meaning of the system under investigation, at first 
we recall that} \elena{several models describing austenite-martensite phase transitions 
have been introduced in the last years 
(see, among the others, \pier{\cite{Fremond, frem2, Rou}} and references therein). 
In~this paper we mainly refer to the Ginzburg-Landau theory describing changes in the internal order structure of the material. 
One of the main advantage of this  approach consists in viewing the phase transition as
a change of the order in the symmetry of the alloys, 
so that just one phase parameter is used instead of vectorial or tensorial parameters 
(see \cite{Fremond} and \cite{Auricchio},~\cite{AuricchioB}). 
The fact that the equation for the phase is scalar represents also a good point for numerical implementation. 
More precisely, we let $\chi$ describe the order structure, 
i.e.\gianni, $\chi=0$ stands for the presence of austenite, 
while different (oriented) variants of martensite are associated to \gianni{$\chi=+1$ and $\chi=-1$}. \ele{Some recent papers deal with this
kind of problem. Let us mention \pcgg{\cite{FaPe} and \cite{bfg2, dfmz}:  
the former is concerned with} a  model for shape memory alloys \pcgg{characterized 
by an intermediate pattern between first and second order phase transitions;
the other two papers focus on} histeretic effects in the solid-solid 
phase transition both for the 1D and 3D cases.}}

Let us point out that in the set $\{\chi=0\}$ \gianni{equation} \eqref{Iseconda} 
postulates an elastic behaviour of the material (as it is $\gamma(0)=0$), 
while if \gianni{$\chi=\pm1$} a transformation stress appears, whose direction 
depends on the orientation of the martensitic variant. 

\elena{Equations are recovered by balance laws and thermodynamic principles by virtue of the following Gibbs free energy functional, 
depending on the state variables $\chi$, $\nabla\chi$, $\theta$, and the stress vector~$\sigma$:
\Beq
  \gibbs(\theta,\chi,\nabla\chi,\sigma)
  = -\cz \theta \log\theta
  - \frac 1{\kappa} \, \sigma \cdot \left( \frac 12 \, \sigma + \gamma(\chi)\ee \right)
  + \frac 1 2 \, \thetac F(\chi) + \alpha(\theta) G(\chi) + \frac 12 \, |\nabla\chi|^2.
  \non
\Eeq
Here, $\theta_c$ represents a critical phase transition temperature.
Let us point out that the constitutive relation in \eqref{Iseconda} comes from the relation between $\sigma$ and deformation (here~$\nabla u$)}
\Beq
  \elena{\nabla u = -\frac {\partial \gibbs}{\partial \sigma} = \frac 1\kappa (\sigma + \gamma(\chi)\ee).}
  \non
\Eeq
\pier{Let us refer to \cite{bfg} and \cite{dfg} for any further detail on the model.}

\pier{Now, let us briefly review some \ele{contributions} related to shape memory alloys models. Indeed, the mathematical
analysis of such models produced a rather impressive literature and received a great deal
of attention in recent years. Some of the authors of this article contributed to study the Fr\'emond and 
other models for shape memory alloys (see, e.g., \cite{b03-1,b03-2,bfl06, c95, cfrs06, cls00, cs93}). 
Concerning related phase transition models, we underline that a model for hydrogen storage 
in metal hydrides has been recently investigated in \cite{bcl}, 
by encountering the difficulty \pcgg{due to the term}  $|\dt\chi|^2 $ 
in the energy balance equation, 
but for a simpler \ele{analytical form of the} other equations. \pcgg{In the system 
studied in \cite{bcl}}
\ele{the presence of the quadratic dissipative contribution 
 $|\dt\chi|^2 $  comes from a \generaliz ed form of the principle 
of virtual powers, accounting for micro-forces and 
micro-motions responsible for the phase transition. 
\pcol{Concerning phase 
change models with microscopic motions,} there is a comprehensive literature originating 
from the Fr\'emond theory \cite{Fremond, frem2}.} We 
quote \cite{BFL}, in which the resulting system of phase field type is 
characterized by the occurrence of $ |\dt \chi |^2$ and other nonlinearities 
which were not present in the classical formulation of phase field 
systems (not accounting for microscopic stresses). Several authors have dealt with 
this kind of problems and various situations have been analyzed. However, mainly for 
analytical difficulties, to our knowledge there is no global in time well-posedness 
result for the complete related system  in the 3D (or 2D) case. 
A global existence result is proved in the 1D setting \cite{LSS,LS} 
or for a non-diffusive phase evolution \cite{CLSS}. Other results 
have been obtained for some regularized versions of the problem~\cite{B1}.}

In this paper, \pier{we mainly focus on  the three-dimensional situation.
However, our results cover the lower-dimensional cases 
$\Omega\subset\erre^d$ with $d=1,2$ and with minor changes 
we hope to be able to improve a little the results if \elena{$d=1,2$}.}

The paper is \organiz ed as follows. 
\elena{In the next \pier{section, we list the assumptions on the data of the problem and state the main existence result. 
\gianni{In the same section, we also sketch the strategy of our existence proof, which is based on a double approximation,
namely, first a \regulariz ation in terms of a parameter~$\eps>0$
that also introduces the viscous contribution $-\eps\Delta\dt\chi$ in~\eqref{Iterza},
and then a time \discretiz ation of the \regulariz ed problem. 
In Section~\ref{Approximating}, we keep $\eps$ fixed and solve such approximating problems.}
Section~\ref{Existenceproof} is devoted to the proof of some uniform estimates, 
independent of the \gianni{parameter~$\eps$}, on the approximating solutions; 
then, the passage to the limit procedure as \gianni{$\eps\Seto0$} is carefully detailed.}}


\section{Assumptions and results}
\label{MainResults}
\setcounter{equation}{0}

The aim of this section is to introduce precise assumptions on the functions and the data 
that enter the mathematical problem under investigation, and state our results.
We assume \ele{the domain} $\Omega\subset\erre^3$ to be a bounded open set
with a smooth boundary~$\Gamma$
and  fix a final time $T\in(0,+\infty)$. We ~set 
\Beq
  Q := \Omega \times (0,T).
  \label{defQ}
\Eeq 
We introduce the notation 
\Bsist
  && V := \Huno, \quad
  H := \Ldue , \quad
  \VV := V^3, \quad
  \HH := H^3 
  \label{defspazia}
  \\
  && W := \Hdue
  \aand
  \Wz := \graffe{v\in W:\ \dn v|_\Gamma = 0}
  \label{defspazib}
\Esist
\Accorpa\Defspazi defspazia defspazib
\Accorpa\Notazioni defQ defspazib
and endow the above spaces with their standard norms,
for which we use a notation like~$\normaV\cpto$.
However, we use the same symbol for the norm
in a space and in any power of~it
and simply write $\norma\cpto_p$ for the usual norm 
in~$L^p(\Omega)$ for $p\in[1,+\infty]$.
Moreover, for such values of~$p$,
the conjugate exponent of~$p$ is denoted by~$p'$.
\pier{We identify $H$ to a subspace of $\Vp$ in the usual way,
i.e., in order that $\<v_1,v_2>=\iO v_1v_2$ for every $v_1\in H$ and $v_2\in V$.} 
Finally, as no confusion can arise,
if $Z$ is any Sobolev space,
we use the same symbol $\<\cpto,\cpto>$ for the duality product 
between \pier{the dual space~$Z^*$ and $Z$ itself.} 

For the structure of our system,
we are given constants and functions
in order that the conditions listed below hold true\pier{:}
\Bsist
  \hskip-1cm && \cz \,,\, \kappa \,,\, \thetac > 0 , \quad
  \ee \in \erre^3 \quad \hbox{with $|\ee|=1$}
  \label{hpconst}
  \\
  \hskip-1cm && \gianni{\alpha:[0,+\infty)\to\erre
  \quad \hbox{is a $C^2$ nonnegative function \pier{such that $\alpha(0)=0$}}}
  \label{hpalpha}
  \\
  \hskip-1cm && F = \Funo + \Fdue 
  \quad \hbox{with} \quad
  \Funo : \erre \to [0,+\infty]
  \aand
  \Fdue : \erre \to \erre
  \label{hpF}
  \\
  \hskip-1cm && \hbox{$\Funo$ is convex, proper, \pier{\lsc\ (l.s.c.)}\, 
  and $\,\Fdue$ is of class $C^2$}
  \label{hpFunodue}
  \\
  \hskip-1cm && \hbox{$G:\erre\to\erre$ is a $C^1$ nonnegative function with $G(0)=0$}.
  \label{hpG}
\Esist
Moreover, we assume the following parabolicity, boundedness, \gianni{and growth} conditions
(where the positive constant $C$ can be the same, without loss of generality):
\gianni{%
\Bsist
  \hskip-1cm && \cz - r\alpha''(r) G(s) \geq \deltaz
  \quad \hbox{for some $\deltaz\in(0,\cz)$ and every $r\geq0$, $s\in\erre$}
  \label{parab}
  \\
  \hskip-1cm &&  |\alpha(r)| + |r \alpha'(r)| + r|\alpha''(r)| \leq C
  \quad \hbox{for some $C>0$ and every $r\geq \pier{1}$}
  \label{hptreD}
  \\
  \hskip-1cm && \hbox{$\Fdue''$, $G$, $G'$ are bounded and $G'$ is \Lip\ continuous.}
  \qquad
  \label{hpbdd}
\Esist
Note that \eqref{hpalpha} and \eqref{hptreD} imply that 
$\alpha$, $\alpha'$ and $\alpha''$ are bounded in~$[0,+\infty)$.} \pier{On the other 
hand, we are not going to use the non-negativity  property of $\alpha $ in our proofs.} 
We set for convenience
\Beq
  \elena{\beta := \partial\Funo \,, \quad
  \pi := \Fdue' \, .}
  \label{abbreviaz}
\Eeq
In \eqref{abbreviaz}, the symbol $\partial\Funo$
denotes the subdifferential of~$\Funo$ (\ele{defined in the sense of convex analysis}, see, e.g., \cite[Ex.~2.1.4, p.~21]{Brezis}). 
\pier{Let $\gamma$ be defined by \eqref{defgamma}. 
In view of \eqref{hpG}, $G$~takes a minimum \gianni{at~$0$}, whence}
\elena{$G'(0)=0$. 
\gianni{Consequently, \eqref{hpG} and \eqref{hpbdd}} yield}
\Beq
  \elena{\hbox{$\gamma$ is a $C^1$ function and $\gamma$, $\gamma^{\pier{\,\prime}}$ are bounded and \Lip\ continuous.}}
  \label{hpgamma}
\Eeq
For the forcing terms and the initial data, we require that
\Bsist
  && \hskip-1em
  R := \RO \in \LQ2
  \aand
  R \geq 0 , \quad
  \BO \in \LQ2
  \aand
  \bG \in \H1{\LxG2}
  \qquad\quad
  \label{regRB}
  \\
  && \hskip-1em
  \pier{\thetaz \in H , \quad
  \ln\thetaz \in \Luno ,} \quad
  \uz \in V ,\quad
  \uz' \in H , \quad
  \chiz \in V , \quad
  \Funo(\chiz) \in \Luno
  \qquad\quad
  \label{regzero}
\Esist
and define $\BG\in\H1\Vp$ and $B\in\L2H+\H1\Vp$ 
(where the sum is meaningful in the sense of the embedding $H\subset\Vp$ 
mentioned at the beginning of the section)
as~follows
\Beq
  \< \BG(t) , v > := \iG \bG(t) v
  \quad \hbox{for every $v\in V$ and $t\in[0,T]$}\pier{, \ \quad}
  B := \BO + \BG \,.
  \label{defB}
\Eeq
\Accorpa\HPtutto hpconst defB
\ele{In particular, let us point out that 
the prescribed sign of $R$ in \eqref{regRB} helps in keeping $\theta>0$,
which complies with thermodynamical laws.} 
Let us come to the equations of our systems.
The presence of the quadratic term $|\dt\chi|^2$ in \eqref{Iprima}
forces the function $\theta$ to be rather irregular.
For that reason, it is convenient to introduce a related \elena{auxiliary} function \elena{$w$}
and present the equation for $\theta$ in a different form, namely
\Beq
  \dt w - \alpha(\theta) \dt \elena{G}(\chi) - \Delta\theta
  = R + |\dt\chi|^2, 
  \quad \hbox{where} \quad
  w := \cz\theta + \bigl( \alpha(\theta) - \theta \alpha'(\theta) \bigr) \elena{G}(\chi) 
  \label{primaw}
\Eeq
complemented with the homogeneous Neumann boundary condition for $\theta$
and the proper initial condition for~$w$ \gianni{(derived from~\eqref{Ibczero})}.
More precisely, we can deal with a weak formulation 
of the \pier{resulting initial-boundary value problem:
so,} we state the problem under investigation in the precise form given below.

\elena{%
\Bdef
\label{Defsol}
\gianni{A sextuplet $(w,\theta,u,\sigma,\chi,\xi)$ of functions} 
is a solution to our initial and boundary value problem if 
for some $q\in(1,3/2)$ the conditions
\Bsist
  && w, \theta \in \L q{\Wx{1,q}} 
  \aand
  \theta > 0
  \quad \aeQ
  \label{regtheta}
  \\
  && \dt w\in \L1{(\Wx{1,q'})^*},
  \quad \hbox{with } \ q'=\frac q{q-1}> 3
  \label{regwt}
  \\
  \spacca
  && u \in \H2\Vp \cap \C1H \cap \C0V
  \label{regu}
  \\
  && \sigma \in \C0\HH 
  \label{regsigma}
  \\
  && \chi \in \H1H \cap \L2\Wz
  \label{regchi}
  \\
  && \xi \in \L2H
  \aand
  \xi \in \beta(\chi)
  \quad \aeQ
  \label{regxi}
\Esist
\Accorpa\Regsoluz regtheta regxi
\pier{are fulfilled along with the following equalities}
\Bsist
  && w = \cz\theta + \bigl( \alpha(\theta) - \theta \alpha'(\theta) \bigr) G(\chi) \quad \aeQ
  \label{defw}
  \\[0.2cm]
  && 
  \<\dt w,v> - \<\alpha(\theta) \dt \elena{G}(\chi),v> + \int_\Omega \nabla\theta\cdot\nabla v
  = \<R + |\dt\chi|^2,v>
  \non
  \\
  && \qquad \hbox{\aeT\ and for every $v\in \Wx{1,q'}$}
  \label{primater}
  \\[0.3cm]
  \spacca
  &&  \sigma = \kappa \nabla u - \gamma(\chi)\ee \quad\aeQ
  \label{defsigma}
  \\
  && \< \dt^2 u , v > + \iO \sigma \cdot \nabla v = \< B , v >
  \quad \hbox{\aeT\ and for every $v\in V$}
  \label{seconda}
  \\[0.2cm]
  \spacca
  && \dt\chi - \Delta\chi
  + \thetac \bigl( \xi + \pi(\chi) \bigr) + \alpha(\theta) G'(\chi) - \sigma\cdot\ee \, \gamma^{\pier{\,\prime}}(\chi) = 0
  \quad \aeQ
  \qquad
  \label{terza}
  \\[0.2cm]
  && w(0)=\wz  \,, \quad
  u(0) = \uz \,, \quad
  \dt u(0) = \uzp \,,
  \aand
  \chi(0) = \chiz
  \label{cauchy}
\Esist 
where
\Beq
  \wz := \cz\thetaz + \bigl( \alpha(\thetaz) - \thetaz \alpha'(\thetaz) \bigr) G(\chiz).
  \label{defwz}
\Eeq
\Accorpa\Pbl defw defwz
\Edef
}

\elena{%
\Brem
\label{Osscond}
We observe that the boundary condition for $\chi$ given in \eqref{Ibc} 
is contained in~\eqref{regchi} (see~\eqref{defspazib}).
\pier{On the other hand, the analogous boundary conditions for $\theta$ and $\sigma$ 
are included in the variational equations \eqref{primater} and \eqref{seconda}.} 
\pier{Next, we note that the right hand side of the equality \eqref{primater}  
makes sense although $R+|\dt\chi|^2$ is just in $\L1\Luno$ 
(cf.~\eqref{regchi}): 
indeed, as $q'>3$, it turns out that $\Wx{1,q'}$ is compactly embedded in $\Linfty$ 
and therefore $\Luno\subset(\Wx{1,q'})^*$. 
Moreover, the first initial condition in \eqref{cauchy} is meaningful in $(\Wx{1,q'})^*$ as well.} 
\Erem
}%

Our existence result reads as follows.

\Bthm
\label{Esistenza}
Assume that \HPtutto\ hold.
Then, there \elena{exists} at least a sextuplet \pier{$(w, \theta, u,
\sigma, \chi  ,\xi)$ \elena{which is a solution to our 
problem in the sense of} Definition~\ref{Defsol}. 
In particular, we have that 
\Bsist
  &&  \ln\theta \in L^\infty(0,T; L^1(\Omega)) \cap \L2V
  \label{regaln}
  \\
  && w, \theta \in L^\infty(0,T; L^1(\Omega)) \cap \LQ{4q/3}
  \label{regtheta2}
\Esist
and \eqref{regtheta} hold for every  $q\in[1,5/4)$.}     
\Ethm

Due to the highly nonlinear character of our problem,
in particular, to the presence of the quadratic term $|\dt\chi|^2$
on the \rhs\ of~\eqref{primater},
our study passes through an approximating system,
depending on \elena{a} parameter $\eps\in(0,1)$,
whose solution is much smoother.
Namely, we \pier{perturb} equation~\eqref{terza} by adding a higher order term 
with $\eps$ in front of it; then, \pier{$\eps$ is intended to go to $0$ in the limit.}
On the other side, we \regulariz e \pier{the data and} \gianni{all the nonlinearities 
(in~particular the subdifferential~$\beta$\ele{, \pcgg{as in the forthcoming}~\eqref{defFeps}})}. 
Let us denote by $\gianni\alphaeps$, $\Feps$, $\Geps$, $\geps$, \ele{$\Beps$},  
\elena{the approximating functions, whose regularity will be specified later on.}

This leads to the approximating problem of finding 
a quintuplet \pier{$(\weps,\thetaeps,\ueps,\sigmaeps,\chieps)$ satisfying
\Bsist
  \hskip-1cm && \weps = \cz\thetaeps + \bigl( \alphaeps(\thetaeps) - (\thetaeps+\eps) \alphaeps'(\thetaeps) \bigr) \Geps(\chieps)\quad\aeQ
  \label{defweps}
  \\[0.2cm]
  \spacca
  \hskip-1cm && \iO \dt\weps \, v
  - \iO \alphaeps(\thetaeps) \dt \Geps(\chieps) \, v
  + \iO \nabla\thetaeps \cdot \nabla v
  = \intQ \bigl( R + |\dt\chieps|^2 \bigr) v 
  \non
  \\
  \hskip-1cm && \qquad \hbox{\aeT\ and for every $v\in V$}
  \qquad
  \label{primaepster}
  \\[0.3cm]
  \hskip-1cm && \sigmaeps = \kappa\nabla\ueps - \geps(\chieps)\ee \quad\aeQ
  \label{defsigmaeps}
  \\[0.2cm]
  \spacca
  \hskip-1cm && \iO \dt^2\ueps \, v + \iO \sigmaeps \cdot \nabla v = \iO \Beps v
  \quad \hbox{\aeT\ and for every $v\in V$}
  \label{secondaeps}
  \\[0.2cm]
  \spacca
  \hskip-1cm && \dt\chieps - \Delta\chieps - \eps \dt\Delta\chieps
  + \thetac \Feps'(\chieps) 
  \non
  \\
  \hskip-1cm && \quad {}
  + \alphaeps(\thetaeps) \elena{\Geps'}(\chieps) - \sigmaeps\cdot\ee \, \gpeps(\chieps) = 0
  \quad \aeQ
  \qquad
  \label{terzaeps}
  \\[0.2cm]
  \hskip-1cm && \weps(0)=\wzeps \,, \quad 
  \ueps(0) = \uzeps \,, \quad
  \dt\ueps(0) = \uzepsp \,,
  \aand
  \chieps(0) = \chizeps
  \label{cauchyeps}
\Esist
\Accorpa\Pbleps defweps cauchyeps
where the unknown functions have to fulfill rather strong regularity conditions, namely
\Bsist
  \hskip-.6cm && \weps \in \H1H \cap \L\infty V
   \label{regweps}
   \\ 
  \hskip-.6cm && \thetaeps \in \H1H \cap \L\infty V \cap \L2{\Wz}
  \gianni{\aand
  \thetaeps \geq 0 \quad \aeQ}
  \qquad
  \label{regthetaeps}
  \\
  \hskip-.6cm && \ueps \in \W{2,\infty}H \cap \W{1,\infty}V
  \label{regueps}
  \\
  \hskip-.6cm && \sigmaeps \in \W{1,\infty}\HH
  \label{regsigmaeps}
  \\
  \hskip-.6cm && \chieps \in \H2V \cap \W{1,\infty}\Wz 
  \label{regchieps}
\Esist
\Accorpa\Regsoluzeps regweps regchieps
In} \eqref{terzaeps}, we simply wrote $\Feps'(\chieps)$
since the function $\Feps$ is constructed in the sequel in order that
its derivative $\Feps'$ is an approximation of the whole sum~$\beta+\pi$.

\Brem
\label{Variantiprimaeps}
Equation \eqref{primaepster} \pier{is the approximating version of~\eqref{primater}:} 
\elena{note that (with respect to \eqref{primater}) 
here we are writing integrals \pier{in place of duality pairings}. 
This is possible thanks to the \pier{further} regularity expressed by \pier{\eqref{regweps} and \eqref{regchieps}, 
which ensures that $\dt\chieps\in\LQ\infty$. 
In this setting, \eqref{primaepster}} can} be replaced by \pier{the variational equation corresponding to~\eqref{Iprima}, that~is} 
\Bsist
  && \iO \bigl( \cz - (\thetaeps+\eps) \alphaeps''(\thetaeps) \Geps(\chieps) \bigr) \dt\thetaeps \, v
  \non
  \\
  && \qquad {}
  - \iO (\thetaeps+\eps) \alphaeps'(\thetaeps) \dt\Geps(\chieps) \, v
  \pier{{}+{}} \iO \nabla\thetaeps \cdot \nabla v
  \qquad
  \non
  \\
  && = \iO \bigl( \Reps + |\dt\chieps|^2 \bigr) v
  \qquad \hbox{\aeT\ and for every $v\in V$} . 
  \label{primaeps}
\Esist
It is worth writing \pier{all such} equations for a future convenience.
Moreover, we \pier{point out} that \elena{\eqref{primaepster} \pier{and}  
\eqref{primaeps}} can also be expressed 
in the strong form of boundary value problems, namely \pier{
\Beq
  \dt\weps
  - \alphaeps(\thetaeps) \dt \Geps(\chieps)
  - \Delta\thetaeps
  = \Reps + |\dt\chieps|^2
  \label{strongprimaepsbis}
\Eeq
and}
\Bsist 
  && \bigl( \cz - (\thetaeps+\eps) \alphaeps''(\thetaeps) G(\chieps) \bigr) \dt\thetaeps  \non
  \\
  && \qquad {}
  - (\thetaeps+\eps) \alphaeps'(\thetaeps) \Geps'(\chieps) \dt\chieps 
  - \Delta\thetaeps
 = \Reps + |\dt\chieps|^2
  \label{strongprimaeps}
\Esist
with the boundary condition $\dn\thetaeps=0$ on $\Gamma\times(0,T)$.
\Erem

\gianni{%
For the approximating functions $\alphaeps:[0,+\infty)\to\erre$ we require that
\Bsist
  && \alphaeps(0) = \alphaeps'(0) = \alphaeps''(0) = 0
  \quad \hbox{for every $\eps\in(0,1)$}
  \label{alphaepszero}
  \\
  && |\alphaeps(r)| + |\alphaeps'(r)| + |(r+\eps)\alphaeps'(r)| + |(r+\eps)\alphaeps''(r)|
  \leq C
  \non
  \\
  && \quad \hbox{for some $C>0$ and \pier{all} $r\geq0$ and $\eps\in(0,1)$}
  \label{bddalphaeps}
  \\
  && \alphaeps(r) \to \alpha(r) , \quad 
  (r+\eps) \alphaeps'(r) \to r \alpha'(r)
  \aand
  (r+\eps) \alphaeps''(r) \to r \alpha''(r)
  \non
  \\
  && \quad \hbox{uniformly \pier{in $[0,+\infty) $} as $\eps\Seto0$}
  \label{convalphaeps}
\Esist
and the forthcoming Proposition~\ref{Esistalphaeps} ensures that such an approximation actually exists.
}%
As far as the approximating function~$\Feps$ is concerned,
we define it along with the related function $\Funoeps,\betaeps:\erre\to\erre$ \pier{as follows}.
We~set
\Beq
  \Funoeps(s) := \min_{s'\in\erre} \Bigl( \frac {(s'-s)^2} {2\eps} + \Funo(s') \Bigr) , \quad
  \betaeps := \Funoeps'
  \aand
  \Feps := \Funoeps + \Fdue .
  \label{defFeps}
\Eeq
\Accorpa\HPalphaepsFeps alphaepszero defFeps
Thus, $\betaeps$ is the Yosida \regulariz ation of~$\beta$
(see, e.g., \pier{\cite[p.~28]{Brezis} and \cite[Prop.~2.11, p.~39]{Brezis}} for its basic properties).
Here, we mention that $\Funoeps$ is convex and \gianni{that}
$\betaeps$ is monotone and \Lip\ continuous.
Moreover, for every $s\in\erre$, \pier{
\Bsist
  \hskip-1cm && \hbox{$0\leq\Funoeps(s)\leq\Funo(s)$, \ $\Funoeps(s)\Neto\Funo(s)$ monotonically as $\eps\Seto0$}
  \label{pier2}
  \\
  \hskip-1cm && \hbox{$|\betaeps(s)|\leq|\beta^0(s)|$, \ $\betaeps(s)$ tends to $\beta^0(s)$ monotonically as $\eps\Seto0$}
  \label{pier3}
\Esist
where $\beta^0(s)$ denotes the element of $\beta(s)$ having minimum modulus
provided $\beta(s)\not=\emptyset.$}

\Brem
\label{Yosida}
In order to make the sequel completely rigorous,
we should have taken an approximation of~$\Funo$ 
that is  more regular than the one given by~\eqref{defFeps}.
For instance, in the next section, we use the pointwise values of~\pier{$\Funoeps''$},
while \pier{$\Funoeps'$} is just \Lip\ continuous.
However, we can quote \cite[Sect.~3]{GiRo},
where the reader can find how to smooth the Yosida \regulariz ation
without loosing its basic properties.
So, we behave later on as $\Feps$ were as smooth as needed.
\Erem

For the approximating functions $\Geps$ and $\geps$, 
we still define the latter through the former as we did for~$\gamma$ (cf.~\eqref{defgamma}), namely
\Beq
  \geps(s) := \Geps(s) \sign (\pier{s})
  \quad \hbox{for $s\in\erre$}.
  \label{defgeps}
\Eeq
\pier{Here,} $\Geps:\erre\to\erre$ is defined in order that
\Bsist
  \hskip-1cm && \hbox{$\Geps$ and $\geps$ are of class $C^3$ and $\Geps^{(k)}$ is bounded for $k=0,\dots,3$}
  \label{regGeps}
  \\
  \hskip-1cm && 0 \leq \Geps(s) \leq \sup G
  \aand
  |\Geps'(s)| \leq \sup |G'|
  \quad \hbox{for every $s\in\erre$}
  \label{disugGeps}
  \\
  \hskip-1cm && \Geps^{(k)} \to G^{(k)}
  \quad \hbox{uniformly on every bounded interval for $k=0,1$}
  \label{convGeps}
\Esist
\Accorpa\HPgeps defgeps convGeps
whence an analogous convergence \pier{follows} for~$\geps$.
Such a $\Geps$ can be obtained this way.
We introduce the function $G^\eps$ defined~by
$G^\eps(r)=G(r-\eps)$ if $r\geq\eps$,
$G^\eps(r)=G(r+\eps)$ if $r\leq-\eps$,
and $G^\eps(r)=0$ if $|r|<\eps$\pier{;}
we \gianni{note} that $G^\eps$ is of class $C^1$ by~\eqref{hpG},
and construct $\Geps$ by convolution
with a kernel supported in~$(-\eps,\eps)$.

\gianni{%
\Brem
\label{Parabeps}
We observe once and for all that 
\eqref{parab}, \accorpa{bddalphaeps}{convalphaeps} and \eqref{disugGeps} imply
\Beq
  \deltastar \leq \cz - (r+\eps) \alphaeps''(r) \Geps(s) \leq \cstar
  \quad \hbox{for \pier{all} $r\geq0$ and $s\in\erre$} 
  \label{parabeps}
\Eeq
where, e.g., $\deltastar=\deltaz/2$ and $\cstar$ is some positive constant, 
provided that $\eps$ is small enough.
This has an important consequence, as we show at once.
We have indeed
\Beq
  \cz - (r+\eps) \alphaeps''(r) \Geps(s) 
  = \partial_r \phieps(r,s)
  \label{drphieps}
\Eeq
where 
\Beq
  \phieps(r,s) := \cz r+\bigl(\alphaeps(r)-(r+\eps)\alphaeps'(r)\bigr)G(s)
  \quad \hbox{for $r\geq0$ and $s\in\erre$}. 
  \label{defphieps}
\Eeq
As $\phieps(0,s)=0$ for every $s$ and \eqref{parabeps} 
means $\deltastar\leq\partial_r\phieps\leq\cstar$,
it follows~that
\Beq
  \deltastar r
  \leq \cz r + \bigl( \alphaeps(r) - (r+\eps) \alphaeps'(r) \bigr) \Geps(s)
  \leq \cstar r
  \quad \hbox{for \pier{all} $r\geq0$ and $s\in\erre$}.
  \label{parabepsbis}
\Eeq
Furthermore, we notice that \eqref{parabepsbis} and the positivity of $\thetaeps$
given in \eqref{regthetaeps} \pier{yield}
\Beq
  \weps \geq \deltastar \, \thetaeps \geq 0
  \aand
  \weps \leq \cstar \, \thetaeps 
  \quad \aeQ .
  \label{segnoweps}
\Eeq
For the same reason, the similar inequalities 
\Beq
  \wzeps \geq \deltastar \, \thetazeps \geq 0
  \aand
  \wzeps \leq \cstar \, \thetazeps 
  \quad \aeO 
  \label{segnowzeps}
\Eeq
hold for the initial data we introduce below.
\Erem
}%

\gianni{For the approximating data of} the $\eps$-problem we assume that
\Bsist
  \hskip-1cm && \BOeps, \,  \BGeps \in \H1H , \quad
  \Beps := \BOeps + \BGeps
  \label{regBeps}
  \\
  \hskip-1cm && \pier{\thetazeps \in V,} \quad 
  \gianni{\hbox{$\thetazeps\geq0$ \ \pier{\aeO}}}, \quad
  \uzeps \in \pier W, \quad
  \uzeps' \in \pier V, \quad
  \chizeps \in \Wz 
  \label{regcauchyeps}
  \\
  \hskip-1cm&& \pier{\wzeps
  := \cz\thetazeps 
  + \bigl( \alphaeps(\thetazeps) - \gianni{(\thetazeps + \eps)} \alphaeps'(\thetazeps ) \bigr) \Geps(\chizeps)}
  \label{defwzeps}
  \\
  \hskip-1cm && \pier{\sigmazeps
  := \kappa\nabla\uzeps - \geps(\chizeps)\ee, \quad \div\sigmazeps \in \HH , \quad \sigmazeps\cdot\nu = 0 
  \quad \hbox{ on } \Gamma}
  \label{regsigmazeps}
\Esist
\Accorpa\Regdatieps regBeps regcauchyeps
and that the following \pier{boundedness and} convergence properties 
are satisfied
\Bsist
  \hskip-1cm && \BOeps \to \BO
  \quad \hbox{strongly in $\L2H$}
  \label{convBOeps}
  \\
  \hskip-1cm && \BGeps \to \BG
  \quad \hbox{strongly in $\H1\Vp$}
  \label{convBGeps}
  \\
  \hskip-1cm && \pier{\thetazeps \to \thetaz
  \quad \hbox{strongly in $H$}
  \aand
  \hbox{$\displaystyle - \iO \ln (\thetazeps + \eps) \leq C$}}
  \label{convthetazeps}
  \\
  \hskip-1cm && \uzeps \to \uz
  \quad \hbox{strongly in $V$}
  \aand
  \uzeps' \to \uz'
  \quad \hbox{strongly in $H$}
  \label{convuzeps}
  \\[0.3cm]  
  \hskip-1cm && \chizeps \to \chiz
  \quad \hbox{strongly in $V$} \pier{, 
  \quad
  \hbox{$\eps^{1/2}\normaH{\Delta\chizeps} \leq C$} \aand}
  \non
  \\
  \hskip-1cm && \hskip6cm \pier{\limsup_{\eps \searrow  
  0}\iO\Funoeps(\chizeps) \leq \iO \Funo (\chiz)          }
  \label{convchizeps}
  \\[0.2cm]
  \hskip-1cm && \pier{\wzeps \to \wz
  \quad \hbox{strongly in $H$}} 
  \label{convwzeps}
  \Esist
\Accorpa\HPdatieps regBeps convwzeps
\pier{as $\eps \searrow 0$, where $C$ denotes a constant independent of~$\eps$.}
\pier{Note that $\BOeps,\,\BGeps$ and $\uzeps'$ 
actually exist, just by density. Moreover, it is not difficult to check  
that \eqref{convwzeps} follows from 
\eqref{defwzeps}, \eqref{convthetazeps} and \eqref{convchizeps}, thanks to the 
uniform boundedness and Lipschitz continuity properties expressed in 
\eqref{bddalphaeps}--\eqref{convalphaeps} and 
\eqref{regGeps}--\eqref{convGeps}.}

\gianni{On the \pier{other hand}, it is not obvious that the remainig requirements can actually be fulfilled.
So, we prove the existence of such approximating data
and construct the approximating functions~$\alphaeps$ as well.
We start from the latter.}

\gianni{%
\Bprop
\label{Esistalphaeps}
There exists a family $\graffe{\alphaeps}_{\eps\in(0,1)}$
of $C^2$ functions $\alphaeps:[0,+\infty)\to\erre$ satisfying
conditions \accorpa{alphaepszero}{convalphaeps}.
\Eprop
%
\Bdim
We first introduce suitable continuous functions $\Deps,\,\Zeps:[0,+\infty)\to\erre$ 
that approximate the Dirac mass at the \pier{origin} and the null function, respectively.
For $\eps\in(0,1)$ we~set
\Bsist
  && \Deps(r) := \frac \leps {2\eps^2} \, r
  \quad \hbox{for $r\in[0,\eps]$}, \qquad
  \Deps(r) := \frac \leps {r+\eps}
  \quad \hbox{for $r\in(\eps,\eps^{1/2})$}
  \non
  \\
  && \Deps(r) := \frac {\leps(2\eps^{1/2}-r)} {\eps+\eps^{3/2}}
  \quad \hbox{for $r\in[\eps^{1/2},2\eps^{1/2}]$}, \qquad
  \Deps(r) = 0
  \quad \hbox{for $r>2\eps^{1/2}$} 
  \non
\Esist
where $\leps>0$ is chosen in order that the conditions
\Beq
  \int_0^{+\infty} \hskip-1em\Deps(s) \, ds = 1 , \quad
  0 \leq \Deps(r) \leq \frac {2\leps} {r+\eps}
  \quad \hbox{for every  $r\geq 0$}
  \aand
  \lim_{\eps\seto0} \leps = 0 
  \label{propDeps}
\Eeq
are satisfied, as we show at once.
Indeed, the second \pier{property in} \eqref{propDeps} 
is obvious for $0\leq r\leq\eps^{1/2}$ and for $\eps\geq2\eps^{1/2}$;
on the other hand, we have for $\eps^{1/2}\leq r\leq2\eps^{1/2}$
\Beq
  0 \leq (r+\eps) \Deps(r)
  \leq (2\eps^{1/2}+\eps) \Deps(r)
  \leq (2\eps^{1/2}+\eps) \, \frac {\leps\,\eps^{1/2}} {\eps+\eps^{3/2}}
  \leq 2\leps \,.
  \non
\Eeq
Moreover, $\Deps$~is continuous.
Furthermore, an elementary computation yields
\Beq
  \int_0^{+\infty} \Deps(r) \, dr
  = \leps \Bigl\{ 
    \frac 14 + \ln(1+\eps^{-1/2}) - \ln 2 + \frac 1 {2(1+\eps^{1/2})}
  \Bigr\}
  \non
\Eeq
so that the proper choice of $\leps$ that guarantees the first 
\pier{condition in} \eqref{propDeps} is obvious.
The same choice clearly ensures the third condition as well.
As far as $\Zeps$ is concerned, we need~that
\Bsist
  && \Zeps(0) = 1, \quad
  |\Zeps(r)| \leq 1 
  \quad \hbox{for $r\geq 0$}, \quad
  \Zeps(r) = 0
  \quad \hbox{for $r\geq\eps $}
  \non
  \\
  && \quad \aand
  \int_0^\eps \Zeps(r) \, dr = 0 .
  \label{propZeps}
\Esist 
This is achieved by setting $\Zeps(r):=\zeta(r/\eps)$,
where $\zeta:[0,+\infty)\to\erre$ is a continuous function satisfying the above properties
written with $\eps=1$.
At this point, by setting just for brevity $\auno:=\alpha'(0)$ and $\adue:=\alpha''(0)$,
we define $\alphaeps:[0,+\infty)\to\erre$~by
\Beq  
  \alphaeps(r) 
  := \alpha(r) - \auno \, r
  + \auno \int_0^r \Bigl( \int_0^s \Deps(\tau) \, d\tau \Bigr) \, ds
  - \adue \int_0^r \Bigl( \int_0^s \Zeps(\tau) \, d\tau \Bigr) \, ds \,.
  \label{defalphaeps}
\Eeq 
Then, $\alphaeps$ is a $C^2$ function.
Moreover, \pier{it turns out that $\alphaeps(0)=0$ and, for $r\geq0$,}
\Bsist
  && \alphaeps'(r)
  = \alpha'(r) - \auno
  + \auno \int_0^r \Deps(s) \, ds
  - \adue \int_0^r \Zeps(s) \, ds
  \label{alphaprimo}
  \\
  && \alphaeps''(r) = \alpha''(r) + \auno \, \Deps(r) - \adue \, \Zeps(r)
  \label{alphasecondo}
\Esist
so that $\alphaeps'(0)=\alphaeps''(0)=0$ and 
\Beq
  \alphaeps'(r) = \alpha'(r)
  \quad \hbox{for $r\geq 2\eps^{1/2}$} .
  \label{idalphaprimo}
\Eeq
\pier{The} representation \eqref{alphaprimo} and properties \eqref{propDeps} and~\eqref{propZeps} ensure that
\Beq
  |\alphaeps'(r) - \alpha'(r)|
  \leq |\auno| \, \Bigl| 1 - \int_0^r \Deps(s) \, ds \Bigr|
  + |\adue| \, \Bigl| \int_0^r \Zeps(s) \, ds \Bigr|
  \leq |\auno| + \eps |\adue|
  \quad \hbox{for $r\geq0$}.
  \non
\Eeq
This, \eqref{idalphaprimo} and our assumption \eqref{hptreD} on~$\alpha'$ imply 
the boundedness and convergence conditions 
\accorpa{bddalphaeps}{convalphaeps} involving~$\alphaeps'$.
Furthermore, the inequality
\Beq
  |\alphaeps(r) - \alpha(r)|
  \leq (|\auno| + |\adue|) \, 2\eps^{1/2}
  \quad \hbox{for $r\geq0$} 
  \non
\Eeq
follows as well.
Thus, also uniform boundedness and convergence for $\alphaeps$ are proved.
Finally, we notice that \eqref{idalphaprimo} implies an analoguous identity for~$\alphaeps''$.
Moreover, \eqref{alphasecondo} and \eqref{propDeps} yield
\Bsist
  && |(r+\eps)\alphaeps''(r) - r\alpha''(r)|
  \leq (r+\eps) |\alphaeps''(r) - \alpha''(r)| + \eps |\alpha''(r)|
  \non
  \\
  && \leq (r+\eps) |\auno \Deps(r) - \adue \Zeps(r)| + \eps |\alpha''(r)|
  \leq 2\pier{|\auno|} \leps + |\adue| (2\eps^{1/2}+\eps) + \eps |\alpha''(r)|
  \non
\Esist
for $0\leq r\leq2\eps^{1/2}$,
whence \accorpa{bddalphaeps}{convalphaeps} \pier{are completely shown}, on account of~\eqref{hptreD}.
\Edim
}%

\Bprop
\label{Esistdatieps}
There exist families $\graffe{\thetazeps}$, $\graffe{\uzeps}$ and $\graffe{\chizeps}$
satisfying the conditions contained in \eqref{regcauchyeps}, \pier{\eqref{regsigmazeps}} 
and the convergence and boundedness properties \accorpa{convthetazeps}{convchizeps}.
\Eprop

\Bdim
\pier{Concerning $\graffe{\thetazeps}$, we can take it as the family of solutions 
$\thetazeps\in \Wz$ to the elliptic problem 
\Beq
  \iO \thetazeps \, v
  + \eps \iO \nabla\thetazeps \cdot \nabla v
  = \iO \thetaz \, v
  \quad \hbox{ for all }\, v\in V.
\label{costthetazeps}
\Eeq
Then, it is not difficult to check that $\thetazeps\to\thetaz$
strongly in~$H$ (weak convergence plus convergence of norms) 
as~well as $\thetazeps\geq0$ in~$\Omega$ (positiveness of $\thetaz$ and maximum principle).
\gianni{In order to show the bound contained in \eqref{convthetazeps},}
it suffices to take $v=-1/(\thetazeps+\eps)$ in \eqref{costthetazeps} 
and observe that $r\mapsto -1/(r +\eps)$ is the derivative of the 
convex function $r\mapsto -\ln(r +\eps)$, $r\geq 0$; 
then, we find out~that
\Bsist
  && - \iO \ln(\thetazeps + \eps) + \iO \ln(\thetaz + \eps) 
  \leq - \iO \frac{\thetazeps - \thetaz}{\thetazeps + \eps}
  \non
  \\
  && \leq - \iO \frac{\thetazeps - \thetaz}{\thetazeps + \eps}
  + \iO \frac 1{(\thetazeps + \eps)^2} \, |\nabla\thetazeps|^2
  = 0 
  \non
\Esist
whence 
\Beq
  - \iO \ln(\thetazeps + \eps)
  \leq - \iO \ln(\thetaz + \eps)
  \leq - \iO \ln\thetaz
  \non
\Eeq
and \eqref{convthetazeps} follows from \eqref{regzero}.}
\gianni{As far as the families $\graffe{\uzeps}$ and $\graffe{\chizeps}$ are concerned,
it is more convenient to construct first the latter and then the former. 
Thus, we proceed as follows.} 
\pier{Let $\chizeps\in\Wz$ solve the elliptic equation
\Beq
  \chizeps - \eps \Delta\chizeps  + \eps \betaeps(\chizeps)
  = \chiz
  \quad \hbox{a.e. in $\Omega$}.
  \label{ellittico0}
\Eeq
Hence, let us test \eqref{ellittico0} by $(\chizeps-\chiz)$ and integrate by parts, 
take advantage of the convexity property
\Beq
  \iO \left( \Funoeps(\chizeps)-\Funoeps(\chiz) \right)
  \leq \iO \betaeps(\chizeps)(\chizeps-\chiz) 
\non
\Eeq
and use the elementary Young inequality to obtain 
\Beq
  \normaH{\chizeps-\chiz}^2 
  + \eps \left(\frac 1 2\normaH{\nabla \chizeps}^2 + \iO\Funoeps(\chizeps)\right )
  \leq \eps \left( \frac 12 \, \normaH{\nabla \chiz}^2 + \iO \Funoeps(\chiz) \right).
  \label{ellittico1}
\Eeq
Now, \pcgg{in view of \eqref{regzero} and \eqref{pier2},} \ele{it follows that
\Beq
  \hbox{$\chizeps \to \chiz$ \ strongly in \ $H$}
  \label{ellittico2}
\Eeq }%
and dividing} by $\eps$ in \eqref{ellittico1} \pcgg{leads~to}
\Beq
  \frac12 \, \norma{\nabla \chizeps}^2_H
  + \int_\Omega\Funoeps(\chizeps) \leq
  \frac 12 \, \normaH{\nabla \chiz}^2 + \iO \Funo(\chiz).  \label{pier11}
\Eeq
Next, we can test \eqref{ellittico0} by $-\Delta\chizeps$ and integrate by parts. 
Using \gianni{Young's inequality once more} and exploiting the monotonicity of~$\Funoeps$, 
we are led~to 
\Beq
  \label{ellittico3}
  \frac 12 \, \gianni{\normaH{\nabla\chizeps}^2} + \eps\normaH{\Delta\chizeps}^2
  \leq \frac 12 \, \normaH{\nabla\chiz}^2 \,.
\Eeq
Thus, we can deduce 
\Beq
  \eps^{1/2}  \normaH{\Delta\chizeps} \ \hbox{ bounded independently of $\eps$}
  \non
\Eeq 
and 
\Beq
  \limsup_{\eps\searrow0} \normaH{\nabla\chizeps}^2
  \leq \normaH{\nabla\chiz}^2 \,. 
  \label{ellittico4}
\Eeq
Then, as $\chizeps \to \chiz$ weakly in $V$ by weak compactness and~\eqref{ellittico2}, 
$\chizeps$~strongly converges to $\chiz$ in $V$ thanks to the convergence of norms, 
which is ensured by~\eqref{ellittico4}. At this point, in view of \eqref{pier11}
we easily recover the property 
\Bsist
&&\limsup_{\eps \searrow 0}\iO\Funoeps(\chizeps) \non \\
&&\quad
 \leq 
\limsup_{\eps \searrow 0} \left(\frac12 \, \norma{\nabla \chizeps}^2_H
  + \int_\Omega\Funoeps(\chizeps) \right) + \lim_{\eps \searrow 0} \left( -  \frac12 \, \norma{\nabla \chizeps}^2_H \right) \leq \iO \Funo (\chiz) 
  \non
\Esist
and \eqref{convchizeps} completely follows. 
Finally, let us arrive at the construction of~$\graffe{\uzeps}$. 
Recalling the definition of $\sigmazeps$ in \eqref{regsigmazeps} and invoking the Lax-Milgram lemma, we can take 
$\uzeps\in V$ as the unique solution of the variational equality   
\Beq
  \iO \uzeps \, v + \eps \iO \sigmazeps \cdot \nabla v
  = \iO u_0 \,  v
  \quad \hbox{for every $v\in V$}. 
  \label{defuzeps}
\Eeq
Moreover, it is not difficult to check that \ $-\eps\div\sigmazeps+\uzeps=\uz$ \ or equivalently
\Beq
  - \eps\kappa \, \Delta\uzeps + \uzeps
  = \uz - \eps \, \geps^{\pier{\,\prime}}(\chizeps) \nabla \chizeps \cdot \ee
  \non
\Eeq
in the sense of distributions over $\Omega$, whence by comparison $\Delta\uzeps\in H$ and consequently 
\Beq
  \sigmazeps \cdot \nu = 0
  \quad \hbox{ or } \quad
  \kappa\,\dn \uzeps
  = \geps(\chizeps)_{|_\Gamma} \, \ee\cdot\nu  
  \non
\Eeq
in the sense of traces on~$\Gamma$. 
Then, \eqref{regsigmazeps} holds and, as 
$\geps(\chizeps)_{\gianni{|_\Gamma}}\,\ee\cdot\nu \in H^{1/2}(\Gamma)$, 
standard elliptic regularity properties ensure that $\uzeps\in W$. 
Taking now $v=\kappa(\uzeps-\uz)$ in \eqref{defuzeps} and setting 
$\sigmaz:=\kappa\nabla\uz-\gamma(\chiz)\,\ee$, we easily deduce that
\Beq
  \frac \kappa\eps \iO |\uzeps - \uz|^2
  + \iO \sigmazeps \cdot \bigl( \sigmazeps - \sigmaz + (\geps(\chizeps) - \gamma(\chiz) ) \, \ee \bigr)
  = 0 
  \non
\Eeq
and, with the help of the elementary Young inequality, 
\Beq
  \frac \kappa\eps \iO |\uzeps - \uz|^2
  + \frac12 \iO |\sigmazeps|^2 
  \leq \frac12 \iO |\sigmaz - (\geps(\chizeps) - \gamma(\chiz) ) \, \ee|^2 . 
  \label{pier4}
\Eeq
Then, in view of \eqref{defgeps}--\eqref{convGeps} and \eqref{convchizeps}, we 
infer that  $\uzeps\to\uz$ strongly in $H$ and 
$\sigmazeps$ weakly converges in $\HH$ to some limit which \gianni{must} coincide with $\sigmaz$~as
\Beq
  \iO \sigmazeps\cdot\vv \to \iO \sigmaz\cdot\vv  
  \quad \hbox{for every $\vv\in\gianni{ (H^1_0(\Omega))^3}$}. 
  \non
\Eeq
Moreover, \eqref{pier4} implies that
\Beq
  \limsup_{\eps\seto0} \normaH{\sigmazeps}^2
  \leq \normaH\sigmaz^2
  \non
\Eeq
for $\geps(\chizeps)-\gamma(\chiz)\to 0$ strongly in~$H$. 
At this point, we can conclude the strong convergence of $\sigmazeps$ to~$\sigmaz$, 
and consequently of $\nabla\uzeps$ to~$\nabla\uz$, in~$\HH$, 
and thus complete the proof of \eqref{convuzeps}. 
\Edim

\pier{Now, we resume at the approximating problem in \Pbleps\ 
and observe that even though it looks} much smoother than the original problem \Pbl, 
it is not obvious that it has at least a solution.
The method we use to prove existence relies on a time \discretiz ation.
For that reason, we introduce a notation.

\Bnot
\label{Gennot}
Let $N$ be a positive integer, $\tau$~a positive parameter and $Z$ a vector space.
Then, we define $\dtau:Z^{N+1}\to Z^N$ as follows:
\Bsist
  && \hbox{for $z=(z_0,z_1,\dots ,z_N)\in Z^{N+1}$ and $w=(w_0,\dots,w_{N-1})\in Z^N$} 
  \non
  \\
  && \dtau z = w
  \quad \hbox{means that} \quad
  w_n := \frac{\znp-\zn} \tau
  \quad \hbox{for $n=0,\dots,N-1$}.
  \qquad
  \label{defdelta}
\Esist
Then, for simplicity, we write $\dzn$ instead of $(\dtau z)_n$
and use the same notation $\dtau$ for different choices of the space~$Z$.
We also can iterate such a procedure
and define, e.g., $\dtau^2z=(\dtau^2\zn)_{n=0}^{N-2}\in Z^{N-1}$.
We~have
\Beq
  \dtau^2\zn := \frac {\dtau\znp - \dtau\zn} \tau
  = \frac {\znpp - 2 \znp + \zn} {\tau^2}
  \quad \hbox{for $n=0,\dots,N-2$} .
  \label{defdeltadue}
\Eeq
\Enot

\gianni{%
The time \discretiz ation scheme we \pier{are introducing}
mainly corresponds to replace the time derivative $\dt$
by the different quotient operator~$\dtau$,
the meaning of $\tau$ being $\tau:=T/N$ from now on.
However, we cannot ensure positivity for the discrete temperature.
For that reason
\Beq
  \hbox{we extend $\alphaeps$ to the whole of $\erre$ by setting $\alphaeps(r)=0$ for $r<0$.}
  \label{estalphaeps}
\Eeq
By \eqref{alphaepszero}, such an extension is a $C^2$ function.
At this point, we are ready to go on.
}%
We define the vectors $(\Rn)_{n=0}^N,(\Bn)_{n=0}^N\in H^{N+1}$
by setting
\Beq
  \Rn := \frac 1\tau \int_{n\tau}^{(n+1)\tau} \!\!\Reps(s) \ds
  \aand
  \Bn := \Beps(n\tau)
  \quad \hbox{for $n=0,\dots,N-1$} 
  \label{defRnBn}
\Eeq
and look for vectors 
$(\thetan)_{n=0}^N$, $(\un)_{\gianni{n=-1}}^N$, $(\chin)_{n=0}^N$,
and $(\sigman)_{n=0}^N$
satisfying the conditions listed below
\Bsist
  \hskip-1cm && \hbox{$\thetaz$, $\uz$ and $\chiz$ are the initial data $\thetazeps$, $\uzeps$ and $\chizeps$, respectively}
  \label{innesco}
  \\
  \hskip-1cm && \pier{\umu := \uz - \tau\uzeps'}
  \label{defumu}
  \\
  \hskip-1cm && \thetan \,,\, \chin \in \Wz
  \pier{\aand
  \un \in V
  \quad \hbox{for $n=1,\dots,N$}}
  \label{regveren}
  \\
  \hskip-1cm && \sigman \in \HH
  \qquad \hbox{for $n=0,\dots,N$}
  \label{regaltren}
  \\
  \hskip-1cm && \sigman = \kappa\nabla\un - \pier{\geps}(\chi_n)\ee
  \qquad \hbox{for $n=0,\dots,N$}
  \label{defaltren}
  \\
  \spacca
  \hskip-1cm && \bigl( \cz - (\thetan+\eps) \alphaeps''(\thetan) \Gn \bigr) \dtau\thetan
  - (\thetan+\eps) \apn \Gpn \dtau\chin 
  - \Delta\thetanp
  \qquad
  \non
  \\
  \hskip-1cm && \quad {}= \Rn + |\dtau\chin|^2 
  \qquad \pier{\aeO,} \, \hbox{ for $n=0,\dots,N-1$}
  \label{priman}
  \\
  \hskip-1cm && \iO \dtau^2\unm \, v 
  + \iO \sigmanp {\cdot} \nabla v = \iO \Bn v
  \quad \forall\, v\in V
  \qquad \hbox{for $n=\pier{0},\dots,N-1$}
  \label{secondan}
  \\
  \hskip-1cm && \dtau\chin - \Delta\chinp - \eps \dtau\Delta\chin
  + \thetac \Feps'(\chinp)
  \non
  \\
  \hskip-1cm && \quad {}
  + \an \Gpn
  - \sigman\cdot\ee \, \gpn
  = 0
  \qquad \pier{\aeO,} \,  \hbox{ for $n=0,\dots,N-1$} .
  \label{terzan}
\Esist
\Accorpa\Pbldiscr innesco terzan

It is clear that all the vectors we are dealing with depend on both $\tau$ and~$\eps$,
even though such a \elena{dependence} is not stressed in the notation.
We also remark that the definitions of the $0^{th}$ components of the unknown vectors
might not render the Cauchy data of the original problem.
For instance, $\chiz$~is now given by~\eqref{innesco}
and thus means~$\chizeps$.
Despite of the ambiguous notation, no confusion can arise in the sequel.
Indeed, we deal with the discrete problem and the original problem
in the next two sections, separately.
Namely, in the former we solve problem \Pbldiscr\ and show that
suitable interpolants of the discrete solutions
converge to a solution of the approximating problem as $\tau$ tends to zero
(for~a subsequence).
In the latter, we let $\eps$ tend to zero
and obtain a solution to the original problem~\Pbl.

\bigskip

Now, we list a number of notations and well-known results we owe to throughout the paper.
First of all, we use the H\"older inequality.
Moreover, we account for the continuous embedding
along with the corresponding Sobolev type inequality \elena{(holding in the three-dimensional case)}
\Beq
  \Wx{1,q} \subset \Lx p
  \aand
  \norma v_p \leq C_{p,q}
  \pier{\norma{v}_{W^{1,p}(\Omega)}}
  \quad \hbox{for every $v\in\Wx{1,q}$}
  \label{gensobolev}
\Eeq
provided that
\Bsist
  && 1\leq p\leq q^*:=\frac {3q}{3-q} \,, \quad
  1\leq p < +\infty , \quad
  1\leq p \leq +\infty
  \non
  \\
  && \quad \hbox{according to} \quad
  q<3 , \ q=3 , \ q>3,
  \label{esponsobolev}
\Esist
\elena{respectively.}
In \eqref{gensobolev}, the constant $C_{p,q}$ 
\elena{depends} only on $\Omega$, $p$ and~$q$.
Moreover, $\Linfty$ can be replaced by $\Cx0$ in \eqref{gensobolev} if $q>3$.
The embedding \eqref{gensobolev} is compact for every allowed $p$ if $q\geq3$,
while compactness is true only if $1\leq p<q^*$ if $q<3$.
In particular
\Beq
  V \subset \Lx p
  \aand
  \norma v_p \leq C_{p,2} \normaV v
  \quad \hbox{for $p\in[1,6]$ and every $v\in V$} \elena{,}
  \label{sobolev}
\Eeq
the embedding being compact if $p<6$.
We also take advantage of the compact embedding 
\Beq
  W \subset \Cx0
  \aand
  \norma v_\infty \leq C \normaW v
  \quad \hbox{for every $v\in W$} 
  \label{sobolevW}
\Eeq
where $C$ depends only on~$\Omega$.
Besides, we account for the Poincar\'e type inequality
\Beq
  \normaV v \leq C \bigl( \normaH{\nabla v} + \norma v_1 \bigr)
  \quad \hbox{for every $v\in V$}.
  \label{poincare}
\Eeq
Again, $C$~depends only on~$\Omega$.
Furthermore, we repeatedly make use of
the elementary identity and inequalities
\Bsist
  \hskip-1cm&& a (a-b) = \frac 12 \, a^2 - \frac 12 b^2 + \frac 12 \, (a-b)^2
  \geq \frac 12 \, a^2 - \frac 12 b^2
  \quad \hbox{for every $a,b\in\erre$}
  \label{elementare}
  \\
  \hskip-1cm&& ab\leq \pier{\lambda}a^2 + \frac 1 {4\pier{\lambda}}\,b^2
  \quad \hbox{for every $a,b\in\erre$ and $\pier{\lambda}>0$}
  \label{young}
\Esist
\Accorpa\Elementari elementare young
(and quote \eqref{young} as the elementary Young inequality),
as~well as of the discrete Gronwall lemma in the following form
(see, e.g., \cite[Prop.~2.2.1]{Jerome}):
for nonnegative real numbers $M$ and $a_n,b_n$, $n=0,\dots,N$,
\Bsist
  a_m \leq M + \somma n0{m-1} b_n a_n 
  \quad \hbox{for $m=0,\dots,N$}
  \qquad \hbox{implies}
  \non
  \\
  a_m \leq M \exp \Bigl( \somma n0{m-1} b_n \Bigr)
  \quad \hbox{for $m=0,\dots,N$}.
  \label{dgronwall}
\Esist
Finally, we set 
\Beq
  Q_t := \Omega \times (0,t)
  \quad \hbox{for $t\in[0,T]$}
  \label{defQt}
\Eeq
and, again throughout the paper,
we use a small-case italic $c$ without subscripts $0,1,\dots$ 
(thus, in contrast with, e.g., $\cz$~in \eqref{hpconst}
and $C$ in~\eqref{sobolevW}, where a capital letter is used)
for~different constants, that
may only depend on~$\Omega$, the final time~$T$, 
the shape of the nonlinearities $\alpha$, $F$, $G$, 
and the properties of the data involved in the statements at hand.
Thus, the values of such constants might change from line to line 
and even in the same formula or chain of inequalities. 
A~notation like~$\cdelta$ signals a constant that depends also on the parameter~$\pier{\lambda}$. 
Finally, we write capital letters (with or without subscripts)
for precise values of \pier{constants} we want to refer~to.


\section{The approximating problem}
\label{Approximating}
\setcounter{equation}{0}

\gianni{In this section, we prove an existence result for the approximating problem \Pbleps.
It is understood that assumptions \HPtutto, \HPalphaepsFeps, \HPgeps\ and \HPdatieps\
on the structure, the approximation and the data are in force;
moreover, by accounting for Remark~\ref{Yosida},
we assume $\Feps'$ to be \Lip\ continuous.}
Here is our existence result.

\Bthm
\label{Esistenzaeps}
Problem \Pbleps\ has at least a solution $\pier{(\weps, \thetaeps,\ueps,\sigmaeps,\chieps )}$ satisfying \Regsoluzeps.
\Ethm

The first step consists in proving the existence of a solution to the discrete problem.

\Bprop
\label{Esistenzan}
Assume Notation~\ref{Gennot}.
Then, there exists $\tau_*>0$, depending only on $\thetac$, \elena{$\pi$, and $\Omega$},
such that the discrete problem \Pbldiscr\ has a unique solution 
\Beq
  \pier{(\thetan)_{n=0}^N, \ \ (\un)_{\gianni{n=-1}}^N, \ \ (\chin)_{n=0}^N, \
  \hbox{ and } \ (\sigman)_{n=0}^N  \quad  \hbox{ if } \, \tau<\tau_* .} 
  \non
\Eeq
\Eprop

\Bdim

\elena{We point out that for the existence proof 
it is sufficient to construct an iterative method for} the first three vectors,
since the \pier{fourth one is simply given by~\eqref{defaltren} in terms of 
$ (\un)_{n=0}^N$ and  $(\chin)_{n=0}^N$.}
\elena{First, note that} $\thetaz$, $\uz$ and $\chiz$ are given by~\eqref{innesco} 
and $\pier{\umu}$ is defined by~\eqref{defumu}.  
We compute the other components by the following steps
(also accounting for the proper boundary condition\pier{s} contained in~\eqref{regveren}): 
\pier{inductively for $n=0,\dots,N-1$
\Bsist
  && \llap{$i)$\quad}
  \hbox{solve \eqref{terzan} for $\chinp$}
  \non
  \\  
  && \llap{$ii)$\quad}
  \hbox{solve \eqref{priman} for $\thetanp$}
  \non
  \\  
  && \llap{$iii)$\quad}
  \hbox{solve \eqref{secondan} for $\unp$.}
  \non
\Esist
We have to prove that each of the above steps yields} a well-posed problem.

\smallskip\noindent
\pier{$i)$}
As $\dtau\chin=(\chinp-\chin)/\tau$,
equation \eqref{terzan} has the form
\Beq
  \frac 1\tau \, \chinp - \Delta\chinp - \frac \eps\tau \, \Delta\chinp
  + \thetac \Feps'(\chinp)
  = f_{1,n}
  \non
\Eeq
where $f_{1,n}\in H$ is known \elena{by virtue of the previous step}.
Hence, the solutions to the corresponding homogeneous Neumann boundary value problem 
are the stationary points of the functional $j_n:V\to\erre$ defined~by
\Beq
  j_n(v) := \iO \Bigl( \frac 1 {2\tau} \, v^2 + \thetac \Feps(v) - f_{1,n} v \Bigr)
  + \Bigl( \frac 12 + \frac \eps {2\tau} \Bigr) \iO |\nabla v|^2 .
  \non
\Eeq
We recall notations \eqref{abbreviaz}, \eqref{defFeps},
the regularity assumptions~\eqref{hpbdd},
and that $\Funoeps$ is convex.
Thus, $\Feps''(s)\geq-\sup|\pi'|$ for every $s\in\erre$,
so that $j_n$ is strictly convex and coercive whenever 
$1/\tau>1/\tau_*:=\thetac\sup|\pi'|$.
Therefore, for $\tau<\tau_*$, the functional $j_n$
has a unique stationary point (namely, a~minimum point)
and the problem to be solved has a unique weak solution $\chinp\in V$.
By accounting for elliptic regularity, we then see that $\chinp\in\Wz$.

\smallskip\noindent
\pier{$ii)$}
We set $a_\eps(r,s):=\cz-(r+\eps)\alphaeps''(r)G(s)$ for $r,s\in\erre$
and $a_n:=a_\eps(\thetan,\chin)$,
and observe that equation \eqref{priman} has the form
\Beq
  \frac 1\tau \, a_n\thetanp
  - \Delta\thetanp
  = f_{2,n}
  \label{primaninduz}
\Eeq
where \pier{$f_{2,n}$} is known as well as~$a_n$,
since $\chinp$ has already been computed. \pier{Note that $f_{2,n}\in H$ because, in particular, 
$\dtau\chin\in\Linfty$ by \eqref{sobolevW}.} 
Moreover, $a_n$~is bounded and satisfies $a_n\geq\deltastar$ \aeO\
thanks to~\eqref{parabeps}.
It follows that the corresponding homogeneous Neumann boundary value problem
has a unique weak solution $\thetanp\in V$ with $\Delta\thetanp\in H$ by comparison.
Elliptic regularity then gives $\thetanp\in\Wz$.

\smallskip\noindent
\pier{$iii)$}
As $\dtau^2\unm=(\unp-2\un+\unm)/\tau^2$ for $1\leq n<N$,
equation \eqref{secondan} has the form
\Beq
  \frac 1 {\tau^2} \iO \unp \, v 
  + \kappa \iO \nabla\unp \cdot v
  = \< f_{3,n} , v >
  \quad \hbox{for every $v\in V$}
  \non
\Eeq
where $f_{3,n}\in\Vp$ is known since $\chinp$ has already been computed \pier{in  step $i)$}.
Hence, the existence of a unique solution $\unp\in V$
is ensured by the Lax-Milgram lemma.
\Edim

As announced in the previous section,
the strategy we use to solve the approximating problem \Pbleps\
is the following.
By using the solution to the discrete problem,
we construct some interpolants 
and prove that they converge to the desired solution as $\tau$ tends to zero
by using compactness methods.
Hence, by keeping $\eps$ fixed,
we prove a number of estimates \pier{in terms of} constants
that might depend on~$\eps$
but are independent of the time \discretiz ation parameter~$\tau$,
at least for $\tau$ small enough
(i.e.,~for $\tau>0$ smaller than some $\tau_\eps>0$ that can depend on~$\eps$).
To start, we assume $\tau\leq1$ and $N\geq2$.
Even though $\eps$ is kept fixed in the whole section,
sometimes we distinguish between $\ceps$ and~$c$, 
according to the general rule explained at the end of Section~\ref{MainResults}.
Moreover, in order to unify some cases, we write sums 
that might have an empty set of \pier{indices}.
It is understood that such sums have to be ignored,
or that they vanish by definition.
Thus, we first introduce the interpolants.
Then, we present some useful preliminary material.
Finally, we start with the true proof of Theorem~\ref{Esistenzaeps}.

\Bnot
\label{Interp}
We use Notation~\ref{Gennot} and recall that $\tau:=T/N$ with $N\geq2$
(without stressing the dipendence of $\tau$ on~$N$).
We~set $I_n:=((n-1)\tau,n\tau)$ for $n=1,\dots,N$
and~define the interpolation maps
from $Z^{N+1}$ into \gianni{spaces} $\W{k,\infty}Z$ as follows:
for $z=(z_0,z_1,\dots ,z_N)\in Z^{N+1}$, 
we associate a further coordinate $\zNp$ defined by
\Beq
  \zNp := 2\zN - \zNm
  \label{defzNp}
\Eeq
so that $\dtau\zN=\dtau\zNm$ and $\dtau^2\zNm=0$,
and set
\Bsist
  && \hskip -2em
  \overz ,\, \underz \in \L\infty Z , \quad
  \hz \in \W{1,\infty}Z
  \aand 
  \tz \in \W{2,\infty}Z
  \qquad
  \label{reginterp}
  \\
  && \hskip -2em
  \overz(t) = \zn
  \aand
  \underz(t) = \znm
  \quad \hbox{for a.a.\ $t\in I_n$, \ $n=1,\dots,N$}
  \label{pwconstant}
  \\
  && \hskip -2em
  \hz(0) = z_0
  \aand
  \dt\hz(t) = \dtau\znm
  \quad \hbox{for a.a.\ $t\in I_n$, \ $n=1,\dots,N$}
  \qquad
  \label{pwlinear}
  \\
  && \hskip -2em
  \tz(0) = z_0 , \enskip
  \dt\tz(0) = \dtau z_0
  \enskip \hbox{and} \enskip
  \dt^2\tz(t) = \dtau^2\znm
  \enskip \hbox{for a.a.\ $t\in I_n$, \ $n=1,\dots,N$}.
  \qquad
  \label{pwquadratic}
\Esist
\Enot

\Brem
\label{Reminterp}
The notation we have used recalls its meaning.
Indeed, the maps defined by \accorpa{pwconstant}{pwlinear} provide
the back/forth piece-wise constant and piece-wise linear interpolants of the discrete vectors, respectively,
since we also have $\hz(n\tau)=\zn$ for every~$n$,
and the function \eqref{pwquadratic} is $C^1$ and piece-wise quadratic.
However, the relation between the latter and the original vector only passes through 
the vector $(\dtau\zn)$ of the difference quotiens,
for we have $\dt\tz(n\tau)=\dtau\zn$ for every~$n$,
while no equality entering the values of $\tz$ and $\zn$ with $n>0$ is true.
\Erem

In order to help the reader,
we collect a number of properties involving the interpolants just introduced.

\Bprop
\label{Propinterp}
With Notation~\ref{Interp}, we have
\Beq
  \dt\hz = \underv
  \aand
  \dt\tz = \hv
  \quad \hbox{if} \quad
  \vn = \dtau\zn \pier{,}
  \quad \hbox{for $n=0,\dots,N$} .
  \label{derhztz}
\Eeq
Moreover, if $Z$ is a normed space, we also have
\Bsist
  && \norma\overz_{\L\infty Z}
  = \max_{n=1,\dots,N} \norma\zn_Z \,, \quad
  \norma\underz_{\L\infty Z}
   = \max_{n=0,\dots,N-1} \norma\zn_Z
  \label{ouLinftyZ}
  \\
  && \norma{\dt\hz}_{\L\infty Z}
  = \max_{0\leq n\leq N-1} \normaZ{\dtau\zn}, \quad
  \norma{\dt^2\tz}_{\L\infty Z}
  = \max_{0\leq n\leq N-2} \normaZ{\dtau^2\zn}
  \qquad
  \label{normadtzLinftyZ}
  \\
  \spacca
  && \norma\overz_{\L2Z}^2
  = \tau \somma n1N |\zn|_Z^2 \,, \quad
  \norma\underz_{\L2Z}^2
  = \tau \somma n0{N-1} |\zn|_Z^2 
  \label{ouLdueZ}
  \\
  \spacca
  && \norma{\dt\hz}_{\L2Z}^2
  = \tau \somma n0{N-1} \normaZ{\dtau\zn}^2 \,, \quad
  \norma{\dt^2\tz}_{\L2Z}^2
  = \tau \somma n0{N-2} \normaZ{\dtau^2\zn}^2 
  \label{normadtzLdueZ}
  \\
  && \norma\hz_{\L\infty Z}
  = \max_{n=1,\dots,N} \max\{\norma\znm_Z,\norma\zn_Z\}
  = \max\{\norma{z_0}_Z,\norma\overz_{\L\infty Z}\}
  \qquad\qquad
  \label{normahzLinftyZ}
  \\
  && \norma\hz_{\L2Z}^2
  \leq \tau \somma n1N \bigl( \norma\znm_Z^2 + \norma\zn_Z^2 \bigr)
  \leq \tau \norma{z_0}_Z^2
  + 2 \norma\overz_{\L2Z}^2 
  \label{normahzLdueZ}
  \\
  \spacca
  && \norma{\overz-\hz}_{\L\infty Z}
  = \max_{n=0,\dots,N-1} \norma{\znp-\zn}_Z
  = \tau \, \norma{\dt\hz}_{\L\infty Z}
  \qquad
  \label{interpLinfty}
  \\
  && \norma{\overz-\hz}_{\L2Z}^2
  = \frac \tau 3 \somma n0{N-1} \norma{\znp-\zn}_Z^2
  = \frac {\tau^2} 3 \, \norma{\dt\hz}_{\L2Z}^2
  \label{interpLdue}
\Esist
and the same identities for the difference $\underz-\hz$. 
Finally
\Bsist
  && \norma{\dt\tz-\dt\hz}_{\L2Z}^2
  = \frac {\tau^2} 3 \, \norma{\dt^2\tz}_{\L2Z}^2
  \label{diffdtz}
  \\
  && \norma{\tz-\hz}_{\L\infty Z}^2
  \leq \frac {T\tau^2} 3 \, \norma{\dt^2\tz}_{\L2Z}^2
  \label{difftzhz}
\Esist
\Eprop

\Bdim
Properties \accorpa{derhztz}{normadtzLdueZ}
and \accorpa{interpLinfty}{interpLdue} are \sfw\ to verify by a direct computation.
Relations \accorpa{normahzLinftyZ}{normahzLdueZ} are a consequence of \eqref{ouLinftyZ} and \eqref{ouLdueZ}
since $\hz(t)$ is a convex combination of $\znm$ and $\zn$ for $t\in I_n$.
Finally, \eqref{diffdtz}~follows from the analogue of \eqref{interpLdue} for $\underv-\hv$
(see~\eqref{derhztz}),
and \eqref{difftzhz} is immediately deduced by representing $\tz-\hz$ as the integral of its derivative
and applying \holder's inequality.
\Edim

We also collect a set of inequalities involving difference quotiens
that are useful in the sequel and prepare an easy lemma.
Consider a vector $(\vn)_{n=0}^N$,
where $\vn:\Omega\to\erre$ are measurable functions
and $f:\erre\to\erre$ is, say, continuous.
Then, the definition of \Lip\ continuity, 
the first and second order Taylor expansions
(around either $\vn(x)$ or $\vnp(x)$ \aaO),
and a standard convex inequality yield
\Beq
  |\dtau f(\vn)| \leq \supess |f'| \, |\dtau\vn|
  \quad \aeO\pier{,}
  \quad \hbox{for $n=0,\dots,N-1$}
  \label{disugdtaulip}
\Eeq
if $f$ is \Lip\ continuous,
\Bsist
  \hskip-1cm&& |\dtau f(\vn) - f'(\vnp) \dtau\vn| 
  \leq \pier{\supess |f''|\,}\tau |\dtau\vn|^2
  \quad \aeO\pier{,}
  \quad \hbox{for $n=0,\dots,N-1$}
  \qquad
  \label{disugdtaulipnpu}
  \\
  \hskip-1cm&& |\dtau f(\vn) - f'(\vn) \dtau\vn| 
  \leq \pier{\supess |f''|}\, \tau |\dtau\vn|^2
  \quad \aeO\pier{,}
  \quad \hbox{for $n=0,\dots,N-1$}
  \label{disugdtaulipn}
  \\
  \hskip-1cm&& |\dtau^2 f(\vn) - f'(\vnp)\dtau^2\vn|
  \leq \supess |f''| \, \bigl( |\dtau\vn|^2 + |\dtau\vnp|^2 \bigr)
  \non
  \\
  \hskip-1cm&& \qquad\quad \aeO\pier{,}
  \quad \hbox{for $n=0,\dots,N-2$}
  \label{discrtaylor}
\Esist
if $f$ is $C^1$ and $f'$ is \Lip\ continuous, and
\Beq
  f'(\vnp) \dtau\vn \geq \dtau f(\vn)
  \quad \aeO\pier{,}
  \quad \hbox{for $n=0,\dots,N-1$}
  \label{disugdtauconv}
\Eeq
if $f$ is $C^1$ and convex.
Even though the notation we have used is self-explaining,
we make it precise, e.g., for~$\dtau f(\vn)$:
the vector we apply $\dtau$ to is $(f(\vn))_{n=0}^N\in Z^{N+1}$,
where $Z$ is the vector space of all measurable functions.
Similarly we behave throughout the section with the solution to the discrete problem,
$f$~being one of the nonlinearities involved in our system.

\Blem
\label{Tecnico}
Let $p\geq1$ be an integer and assume that
$(\zn)_{n=0}^p\in H^{p+1}$ and $(\fn)_{n=0}^{p-1}\in H^p$ satisfy
\Beq
  \frac 1\eps \, \znp + \dtau\zn = \fn 
  \quad \aeO
  \quad \hbox{for $n=0,\dots,p-1$}.
  \label{hptecnico}
\Eeq
Then, we have for $m=0,\dots,p-1$
\Bsist
  && \frac \tau\eps \somma n0m \normaH\znp^2
  + \normaH\zmp^2
  \leq \normaH{z_0}^2
  + \eps\tau \somma n0m \normaH\fn^2
  \label{tstecnicoA}
  \\
  && \normaH{\dtau\zm}^2
  \leq 2 \normaH\fm^2 
  + \ceps \, \normaH{z_0}^2
  + \ceps \, \tau \somma n0m \normaH\fn^2 \,.
  \label{tstecnicoB}
\Esist
\Elem

\Bdim
By multiplying \eqref{hptecnico} by~$2\tau\znp$, integrating over~$\Omega$
and owing to the elementary inequalities~\Elementari,
we easily obtain
\Beq
  \frac {2\tau} \eps \iO |\znp|^2
  + \iO |\znp|^2
  - \iO |\zn|^2
  \leq \frac \tau\eps \iO |\znp|^2
  + \eps\tau \iO |\fn|^2 .
  \non
\Eeq
By rearranging and summing over $n=0,\dots,m\leq p-1$,
we trivially deduce~\eqref{tstecnicoA}.
Now, by comparison in~\eqref{hptecnico} written with $n=m$, we~have
\Beq
  \normaH{\dtau\zm}^2
  \leq 2 \normaH\fm^2
  + \frac 4 {\eps^2} \, \normaH\zmp^2
  \non
\Eeq
so that \eqref{tstecnicoB} follows from~\eqref{tstecnicoA}.
\Edim

Next, it is convenient to collect a number of estimates 
involving the forcing terms and the initial data of the discrete problem.
We recall definitions \Regdatieps\ and \accorpa{defRnBn}{defumu}.

\Blem
\label{Stimeinnesco}
\pier{If $(R_n)_{n=0}^N$ and $(B_n)_{n=0}^N$ are specified by \eqref{defRnBn}
and if \eqref{innesco}--\eqref{terzan}} \gianni{are in force, we have that} 
\Bsist
  &
  \tau \displaystyle\somma n0{N-1} \normaH\Rn^2
  \leq \norma\Reps_{\L2H}^2
  \pier{{}\leq{}} c &
  \label{stimeRn}
  \\
  &  \displaystyle
  \max_{0\leq n\leq N} \normaH\Bn
  \leq \norma\Beps_{\L\infty H}
  \pier{{}\leq{}} \ceps\,, \quad
  \tau \somma n0{N-1} \normaH{\dtau\Bn}^2
  \leq \norma{\dt\Beps}_{\L2H}^2
  \pier{{}\leq{}}  \ceps &
  \label{stimeBn}
  \\
  & \pier{\normaV\thetaz \pier{{}\leq{}}  \pier{\ceps} \,, \quad
  \normaV{\dtau \pier{\umu}}
  = \normaV{\uzeps'} \pier{{}\leq{}} \ceps}, \quad
  \pier{ \normaH\sigmaz \leq \ceps \,, \quad
  \normaW\chiz \pier{{}\leq{}} \ceps \,.}
  &  \hskip1cm
  \label{stimezero}
\Esist
Moreover, for $\tau$ small enough, we have 
\Beq
  \normaW{\dtau\chiz} \leq \ceps \,, \enskip
  \pier{\normaH{\dtau^2\umu} \leq \ceps \,, \enskip
  \normaH{\dtau\sigmaz} \leq \ceps \,.}
  \label{stimeinnesco}
\Eeq
\Elem

\Bdim
\pier{Estimates \accorpa{stimeRn}{stimezero}
follow from \Regdatieps, \accorpa{defRnBn}{defumu}
and \eqref{defaltren} with $n=0$, due to 
the boundedness of~$\geps$.}
\pier{In order to prove}~\eqref{stimeinnesco},
we first estimate $\normaV\chiuno$
by testing \eqref{terzan} written \pier{for} $n=0$ by~$\tau\chiuno$.
As \gianni{$\alphaeps(\thetaz)\in H$ (cf.~\eqref{bddalphaeps})},
\pier{owing} to the \Lip\ continuity of the nonlinearities 
and the elementary inequalities \Elementari,
we easily obtain
\Bsist
  && \frac 12 \iO |\chiuno|^2
  + \tau \iO |\nabla\chiuno|^2
  + \frac \eps 2 \iO |\nabla\chiuno|^2
  \non
  \\
  && \leq \frac 12 \iO |\chiz|^2
  + \frac \eps 2 \iO |\nabla\chiz|^2
  + \ceps \, \tau \iO ( 1 + |\sigma_0|^2 + |\chiuno|^2 ) 
  \non
\Esist
and conclude that  
\Beq
  \normaV\chiuno \leq \ceps
  \label{stimachiuno}
\Eeq 
for $\tau$ small enough.
At this point, we rewrite \eqref{terzan} with $n=0$ in the form
\Bsist
  && \frac 1\eps \, (\chiuno-\eps\Delta\chiuno)
  + \dtau (\chiz-\eps\Delta\chiz)
  \non
  \\
  && = \fz
  := \frac 1\eps \, \chiuno
  - \thetac \Feps'(\chiuno)
  - \alphaeps(\thetaz)\Geps'(\chiz)
  + \sigmaz \cdot \ee \geps^{\pier{\,\prime}}(\chiz)
  \label{applTecnico}
\Esist
and notice that $\normaH\fz\leq\ceps$
due to the \Lip\ continuity of~$\Feps'$,
the boundedness of the other nonlinearities, and estimate \eqref{stimachiuno} just obtained.
\pier{Now,} we can apply Lemma~\ref{Tecnico} with $p=1$
and $\zn=\chin-\eps\Delta\chin\in H$ for $n=0,1$\pier{;
thus, we deduce that}
$\normaH{z_1}\leq\ceps$ and $\normaH{\dtau z_0}\leq\ceps$.
Then, the desired estimates follow by elliptic regularity \pier{because $z_0\in H$}.
As a by-product, we have an improvement of~\eqref{stimachiuno},
namely, $\normaW\chiuno\leq\ceps$.
Let us come to the second and third \pier{properties in}~\eqref{stimeinnesco}.
\pier{We take \eqref{secondan} written for $n=0$ and subtract to both sides  
the term $\iO\sigmaz\cdot\nabla v$, then choose $v=\kappa\dtau^2\umu$ 
finding
\Beq
  \kappa \iO |\dtau^2\umu|^2
  + \iO \dtau\sigmaz \cdot \nabla\kappa\tonde{\dtau\uz - \dtau\umu}
  = \kappa \iO \Bz \dtau^2\umu
  - \iO \sigmaz \cdot \nabla(\kappa\dtau^2\umu).
  \label{pier5}
\Eeq
Next, we recall that 
$\nabla(\kappa\dtau\uz)=\dtau\sigmaz+\dtau\geps(\chiz)\ee$ (cf.~\eqref{defaltren})  
and that $\sigmaz$ is nothing but the vector $\sigmazeps$ defined in \eqref{regsigmazeps}, 
so that we can integrate by parts in the last integral of~\eqref{pier5}. 
\gianni{Moreover, \pier{we have} $\dtau\umu=\uzeps'$ by~\eqref{defumu}.
Hence, from \eqref{pier5} it follows that}
\Bsist
  \hskip-2cm && \kappa \iO |\dtau^2\umu|^2
  + \iO |\dtau\sigmaz|^2 
  \non 
  \\
  \hskip-2cm &&
  = - \iO \dtau\sigmaz \cdot (\dtau \geps(\chiz) \ee 
  - \nabla (\kappa\uzeps'))
  + \kappa \iO (\Bz + \div\sigmazeps) \, \dtau^2\umu \,.
  \label{pier6}
\Esist
Now, we want to apply the Young inequality~\eqref{young} 
in the two integrals on the right hand side of~\eqref{pier6}. 
For the treatment of $\dtau \geps(\chiz)$
we invoke \eqref{disugdtaulip} and the boundedness of $\geps^{\pier{\,\prime}}$ along with the 
control $\normaH{\dtau\chiz}\leq\ceps$. 
Then, in view of \eqref{regcauchyeps}, \eqref{stimeBn}, \eqref{regsigmazeps} as well, 
we can proceed and deduce~that  
\Beq
  \frac \kappa 2 \iO |\dtau^2\umu|^2
  + \frac 12 \iO |\dtau\sigmaz|^2 \leq \ceps .
  \non
\Eeq
Consequently, \eqref{stimeinnesco} is completely proved.}
\Edim

At this point, we can start the true proof of Theorem~\ref{Esistenzaeps}.

\step
First a priori estimate

\pier{We choose $v=\kappa\dtau\un$ in \eqref{secondan},
and observe that \eqref{defaltren} yields 
$\nabla v=\dtau\sigman+\dtau\gn\ee$.
Hence, for $0\leq n<N$ we~have}
\Beq
  \frac \kappa\tau \iO (\dtau\un - \dtau\unm) \dtau\un
  + \iO \sigmanp \cdot \dtau\sigman
  = \gianni\kappa \iO \Bn \dtau\un
  - \iO \sigmanp \cdot \ee \, \dtau\gn .
  \non
\Eeq
By accounting for the \holder\ and elementary inequalities \Elementari,
we obtain
\Bsist
  && \frac \kappa {2\tau} \iO \bigl( |\dtau\un|^2 - |\dtau\unm|^2 \bigr)
  + \frac 1 {2\tau} \iO \bigl( |\sigmanp|^2 - |\sigman|^2 \bigr)
  \non
  \\
  && \leq \iO \bigl( 
  \gianni\kappa |\Bn|^2 + \gianni\kappa |\dtau\un|^2 + \cdelta |\sigmanp|^2 
  \bigr)
  + \pier{\lambda}\iO |\dtau\gn|^2 
  \non
  \\
  && \leq \iO \bigl( 
  \gianni\kappa |\Bn|^2 + \gianni\kappa |\dtau\un|^2 + \cdelta |\sigmanp|^2
  \bigr)
  + \pier{\lambda}\sup|\gamma^{\pier{\,\prime}}| \iO |\dtau\chin|^2
  \non
\Esist
for every $\pier{\lambda}>0$,
the last inequality by~\eqref{disugdtaulip}.
Then, we choose $\pier{\lambda}$ such that $\pier{\lambda}\sup|\gamma^{\pier{\,\prime}}|\leq1/4$,
multiply the inequality we get by~$\tau$ and sum over $n=0,\dots,m$,
where $0\leq m<N$.
Hence, by accounting for~\eqref{stimeBn} and with some vanishing empty sum if $m=0$,
we~have
\Bsist
  && \frac \kappa 2 \iO |\dtau\um|^2
  + \frac 12 \iO |\sigmamp|^2
  \non
  \\
  && \leq \frac \kappa 2 \iO |\dtau\umu|^2
  + \frac 12 \iO |\sigmaz|^2
  + \ceps 
  + \tau \iO \bigl( \gianni\kappa |\dtau\um|^2 + c |\sigmamp|^2 \bigr)
  \non
  \\
  && \quad {}
  + \tau \somma n0{m-1} \iO \bigl( \gianni\kappa |\dtau\un|^2 + c |\sigmanp|^2 \bigr)
  + \frac \tau 4 \somma n0m |\dtau\chin|^2
  \non
\Esist
whence also
\Bsist
  && \frac \kappa 3 \iO |\dtau\um|^2
  + \frac 13 \iO |\sigmamp|^2
  \non
  \\
  && \leq \frac \kappa 2 \iO |\dtau\umu|^2
  + \frac 12 \iO |\sigmaz|^2
  + \ceps
  \non
  \\
  && \quad {}
  + \tau \somma n0{m-1} \iO \bigl( \gianni\kappa |\dtau\un|^2 + c |\sigmanp|^2 \bigr)
  + \frac \tau 4 \somma n0m |\dtau\chin|^2
  \non
\Esist
for $\tau$ small enough and $1\leq m<N$.
By Lemma~\ref{Stimeinnesco} \gianni{(see~\eqref{stimezero})},
we can upgrade such an inequality as~follows
\Beq
  \frac \kappa 3 \iO |\dtau\um|^2
  + \frac 13 \iO |\sigmamp|^2
  \leq \ceps
  + \tau \somma n0{m-1} \iO \bigl( |\dtau\un|^2 + c |\sigmanp|^2 \bigr)
  + \frac \tau 4 \somma n0m |\dtau\chin|^2
  \label{primastimaepsA}
\Eeq
for $0\leq m<N$.
Next, \gianni{we add $\chinp$ to both sides of~\eqref{terzan},
multiply the resulting equality by~$\dtau\chin$, integrate over~$\Omega$
by accounting \gianni{for} \eqref{defaltren} and~\eqref{defFeps}, and rearrange.
Owing} to the boundedness of the \gianni{involved nonlinear functions}
and to the elementary Young inequality~\eqref{young},
we \pier{infer that}
\Bsist
  && \iO |\dtau\chin|^2
  + \gianni{\iO \chinp \, \dtau\chin}
  + \iO \nabla\chinp \cdot \nabla\dtau\chin
  + \eps \iO |\dtau\nabla\chin|^2
  + \gianni{\thetac \iO \Funoeps'(\chinp) \dtau\chin}
  \non
  \\
  && = - \iO \an \Gpn \dtau\chin 
  + \iO \sigman \cdot \ee \gpn \dtau\chin
  + \gianni{\iO \bigl( \chinp - \thetac \Fdue'(\chinp) \bigr) \dtau\chin} 
  \non
  \\
  && \leq \ceps \iO (1 + |\sigman| + \gianni{|\chinp|}) |\dtau\chin|
  \non
  \\
  && \leq \pier{\frac 12} \iO |\dtau\chin|^2
  + \ceps \iO |\sigman|^2 
  + \gianni{\ceps \iO |\chinp|^2}
  + \ceps \,.
  \non
\Esist
\gianni{On the other hand,
by applying \eqref{disugdtauconv} to $\Funoeps$, 
we obtain
\Beq
  \iO \Funoeps'(\chinp) \dtau\chin 
  \geq \iO \dtau\Funoeps(\chin) \,.
  \non
\Eeq
}%
Hence, by combining and \pier{applying} the elementary inequality~\eqref{elementare},
we derive~that
\Bsist
  && \frac 12 \iO |\dtau\chin|^2
  + \gianni{\frac 1 {2\tau} \iO |\chinp|^2
  - \frac 1 {2\tau} \iO |\chin|^2}
  + \frac 1 {2\tau} \iO |\nabla\chinp|^2
  - \frac 1 {2\tau} \iO |\nabla\chin|^2
  \non
  \\
  && \quad {}
  + \eps \iO |\dtau\nabla\chin|^2
  + \frac \thetac \tau \iO \gianni{\bigl( \Funoeps(\chinp) - \Funoeps(\chin) \bigr)}
  \leq \ceps \iO |\sigman|^2
  + \gianni{\ceps \iO |\chinp|^2}
  + \ceps 
  \non
\Esist
for $\tau$ small enough.
Now, we multiply by $\tau$ and sum over $n=0,\dots,m$, where $0\leq m<N$,
\pier{obtaining}
\Bsist
  && \frac \tau 2 \somma n0m \iO |\dtau\chin|^2
  + \gianni{\frac 1 2 \iO |\chimp|^2}
  + \frac 1 2 \iO |\nabla\chimp|^2
  \non
  \\
  && \quad {}
  + \eps\tau \somma n0m \iO |\dtau\nabla\chin|^2
  + \thetac \iO \gianni\Funoeps(\chimp)
  \non
  \\
  && \leq {}\gianni{\frac 1 2 \iO |\chiz|^2} 
  + \frac 1 2 \iO |\nabla\chiz|^2
  + \thetac \iO \gianni\Funoeps(\chiz)
  \non
  \\
  && \quad {}
  + \ceps \tau \somma n0m \iO |\sigman|^2
  + \ceps \tau \somma n0m \iO |\chinp|^2
  + \ceps \,.
  \non
\Esist
\gianni{We note that $\normaH\sigmaz\leq\ceps$ and $\normaV\chiz\leq\ceps$ by~\eqref{stimezero}.
Moreover, $\Funoeps(\chimp)$ is nonnegative (see~\eqref{pier2})
and $\chiz=\chizeps$ (by~\eqref{innesco}), whence $\Funoeps(\chiz)$ is~independent of~$\tau$.
Finally, we can absorb the term $\ceps\tau|\chimp|^2$ that appears in the last sum
by the corresponding one on the \lhs\ just by assuming that $\tau$ is small enough.}
So, we improve the above inequality and sum it to~\eqref{primastimaepsA}.
We obtain
\Bsist
  && \frac \kappa 3 \iO |\dtau\um|^2
  + \frac 13 \iO |\sigmamp|^2
  \non
  \\
  && \quad {}
  + \frac \tau 4 \somma n0m \iO |\dtau\chin|^2
  + \gianni{\frac 1 4 \iO |\chimp|^2}
  + \frac 1 2 \iO |\nabla\chimp|^2
  + \eps\tau \somma n0m \iO |\dtau\nabla\chin|^2
  \non
  \\
  && \leq \ceps
  + \gianni{\ceps\,\tau \somma n0{m-1} \iO \bigl( |\dtau\un|^2 + |\sigmanp|^2 + |\chinp|^2 \bigr)}.
  \non
\Esist
\gianni{At this point,
we can apply the discrete Gronwall lemma \eqref{dgronwall} and \pier{deduce}~that}
\Bsist
  && \normaH{\dtau\um}^2
  + \normaH{\sigmamp}^2
  + \normaV\chimp^2
  \non
  \\
  && \quad {}
  + \tau \somma n0m \normaH{\dtau\chin}^2
  + \tau \somma n0m \normaH{\dtau\nabla\chin}^2
  \leq \ceps
  \qquad \hbox{for $0\leq m<N$} .
  \label{primastimadiscr}
\Esist
On the other hand, $\kappa\nabla\ump=\sigmamp+\geps(\chimp)\ee$,
whence also a bound for~$\normaH{\nabla\ump}$ \gianni{follows}.
In terms of the interpolants
(see~Notation~\ref{Interp} and Remark~\ref{Reminterp}),
this means~that
\Beq
  \norma\overu_{\L\infty V}
  + \norma{\dt\hu}_{\L\infty H}
  + \norma\oversigma_{\L\infty H}
  + \norma\overchi_{\L\infty V}
  + \norma{\dt\hchi}_{\L2V}
  \leq \ceps \,.
  \non
\Eeq
We infer that 
$\hu$ is bounded in $\L\infty H$ 
since so is $\dt\hu$ and $\normaH\uz\leq\ceps$.
By also accounting for Proposition~\ref{Propinterp} and~\eqref{stimezero}, 
we can conclude~that
\Bsist
  && \norma\overu_{\L\infty V}
  + \norma\underu_{\L\infty V}
  + \norma\hu_{\W{1,\infty}H}
  \non
  \\
  && \quad {}
  + \norma\oversigma_{\L\infty H}
  + \norma\undersigma_{\L\infty H}
  \non
  \\
  && \quad {}
  + \norma\overchi_{\L\infty V}
  + \norma\underchi_{\L\infty V}
  + \norma\hchi_{\H1V}
  \leq \ceps 
  \label{primastimaeps}
  \\
  \noalign{\medskip}
  && \norma{\overu-\hu}_{\L\infty H}
  + \norma{\underu-\hu}_{\L\infty H}
  \non
  \\
  && \quad {}
  + \norma{\overchi-\hchi}_{\L2V}
  + \norma{\underchi-\hchi}_{\L2V}
  \leq \ceps \, \tau .
  \label{primadiffinterp}
\Esist

\step
Consequence

We set for convenience $\zn:=\chin-\eps\Delta\chin$ for $0\leq n\leq N$.
Then $\zn\in H$ for every $n$ and \eqref{terzan} can be rewritten
in the form \eqref{hptecnico} with $p=N$
and $\normaH\fn\leq\ceps$ for every~$n$,
thanks to~\eqref{primastimadiscr} and the properties of the nonlinearities.
By Lemma~\ref{Tecnico}, we deduce~that
\Beq
  \normaH\zmp + \normaH{\dtau\zm}
  \leq \ceps 
  \quad \hbox{for $m=0,\dots,N-1$} .
  \label{stimaetan}
\Eeq
As $\zm=\chim-\eps\Delta\chim$ and $\chim\in\Wz$,
standard elliptic regularity results yield
\Beq
  \normaW\chimp + \normaW{\dtau\chim}
  \leq \ceps 
  \quad \hbox{for $m=0,\dots,N-1$} 
  \label{stimaWn}
\Eeq
whence also (by the continuous embedding~\eqref{sobolevW})
\Beq
  \norma\chimp_\infty + \norma{\dtau\chim}_\infty
  \leq \ceps 
  \quad \hbox{for $m=0,\dots,N-1$} .
  \label{stimaLinftyn}
\Eeq
In terms of interpolants, the above estimates read
\Beq
  \norma\overchi_{\LQ\infty} + \norma{\gianni\dt\hchi}_{\LQ\infty}
  \leq c \bigl( \norma\overchi_{\L\infty W} + \norma{\gianni\dt\hchi}_{\L\infty W} \bigr)
  \leq \ceps 
  \label{stimachi}
\Eeq
and the similar ones obtained by replacing $\overchi$ by $\underchi$ hold true as well.

\step
Second a priori estimate

We add $\thetanp$ to both sides of~\eqref{priman} for convenience.
Then, we multiply by $\dtau\thetan$ and integrate over~$\Omega$.
Thanks to the parabolicity and elementary inequalities \eqref{parabeps} and~\eqref{elementare},
we~obtain for $0\leq n<N$
\Bsist
  && \deltastar \iO |\dtau\thetan|^2
  + \frac 1 {2\tau} \iO \bigl( |\thetanp|^2 + |\nabla\thetanp|^2 \bigr)
  - \frac 1 {2\tau} \iO \bigl( |\thetan|^2 + |\nabla\thetan|^2 \bigr)
  \non
  \\
  && \leq \iO \bigl(
    \thetanp + \gianni{(\thetan+\eps)} \apn G'(\chin) \dtau\chin + \Rn + |\dtau\chin|^2
  \bigr) \dtau\thetan \,.
  \non
\Esist
Due to the boundedness of all the nonlinear functions involved and to~\eqref{stimaLinftyn},
the \rhs\ of the above inequality is \gianni{bounded~by}
\Bsist
  && \gianni\ceps \iO \bigl( |\thetanp| + |\thetan| + |\Rn| + 1 \bigr) |\dtau\thetan| 
  \non
  \\
  && \leq \frac \deltastar 2 \iO |\dtau\thetan|^2
  + \ceps \iO \bigl( |\thetanp|^2 + |\thetan|^2 + |\Rn|^2 + 1 \bigr).
  \non
\Esist
By combining, multiplying by~$\tau$, summing over $n=0,\dots,m$ with $0\leq m<N$,
and owing to~\eqref{stimeRn}, we deduce~that
\Beq
  \frac \deltastar 2 \, \tau \somma n0m \iO |\dtau\thetan|^2
  + \frac 12 \, \normaV\thetamp^2
  \leq \frac 12 \, \normaV\thetaz^2
  + \ceps \tau \iO |\thetamp|^2
  + \ceps \tau \somma n0m \iO |\thetam|^2
  + \ceps \,.
  \non
\Eeq
Now, \gianni{we absorb the term on the \rhs\ that involves $\thetamp$ by the \lhs\
provided $\tau$ is small enough,
and then apply the discrete Gronwall lemma~\eqref{dgronwall}.
Thus, we conclude~that}
\Beq
  \tau \somma n0m \iO |\dtau\thetan|^2
  + \normaV\thetamp^2
  \leq \ceps
  \quad \hbox{for $m=0,\dots,N-1$} .
  \label{secondastimadiscr}
\Eeq
For the interpolants, this implies that
\Beq
  \norma{\dt\htheta}_{\L2H} + \norma\htheta_{\L\infty V} + \norma\overtheta_{\L\infty V}
  \leq \ceps \,.
  \non
\Eeq
\gianni{By representing $\htheta$ by means of its initial value~$\thetaz$
and its derivative $\dt\htheta$,
observing that $\normaV\thetaz\leq\ceps$ by~\eqref{stimezero}},
and owing to Proposition~\ref{Propinterp},
we conclude~that
\Bsist
  && \norma\htheta_{\H1H\cap\L\infty V}
  + \norma\overtheta_{\L\infty V}
  + \norma\undertheta_{\L\infty V}
  \leq \ceps 
  \label{secondastimaeps}
  \\
  && \norma{\overtheta-\htheta}_{\L2H}
  + \norma{\undertheta-\htheta}_{\L2H}
  \leq \ceps \, \tau \,.
  \label{secondadiffinterp}
\Esist

\step
Consequence

From \eqref{secondastimadiscr} and the previous estimates,
we derive a bound for a higher norm of~$(\thetan)$ by \gianni{comparing terms} in~\eqref{priman}.
\gianni{Indeed, as the terms in front of $\dtau\thetan$ and $\dtau\chin$
are bounded by the properties of the approximating nonlinearities (cf.~\eqref{bddalphaeps} and~\eqref{disugGeps}),  
estimates \eqref{stimaLinftyn}, \eqref{secondastimadiscr}, \eqref{stimeRn}
and elliptic regularity immediately imply~that
\Beq
  \tau \somma n0{N-1} \normaH{\Delta\thetanp}^2 \leq \ceps
  \aand
  \tau \somma n0{N-1} \normaW\thetanp^2 \leq \ceps \,.
  \non
\Eeq
In terms of the interpolant~$\overtheta$, this reads
\Beq
  \norma\overtheta_{\L2W} \leq \ceps \,.
  \label{stimathetabis}
\Eeq
}%

\step
Third a priori estimate

By setting for convenience 
\Beq
  \etan := \dtau\un
  \quad \hbox{for $n=-1,\dots,N-1$}
  \label{defetan}
\Eeq 
and \gianni{recalling \eqref{secondan}},
we see that 
\Beq
  \iO \dtau\etanm v 
  + \iO \sigmanp \cdot \nabla v
  = \iO \Bn v
  \quad \hbox{for every $v\in V$ \ and \ $n=0,\dots,N-1$}.
  \label{secondaetan}
\Eeq
We perform a discrete differentiation with respect to time, i.e.,
we take the difference between \eqref{secondaetan} 
written with $n+1$ in place of $n$ and \eqref{secondaetan} itself
and divide by~$\tau$.
Then, we choose $v=\kappa\dtau\etan$ as a test function
and obtain for $n=0,\dots,N-2$
\Beq
  \frac \kappa \tau \iO \bigl( \dtau\etan - \dtau\etanm \bigr) \dtau\etan
  +  \iO \dtau\sigmanp \cdot \nabla (\kappa\dtau\etan)
  = \kappa \iO \dtau\Bn \, \dtau\etan \,.
  \non
\Eeq
On the other hand, \eqref{defaltren} yields
\Bsist
  && \nabla (\kappa\dtau\etan)
  = \dtau (\kappa\nabla\etan)
  = \frac {\kappa\nabla\dtau\unp - \kappa\nabla\dtau\un} \tau
  \non
  \\
  && = \frac {\dtau\sigmanp - \dtau\sigman} \tau
  + \frac {\dtau \gianni\gpnp) - \dtau \gpn} \tau \, \ee
  = \frac {\dtau\sigmanp - \dtau\sigman} \tau
  + \dtau^2 \gpn \, \ee
  \non
\Esist
so that the above equality becomes
\Bsist
  && \frac \kappa \tau \iO \bigl( \dtau\etan - \dtau\etanm \bigr) \dtau\etan
  + \frac 1\tau \iO \dtau\sigmanp \cdot (\dtau\sigmanp - \dtau\sigman)
  \non
  \\
  && = \kappa \iO \dtau\Bn \, \dtau\etan
  - \iO \dtau\sigmanp \cdot 
    \ee \, \dtau^2 \gpn .
  \non
\Esist
At this point, \gianni{in view of} the elementary inequality~\eqref{elementare},
\gianni{we infer~that}
\Bsist
  && \frac \kappa {2\tau} \iO |\dtau\etan|^2
  - \frac \kappa {2\tau} \iO |\dtau\etanm|^2
  + \frac 1 {2\tau}  \iO |\dtau\sigmanp|^2
  - \frac 1 {2\tau}  \iO |\dtau\sigman|^2
  \non
  \\
  && \leq c \iO |\dtau\Bn|^2
  + c \iO |\dtau\etan|^2
  + \iO |\dtau\sigmanp|^2
  + \iO |\dtau^2 \gpn|^2 .
  \non
\Esist
Now, we multiply by $\tau$ and sum over $n=0,\dots,m$ with $0\leq m\leq N-2$.
By accounting for~\eqref{stimeBn}, we obtain
\Bsist
  && \frac \kappa 2 \iO |\dtau\etam|^2
  + \frac 12  \iO |\dtau\sigmamp|^2
  \non
  \\
  && \leq \frac \kappa 2 \iO |\dtau\eta_{-1}|^2
  + \frac 12  \iO |\dtau\sigmaz|^2
  + \ceps
  \non
  \\
  && \quad {}
  + c \, \tau \somma n0m \iO |\dtau\etan|^2
  + \tau \somma n0m \iO |\dtau\sigmanp|^2
  + \tau \somma n0m \iO |\dtau^2 \gpn|^2 .
  \non
\Esist
Now, we compensate the terms on the \rhs\ that involve $\etam$ and $\sigmamp$
with the \lhs\ by assuming $\tau$ small enough.
Hence, we conclude that for $m=0,\dots,N-2$
(with some vanishing empty sums if~$m=0$)
\Bsist
  && \iO |\dtau\etam|^2
  + \iO |\dtau\sigmamp|^2
  \non
  \\
  && \leq c \iO |\dtau\eta_{-1}|^2
  + c \iO |\dtau\sigmaz|^2
  + \ceps
  \non
  \\
  && \quad {}
  + c \, \tau \somma n0{m-1} \iO |\dtau\etan|^2
  + c \, \tau \somma n0{m-1} \iO |\dtau\sigmanp|^2
  + c \, \tau \somma n0m \iO |\dtau^2 \gpn|^2 .
  \label{quasistimadtsecondan}
\Esist
The first two terms \gianni{on the \rhs} are estimated by~\eqref{stimeinnesco}.
On the other hand, we can apply \eqref{discrtaylor}
with $f=\geps^{\pier{\,\prime}}$ and $\gianni\vn=\chin$ and take advantage of~\eqref{stimaLinftyn} this~way
\Beq
  |\dtau^2 \gpn| 
  \leq \ceps \bigl( |\dtau^2\chin| + |\dtau\chin|^2 + |\dtau\chinp|^2 \bigr)
  \leq \ceps \bigl( |\dtau^2\chin| + 1 \bigr)
  \quad \aeO .
  \non
\Eeq
Hence, inequality~\eqref{quasistimadtsecondan} becomes
\Bsist
  && \iO |\dtau\etam|^2
  + \iO |\dtau\sigmamp|^2
  \non
  \\
  && \leq \ceps
  + c \, \tau \somma n0{m-1} \iO |\dtau\etan|^2
  + c \, \tau \somma n0{m-1} \iO |\dtau\sigmanp|^2
  + C \tau \somma n0m \iO |\dtau\zetan|^2 
  \label{stimadtsecondan}
\Esist
for $m=0,\dots,N-2$,
where we have set for convenience
\Beq
  \zetan := \dtau\chin 
  \label{defzetan}
\Eeq
and marked the constant in front of the last sum by using the capital letter~$C$
for a future reference.
Now, we stop for a while and suitably test the equation obtained
by differentiating \eqref{terzan} in the discrete sense.
Namely, we write \eqref{terzan} with $n+1$ in place of~$n$, take the difference
of the equality we obtain and \eqref{terzan} itself and divide by~$\tau$.
By keeping the notation~\eqref{defzetan}, we multiply by $\dtau\zetan$
and~have for $n=0,\dots,N-2$
\Bsist
  && \iO \bigl( |\dtau\zetan|^2 + \eps |\nabla\dtau\zetan|^2 \bigr)
  + \frac 1\tau \iO \nabla\zetanp \cdot (\nabla\zetanp - \nabla\zetan)
  \non
  \\
  && = - \thetac \iO \dtau\Feps'(\chinp) \, \dtau\zetan
  - \iO \dtau \bigl( \an G'(\chin) \bigr) \, \dtau\zetan
  + \iO \dtau \bigl( \sigman \cdot \ee \gpn \bigr) \, \dtau\zetan \,.
  \non
\Esist
\gianni{Using} the elementary inequality \eqref{elementare} on the \lhs\
and \gianni{exploiting} the boundedness of the nonlinearities, 
inequality \eqref{disugdtaulip} and \gianni{estimates \eqref{primastimadiscr}, \eqref{stimaLinftyn}} on the \rhs,
we easily deduce~that
\Bsist
  && \iO \bigl( |\dtau\zetan|^2 + \eps |\nabla\dtau\zetan|^2 \bigr)
  + \frac 1 {2\tau} \iO |\nabla\zetanp|^2 
  - \frac 1 {2\tau} \iO |\nabla\zetan|^2 
  \non
  \\
  && \leq \ceps \iO \gianni{\bigl(
    |\dtau\chinp| + |\dtau\thetan| + |\dtau\chin| + |\dtau\sigman| + |\sigman| \, |\dtau\chin|
  \bigr) |\dtau\zetan|}
  \non
  \\
  && \leq \frac 12 \iO |\dtau\zetan|^2
  + \ceps \iO \bigl( 1 + |\dtau\thetan|^2 + |\dtau\sigman|^2 \bigr) .
  \non
\Esist
Now, we rearrange and multiply by~$(\gianni 2C+1)\tau$,
where $C$ is the marked constant \gianni{in}~\eqref{stimadtsecondan}.
Then, we sum over $n=0,\dots,m$ with $0\leq m\leq N-2$,
use~\eqref{secondastimadiscr},
and observe that $\zeta_0=\dtau\chiz$ is bounded in $V$ by the first of~\eqref{stimeinnesco}.
We deduce~that
\Beq
  (C+1) \tau \, \gianni{\Graffe{
    \somma n0m \iO \bigl( |\dtau\zetan|^2
    + 2\eps |\nabla\dtau\zetan|^2 \bigr)
    + \iO |\nabla\zetamp|^2 
  }}
  \leq \ceps + \ceps \tau \somma n0m \iO |\dtau\sigman|^2 .
  \non
\Eeq
Now, we add \gianni{this} inequality to \eqref{stimadtsecondan}
and obtain for $0\leq m\leq N-2$
\Bsist
  \hskip-.8cm && \iO |\dtau\etam|^2
  + \iO |\dtau\sigmamp|^2
  + (C+1) \tau \, \gianni{\Graffe{
    \somma n0m \iO \bigl( |\dtau\zetan|^2
    + 2\eps |\nabla\dtau\zetan|^2 \bigr)
    + \iO |\nabla\zetamp|^2 
  }}
  \non
  \\
  \hskip-.8cm && \leq \ceps
  + c \, \tau \somma n0{m-1} \iO |\dtau\etan|^2
  + c \, \tau \somma n0{m-1} \iO |\dtau\sigmanp|^2
  \pier{{}+ C \tau \somma n0m \iO |\dtau\zetan|^2 
  + \ceps \tau \somma n0m \iO |\dtau\sigman|^2 .}
  \non
\Esist
At this point, we rearrange, apply the discrete Gronwall lemma~\eqref{dgronwall},
replace $\etan$ and $\zetan$ by their values (see~\eqref{defetan} and~\eqref{defzetan})
and conclude~that
\Beq
  \normaH{\dtau^2\um}^2
  + \normaH{\dtau\sigmamp}^2
  + \tau \somma n0m \normaV{\dtau^2\chin}^2
  + \normaH{\nabla\dtau\chimp}^2
  \leq \ceps
  \label{terzastimadiscr}
\Eeq
for $0\leq m\leq N-2$.
As $\kappa\dtau\nabla\ump=\dtau\sigmamp+\dtau\geps^{\pier{\,\prime}}(\chimp)$ and \eqref{stimaLinftyn} holds,
we easily deduce that
$\normaH{\dtau\nabla\ump}\leq\ceps$.
By also accounting for~\eqref{stimezero}, 
we conclude~that
\Beq
  \normaV{\dtau\um} \leq \ceps
  \quad \hbox{for $m=0,\dots,N-1$}.
  \label{stimedtu}
\Eeq
In other words, all this reads
\Bsist
  && \norma{\dt^2\tu}_{\L\infty H}
  + \norma{\dt\hsigma}_{\L\infty H}
  + \norma{\dt\hu}_{\L\infty V}
  \non
  \\
  && \quad {}
  + \norma{\dt^2\tchi}_{\L2V}
  + \norma{\nabla\dt\hchi}_{\L\infty H}
  \leq \ceps \,.
  \label{terzastimaeps}
\Esist
Moreover, as both $\tu(0)=\uz$ and $\dt\tu(0)=\dtau\uz$ 
are bounded in~$W$, thus in~$H$, by~\eqref{stimezero},
we derive a bound for $\tu$ itself in~$\W{2,\infty}H$.
A~similar argument yields an estimate for $\tchi$ in~$\H2V$
and Proposition~\ref{Propinterp} and \eqref{primastimaeps} imply
bounds for the different interpolants of the vector~$(\sigman)$.
We collect \gianni{here} some of the consequences we can derive this way.
\gianni{They are \pier{useful} in the sequel:}
\Bsist
  && \norma\tu_{\W{2,\infty}H}
  + \norma\hsigma_{\W{1,\infty}H}
  + \norma\tchi_{\H2V}
  + \norma\hchi_{\W{1,\infty}V}
  \leq \ceps
  \label{daterzastimaeps}
  \\
  \noalign{\medskip}
  && \norma{\tu-\hu}_{\L\infty V}
  + \norma{\oversigma-\undersigma}_{\L\infty H}
  \non
  \\
  && \quad {}
  + \norma{\dt\tchi-\dt\hchi}_{\L2V}
  + \norma{\hchi-\overchi}_{\L\infty V}
  + \norma{\hchi-\underchi}_{\L\infty V}
  \leq \ceps \, \tau.
  \qquad
  \quad
  \label{terzadiffinterp}
\Esist

\step
Conclusion

First of all, we rewrite the equations of the discrete problem as follows
\Bsist
  && \gianni{\bigl( \cz - (\undertheta+\eps) \alphaeps''(\undertheta) \Geps(\underchi) \bigr)} \dt\htheta
  - (\undertheta+\eps) \alphaeps'(\undertheta) \Gpeps(\underchi) \dt\hchi
  - \Delta\overtheta
  \non
  \\
  && \qquad {}= \underR + |\dt\hchi|^2
  \qquad \aeQ
  \label{primainterp}
  \\[0.2cm]
  && \undersigma = \kappa\nabla\underu - \geps^{\pier{\,\prime}}(\underchi) \ee \,, \enskip
  \oversigma = \kappa\nabla\overu - \geps^{\pier{\,\prime}}(\overchi) \ee 
  \quad \aeQ
  \label{defaltreinterp}
  \\[0.2cm]
  && \iO \dt^2\tu \, v
  + \iO \oversigma \cdot \nabla v
  = \iO \underB v
  \quad \hbox{\gianni{\aeT\ and for every $v\in V$}}
  \qquad
  \label{secondainterp}
  \\[0.2cm]
  && \dt\hchi - \Delta\overchi - \eps \Delta\dt\hchi
  + \thetac \Feps'(\overchi)
  \non
  \\
  && \qquad {}
  + \alphaeps(\undertheta) \Geps''(\underchi)
  - \undersigma \cdot \ee \, \geps^{\pier{\,\prime}}(\underchi)
  = 0
  \quad \aeQ 
  \label{terzainterp}
\Esist
and observe that the proper boundary conditions for $\oversigma$ 
is contained in \eqref{secondainterp} in a weak sense,
while the homogeneous Neumann boundary conditions  
for $\overtheta$ and $\overchi$
follow from $\thetan,\chin\in\Wz$ for every~$n$.
Our aim is to let $\tau$ tend to zero in such a problem by compactness methods.
In the sequel, it is understood that the convergence we derive
always holds for a subsequence, even though we never mention this fact.
So, by the a~priori estimates \eqref{primastimaeps},
\eqref{stimachi} as well as its analogue involving~$\underchi$,
\eqref{secondastimaeps}, \eqref{stimathetabis} and \accorpa{terzastimaeps}{daterzastimaeps},
we deduce that all the interpolants we are interested in converge weakly or weakly star
to some limits in the proper topologies.
Moreover, the estimates of the differences given by 
\eqref{primadiffinterp}, \eqref{secondadiffinterp} and \eqref{terzadiffinterp}
imply that some of the weak limits coincide.
Hence, we~have
\Bsist
  && \hskip-3em
  \htheta \to \theta
  \aand 
  \overtheta , \, \undertheta \to \theta
  \qquad \hbox{weakly star in}
  \non
  \\
  && \hskip-3em
  \qquad \hbox{$\H1H\cap\L\infty V\cap\L2{\pier{W}}$ and $\L\infty V$, resp.}
  \label{convtheta}
  \\
  && \hskip-3em
  \tu \to u , \quad
  \hu \to u
  \aand
  \overu , \, \underu \to u
  \qquad \hbox{weakly star in}
  \non
  \\
  && \hskip-3em
  \qquad \hbox{$\W{2,\infty}H$, $\W{1,\infty}V$ and $\L\infty V$, resp.}
  \qquad
  \label{convu}
  \\
  && \hskip-3em
  \hsigma \to \sigma
  \aand
  \oversigma , \, \undersigma \to \sigma
  \qquad \hbox{weakly star in}
  \non
  \\
  && \hskip-3em
  \qquad \hbox{$\W{1,\infty}\HH$ and $\L\infty\HH$, resp.}
  \label{convsigma}
  \\
  && \hskip-3em
  \tchi \to \chi , \quad
  \hchi \to \chi
  \aand
  \overchi , \, \underchi \to \chi
  \qquad \hbox{weakly star in}
  \non
  \\
  && \hskip-3em
  \qquad \hbox{$\H2V$, $\W{1,\infty}V\cap\L\infty W$ and $\L\infty W$, resp.}
  \qquad
  \label{convchi}
\Esist
and $\hchi,\, \overchi, \, \underchi$ converge to $\chi$ weakly star in $\LQ\infty$ as well
because of the continuous embedding $W\subset\Linfty$ (see also~\eqref{stimachi}).
The quadruplet $(\theta,u,\pier{\sigma, \chi})$ 
(we~avoid writing the subscript $\eps$ for simplicity)
is~a candidate to \gianni{satisfy \Regsoluzeps\ and be a solution to problem \Pbleps},
where we forget about $w$ and consider the initial-boundary value problem for \eqref{strongprimaeps} 
in place of the variational equation~\eqref{primaepster}.
We prove that this actually is the case.
The regularity requirements and the Cauchy conditions \eqref{cauchyeps} and \gianni{$\thetaeps(0)=\thetazeps$} are clearly verified.
The homogeneous Neumann boundary conditions for $\theta$ and~$\chi$
are satisfied as well,
since the trace operator $\dn$ is continuous from $\L2{\pier{W}}$ into $\L2{\pier{L^2(\Gamma)}}$.
Thus, it  remains to identify the limits of the nonlinear terms.
To this end, some strong convergence is useful
and we can derive what we need first by accounting for \cite[Sect.~8, Cor.~4]{Simon} 
and the \gianni{compact embeddings~\eqref{sobolev} and~\eqref{sobolevW}},
then by \pier{recalling} the estimates on the differences between the interpolants.
The following is sufficient for the sequel
\Bsist
  && \htheta \to \theta
  \aand
  \dt\tchi \to \dt\chi
  \non
  \\
  && \quad \hbox{strongly in $\C0{\Lx q}$ for $1\leq q<6$}
  \label{strongA}
  \\[0.2cm]
  && \hchi \to \chi
  \quad \hbox{strongly in $\C0{\Cx0}=\CQ$}.
  \label{strongB}
\Esist
Indeed, the second \pier{convergence in} \eqref{strongA}, \gianni{\eqref{stimachi}}
and \eqref{diffdtz} imply that
\Bsist
  && \dt\hchi \to \dt\chi
  \quad \hbox{strongly in $\L2H$}
  \non
  \\
\noalign{\noindent\pier{whence}} 
  && |\dt\hchi|^2 \to |\dt\chi|^2
  \quad \hbox{strongly in $\LQ p$ for every $p<+\infty$}
  \non
\Esist 
due to uniform boundedness in~$\LQ\infty$.
Moreover, we infer that
\Beq
  \undertheta \to \theta
  \aand
  \underchi , \,\overchi \to \chi
  \quad \hbox{strongly in $\L2H$ and in $\L\infty{\Lx6}$, resp.}
  \non
\Eeq
the former by \pier{\eqref{strongA} and~\eqref{secondadiffinterp}}, 
the latter by \eqref{strongB}, \eqref{terzadiffinterp}
and the countinous embedding $V\subset\Lx6$.
\gianni{This and the induced convergence almost everywhere 
imply a proper convergence for the nonlinear terms.
For instance, we have 
$(\undertheta+\eps)\alphaeps''(\undertheta)\to(\theta+\eps)\alphaeps''(\theta)$
and $\Geps(\underchi)\to\Geps(\chi)$
strongly in $\LQ p$ for every $p\in[1,+\infty)$. 
Indeed, $(\undertheta+\eps)\alphaeps''(\undertheta)$ and $\Geps(\underchi)$ are bounded in $\LQ\infty$
since $r\mapsto(r+\eps)\alphaeps''(r)$ and $\Geps$ are bounded functions
(see \eqref{bddalphaeps} and~\eqref{disugGeps}).
}%
Thus, by also accounting for the \holder\ inequality, we immediately see that
(a~better convergence holds true indeed)
\Beq
  (\undertheta+\eps) \alphaeps''(\undertheta) \Geps(\underchi) \dt\htheta
  \to (\theta+\eps) \alphaeps''(\theta) \Geps(\chi) \dt\theta
  \quad \hbox{weakly in $\LQ1$} .
  \non
\Eeq
As the other nonlinear terms and products in system \accorpa{primainterp}{terzainterp}
can be dealt with in a similar and even simpler way,
we conclude that the quadruplet $(\theta,u,\chi,\sigma)$ we have constructed
satisfies \eqref{strongprimaeps}, \eqref{defsigmaeps} and~\eqref{terzaeps},
as well as an integrated form of~\eqref{secondaeps}, namely
\Beq
  \intQ \dt^2 u \, v
  + \intQ \sigma \cdot \nabla v
  = \intQ \Beps v
  \quad \hbox{for every $v\in\L2V$}
  \non
\Eeq
which is equivalent to \eqref{secondaeps} itself.
\gianni{%
It remains to show that the function $\theta=\thetaeps$ we have constructed is nonnegative.
More generally, we can show that the same properties holds for every solution to the approximating system,
provided that the function $\alphaeps$ is extended by~$0$ on the negative half-line (cf.~\eqref{estalphaeps})
in~order that the approximating problem is meaningful without assumptions on the sign of temperature.
We write equation \eqref{primaeps} at the time $s$ with $v=-\thetaeps^-(s)$,
where $(\cdot)^-$ denotes the negative part.
\gianni{Notice} that such a choice of $v$ yields 
\Beq
  \alphaeps''(\thetaeps(s)) v 
  = \alphaeps'(\thetaeps(s)) v 
  = 0
  \aand
  \bigl(\Reps(s) + |\dt\chi(s)|^2 \bigr) v \leq 0
  \non
\Eeq
since $\alphaeps$ vanishes on $(-\infty,0]$ and $\Reps$ is nonnegative.
Hence, after integrating over $(0,t)$ with respect to~$s$, where $t\in(0,T)$ is arbitrary,
and recalling that $\thetazeps\geq0$,
we obtain
\Beq
  \frac \cz 2 \iO |\thetaeps^-(t)|^2 
  + \intQt |\nabla\thetaeps^-|^2
  \leq \frac \cz 2 \iO |\thetazeps^-|^2 
  = 0 .
  \non
\Eeq
Therefore, $\thetaeps^-=0$, whence $\thetaeps\geq0$,
and the proof is complete.\QED
}%


\section{The existence result}
\label{Existenceproof}
\setcounter{equation}{0}

In this section, we prove Theorem~\ref{Esistenza}
by using compactness techniques as before
and monotonicity arguments in addition.
\gianni{We prepare} \elena{a useful energy \pier{equality} 
for equations \accorpa{defsigma}{seconda} and~\eqref{cauchy}.}

\Blem
\label{Energia}
Assume that
\Bsist
  && u \in \W{2,1}\Vp \cap \W{1,\infty}H \cap \L\infty V 
  \label{pocareg}
  \\
  && \sigma \in \L\infty\HH
  \aand
  \chi \in \H1H
  \label{pocaregbis}
\Esist
satisfy \accorpa{defsigma}{seconda} and~\eqref{cauchy}.
Then, \gianni{$u$ and $\sigma$ satisfy \accorpa{regu}{regsigma}}
and the identity
\Bsist
  && \frac \kappa 2 \iO |\dt u(t)|^2
  - \frac \kappa 2 \iO |\uz'|^2
  + \frac 12 \iO |\sigma(t)|^2
  - \frac 12 \iO |\sigma(0)|^2
  + \intQt \sigma \cdot \ee \, \dt\gamma(\chi) 
  \non
  \\
  && = \kappa \intQt \BO \dt u
  + \kappa \< \BG(t) , u(t) > 
  - \kappa \< \BG(0) , \uz >
  - \kappa \iot \< \dt\BG(s) , u(s) > \ds
  \label{energia}
\Esist
holds true for every $t\in[0,T]$.
\Elem

\Bdim
\gianni{On account of \eqref{defsigma}, we write \eqref{seconda} in the form}
\Beq
  \< \dt^2 u , v >
  + \kappa \iO \nabla u \cdot \nabla v
  = \< B + \gammastar , v >
  \quad \hbox{\gianni{\aeT\ and for every $v\in V$}}
  \label{astratta}
\Eeq
where $\gammastar\in\H1\Vp$ is defined by
$\<\gammastar(t),v>:=\iO\gamma(\chi(t))\ee\cdot\nabla v$
\aat\ and every $v\in V$.
\gianni{Let us} read \eqref{astratta} as an abstract second order equation
with a given \rhs.
Then, the Cauchy problem obtained by complementing \eqref{astratta}
with the first two initial conditions \eqref{cauchy}
has a unique solution $u$ satisfying~\eqref{regu}
and a unique \generaliz ed solution in a class of functions
satisfying regularity requirements that are weaker than~\eqref{pocareg}
(see, e.g., \cite[\gianni{Thms.}~3.3 and~4.4]{Baiocchi} or~\cite{DauLio}).
Hence, \eqref{regu}~follows and \eqref{regsigma} is a trivial consequence,
on account of \eqref{defsigma} and~\eqref{pocaregbis}.
In particular, \eqref{energia}~actually is meaningful for every~$t$.
We also observe that \eqref{energia} can be formally obtained 
by choosing $v=\dt u(s)$ in \eqref{astratta} written at the time $s$ 
and then integrating over $(0,t)$ with respect to~$s$.
However, such a choice of $v$ is not allowed due to a lack of regularity.
Therefore, for $\pier{\lambda}>0$,
we introduce the solution $\vdelta$ 
of the time dependent elliptic problem
\Beq
  \vdelta(t) \in \Wz
  \aand
  \vdelta(t) - \pier{\lambda}\Delta\vdelta(t) = u(t)
  \quad \aeO ,
  \quad \hbox{for every $t\in[0,T]$}
  \non
\Eeq
and perform the above formal argument by replacing $u$ by~$\vdelta$
and observing that $\vdelta$ is much smoother that~$u$.
For our purpose, it is sufficient to notice that $\vdelta,\dt\vdelta\in\L\infty\Wz$
and that the following convergence holds true as $\pier{\lambda}\Seto0$ (see, e.g., \cite[Appendix]{CGG})
\Bsist
  && \vdelta(t) \to u(t)
  \quad \hbox{strongly in $V$ for every $t\in[0,T]$}
  \non
  \\
  && \vdelta \to u
  \aand
  \dt\vdelta \to \dt u
  \quad \hbox{strongly in $\L2V$ and $\L2H$, respectively}
  \non
  \\
  && \iot \< \dt^2 u(s) , \dt\vdelta(s) > \ds
  \to \frac 12 \iO |\dt u(t)|^2 - \frac 12 \iO |\dt u(0)|^2
  \quad \hbox{for every $t\in[0,T]$}
  \non
  \\
  && \intQt \nabla u \cdot \nabla\dt\vdelta
  \to \frac 12 \iO |\nabla u(t)|^2 - \frac 12 \iO |\nabla u(0)|^2 
  \quad \hbox{for every $t\in[0,T]$} .
  \non
\Esist
So, we test \eqref{astratta} by $\dt\vdelta$ and take the limit as $\pier{\lambda}\Seto0$.
By the above formulas, the limit of the \lhs\ of the equality we obtain~is
\Beq
  \frac 12 \iO |\dt u(t)|^2 - \frac 12 \iO |\uz'|^2
  + \frac \kappa 2 \iO |\nabla u(t)|^2 - \frac \kappa  2 \iO |\nabla\uz|^2 .
  \label{limlhs}
\Eeq
On the other hand, for $\pier{\lambda}>0$, the \rhs\ of the same equality can be written~as
\Bsist
  && \iot \< B(s) + \gammastar(s) , \dt\vdelta(s) > \ds
  = \intQt \BO \, \dt\vdelta
  + \iot \< \BG(s) + \gammastar(s) , \dt\vdelta(s) > \ds
  \non
  \\
  \spacca
  && = \intQt \BO \, \dt\vdelta
  + \< \BG(t) + \gammastar(t) , \vdelta(t) >
  - \< \BG(0) + \gammastar(0) , \vdelta(0) >
  \non
  \\
  && \quad {}
  - \iot \< \dt\BG(s) + \dt\gammastar(s) , \vdelta(s) > \ds
  \non
  \\
  \spacca
  &&  = \intQt \BO \, \dt\vdelta
  + \< \BG(t) , \vdelta(t) >
  - \< \BG(0) , \vdelta(0) >
  - \iot \< \dt\BG(s) , \vdelta(s) > \ds
  \non
  \\
  && \quad {}
  + \iO \gamma(\chi(t)) \ee \cdot \nabla\vdelta(t) 
  - \iO \gamma(\chiz) \ee \cdot \nabla\vdelta(0) 
  - \intQt \dt\gamma(\chi) \, \ee \cdot \nabla\vdelta \,.
  \non
\Esist 
Hence, its limit as $\pier{\lambda}\Seto0$ has to coincide with~\eqref{limlhs}.
\gianni{Multiplying} by~$\kappa$, we thus obtain
\Bsist
  && \frac \kappa 2 \iO |\dt u(t)|^2 - \frac \kappa 2 \iO |\uz'|^2
  + \frac {\kappa^2} 2 \iO |\nabla u(t)|^2 - \frac {\kappa^2} 2 \iO |\nabla\uz|^2 
  \non
  \\
  \spacca
  &&  = \kappa \intQt \BO \, \dt u
  + \kappa \< \BG(t) , u(t) >
  - \kappa \< \BG(0) , \uz >
  - \kappa \iot \< \dt\BG(s) , u(s) > \ds
  \non
  \\
  && \quad {}
  + \kappa \iO \gamma(\chi(t)) \ee \cdot \nabla u(t) 
  - \kappa \iO \gamma(\chiz) \ee \cdot \nabla\uz 
  - \kappa \intQt \dt\gamma(\chi) \, \ee \cdot \nabla u .
  \label{perenergia}
\Esist
On the other hand, by recalling the definition \eqref{defsigma} of $\sigma$ and that $|\ee|=1$, we have
\Beq
  |\sigma|^2
  = \kappa^2 |\nabla u|^2 + |\gamma(\chi)|^2 - 2 \kappa \gamma(\chi) \ee \cdot \nabla u \non 
\Eeq
and
\Beq
  \sigma \cdot \ee \dt\gamma(\chi)
  = \kappa \dt\gamma(\chi) \ee \cdot \nabla u
  - \gamma(\chi) \dt\gamma(\chi) .
  \non
\Eeq
Therefore, we deduce that
\Bsist
  && \frac 12 \iO |\sigma(t)|^2
  - \frac 12 \iO |\sigma(0)|^2
  + \intQt \sigma \cdot \ee \, \dt\gamma(\chi) 
  \non
  \\
  && = \frac {\kappa^2} 2 \iO |\nabla u(t)|^2
  + \frac 12 \iO |\gamma(\chi(t))|^2
  - \kappa \iO \gamma(\chi(t)) \ee \cdot \nabla u(t)
  \non
  \\
  && \quad {}
  - \frac {\kappa^2} 2 \iO |\nabla\uz|^2
  - \frac 12 \iO |\gamma(\chiz)|^2
  + \kappa \iO \gamma(\chiz) \ee \cdot \nabla\uz
  \non
  \\
  && \quad {}
  + \kappa \intQt \dt\gamma(\chi) \ee \cdot \nabla u
  - \intQt \gamma(\chi) \dt\gamma(\chi)
  \non
  \\
  \spacca
  && = \frac {\kappa^2} 2 \iO |\nabla u(t)|^2 - \frac {\kappa^2} 2 \iO |\nabla\uz|^2 
  \non
  \\
  && \quad {}
  - \kappa \iO \gamma(\chi(t)) \ee \cdot \nabla u(t) 
  + \kappa \iO \gamma(\chiz) \ee \cdot \nabla\uz 
  + \kappa \intQt \dt\gamma(\chi) \, \ee \cdot \nabla u .
  \non
\Esist
By adding this to \eqref{perenergia}, we obtain~\eqref{energia}.
\Edim

\Brem
\label{Energiaeps}
An analogous identity holds for the approximating problem, namely
\Bsist
  && \hskip-1em
  \frac \kappa 2 \iO |\dt\ueps(t)|^2
  - \frac \kappa 2 \iO |\uzeps'|^2
  + \frac 12 \iO |\sigmaeps(t)|^2
  - \frac 12 \iO |\sigmaeps(0)|^2
  + \intQt \sigmaeps \cdot \ee \, \dt\geps(\chieps) 
  \non
  \\
  && \hskip-1em
  = \kappa \intQt \BOeps \dt\ueps
  + \kappa \< \BGeps(t) , \ueps(t) > 
  - \kappa \< \BGeps(0) , \uzeps >
  - \kappa \iot \< \dt\BGeps(s) , \ueps(s) > \ds
  \qquad\quad
  \label{energiaeps}
\Esist
for every $t\in[0,T]$,
and the correponding proof is much simpler.
Indeed, one can test equation \eqref{secondaeps} \pier{directly by~$\dt\ueps$}, 
since the solution and the data are \gianni{smoother}.
\Erem

\gianni{At this point, we recall \eqref{regthetaeps} and Remark~\ref{Parabeps},
in particular that both $\thetaeps$ and $\weps$ are nonnegative,
and start estimating.}

\step
First a priori estimate

Our strategy consists in suitably testing all the equations of the system
and then summing up.
We first take $v=1$ in \eqref{primaepster} and integrate over~$(0,t)$ with $t\in(0,T)$.
\gianni{As $R\in\LQ2$, and \eqref{bddalphaeps}, \eqref{disugGeps} and \eqref{convwzeps} hold, 
using the \holder\ and Young inequalities, we obtain
\Beq
  \iO \weps(t)
  = \iO \wzeps
  + \intQt \alphaeps(\thetaeps) \dt\Geps(\chieps)
  + \intQt \bigr( \Reps + |\dt\chieps|^2 \bigr)
  \leq \frac 32 \intQt |\dt\chieps|^2 
  + c \,.
  \label{testprimabis}
\Eeq
}%
Next, we note that $-1/(\thetaeps+\eps)$ is meaningful and belongs to $\L2V$.
Hence, its values at $s\in(0,T)$ can be chosen as a test function in \eqref{primaeps} written at the time~$s$.
By integrating over~$(0,t)$ with respect to~$s$ and rearranging, we have
\Bsist
  && - \cz \iO \ln(\thetaeps(t)+\eps)
  + \intQt |\nabla\ln(\thetaeps+\eps)|^2
  + \intQt \frac {|\dt\chieps|^2} {\thetaeps+\eps}
  \non
  \\
  && = - \cz \iO \ln(\pier{\thetazeps}  +\eps) 
  \pier{{}-{}} \intQt \bigl\{
    \alphaeps''(\thetaeps) \Geps(\chieps) \dt\thetaeps
    + \alphaeps'(\thetaeps) \dt\Geps(\chieps)
  \bigr\}
  - \intQt \frac \Reps {\thetaeps+\eps} 
  \non
\Esist
\gianni{and observe that the last integral on the \lhs\ is nonnegative.}
On the other hand, \pier{we have that \,
$- \cz \iO \ln(\pier{\thetazeps}  +\eps)\leq c$ \, by \eqref{convthetazeps}} 
and $\Reps\geq0$ by~\accorpa{regRB}{regzero}, 
\pier{and the second integrand} on the \rhs\ can be written as
$\dt\{\alphaeps'(\thetaeps)\Geps(\chieps)\}$.
Moreover, \pier{\eqref{bddalphaeps}}
 and \eqref{disugGeps} hold,
so that both $\alphaeps'$ and $\Geps$ are uniformly bounded.
Hence, the above equality implies
\gianni{%
\Beq
  - \cz \iO \ln(\thetaeps(t)+\eps)
  + \intQt |\nabla\ln(\thetaeps+\eps)|^2
  \leq c .
  \label{testprima}
\Eeq
}%
Now, we write \eqref{secondaeps} at the time~$s$, 
choose $v=\gianni2\kappa\dt\ueps(s)\in V$ as a test function
and observe that $\nabla v=\gianni2\dt\sigmaeps(s)+\gianni2\geps^{\pier{\,\prime}}(\chieps(s))\dt\chieps(s)$
by~\eqref{defsigmaeps}.
Then, we integrate over $(0,t)$ with respect to~$s$
and add the same \gianni{term $2\!\intQt\ueps\dt\ueps$} to both sides for convenience.
\gianni{As the norms $\normaH{\uzeps'}$, $\normaV\uzeps$ and $\normaH{\sigmaeps(0)}$ 
of the initial values
are bounded (see \accorpa{convuzeps}{convchizeps} and~\eqref{regsigmazeps}),
we obtain 
\Bsist
  && \kappa \iO |\dt\ueps(t)|^2
  + \iO |\ueps(t)|^2
  + \iO |\sigmaeps(t)|^2
  + 2 \intQt \sigmaeps \cdot \ee \, \geps^{\pier{\,\prime}}(\chieps) \dt\chieps
  \non
  \\
  && = \kappa \iO |\uzeps'|^2
  + \iO |\uzeps|^2
  + \iO |\sigmaeps(0)|^2
  + \intQt (2\kappa\Beps+2\ueps) \dt\ueps 
  \non
  \\
  && \leq c
  + \intQt |\dt\ueps|^2
  + \intQt |\ueps|^2
  + 2\kappa \intQt \Beps \, \dt\ueps \,.
  \non
\Esist
}%
We recall that $\Beps=\BOeps+\BGeps$ (see~\eqref{regBeps})
and that $\BOeps$ and $\BGeps$ are bounded in $\L2H$ and in $\H1\Vp$, 
respectively (cf.~\accorpa{convBOeps}{convBGeps})\pier{. Hence,
\gianni{for every $\lambda>0$} we have that}
\gianni{%
\Bsist
  && \intQt \Beps \, \dt\ueps 
  = \intQt \BOeps \, \dt\ueps
  \non
  \\
  && \quad {}
  + \iO \BGeps(t) \ueps(t)
  - \iO \BGeps(0) \uzeps
  - \iot \< \dt\BGeps(s) , \ueps(s) > \ds
  \non
  \\
  && \leq \intQt |\dt\ueps|^2
  + \lambda \normaV{\ueps(t)}^2
  + c \iot \normaV{\ueps(s)}^2 \ds
  + \cdelta \,.
  \non
\Esist
Therefore, by combining the last two inequalities, we deduce that 
\Bsist
  && \kappa \iO |\dt\ueps(t)|^2
  + \iO |\ueps(t)|^2
  + \iO |\sigmaeps(t)|^2
  + 2 \intQt \sigmaeps \cdot \ee \, \geps^{\pier{\,\prime}}(\chieps) \dt\chieps
  \non
  \\
  && \leq c \intQt \bigl( |\dt\ueps|^2 + |\ueps|^2 + |\nabla\ueps|^2 \bigr)
  + \lambda \iO |\ueps(t)|^2
  + \lambda \iO |\nabla\ueps(t)|^2
  + \cdelta \,.
  \qquad
  \label{testseconda}
\Esist
}%
Next, 
\gianni{we add $\chieps$ to both sides of~\eqref{terzaeps},
multiply the equality we get by $\gianni2\dt\chieps$, rearrange, and integrate over~$Q_t$.
Using the uniform boundedness of~$\alphaeps$ given by~\eqref{bddalphaeps}, \pier{the Lipschitz 
continuity} of~$\Fdue'$, \eqref{disugGeps} and~\eqref{convchizeps},
we \pier{infer}
\Bsist
  && 2\intQt |\dt\chieps|^2
  + 2\eps \intQt |\nabla\dt\chieps|^2
  +  \iO |\nabla\chieps(t)|^2
  + \iO |\chieps(t)|^2
  + 2\thetac \iO \Funoeps(\chieps(t))
  \non
  \\
  && = \iO \bigl( |\nabla\chizeps|^2 + |\chizeps|^2 + 2\thetac \Funoeps(\chizeps) \bigr)
  \non
  \\
  && \quad {}
  + 2 \intQt \bigl( \chieps - \thetac \Fdue'(\chieps) - \alphaeps(\thetaeps) \Geps'(\chieps) \bigr) \dt\chieps
  + 2 \intQt \sigmaeps \cdot \ee \, \geps^{\pier{\,\prime}}(\chieps) \dt\chieps 
  \non
  \\
  && \leq c + \frac 14 \intQt |\dt\chieps|^2
  + c \intQt |\chieps|^2
  + 2 \intQt \sigmaeps \cdot \ee \, \geps^{\pier{\,\prime}}(\chieps) \dt\chieps \,.
  \label{testterza}
\Esist
}%
\gianni{Finally}, by rearranging \eqref{defsigmaeps} and squaring, applying the elementary Young inequality
and recalling that $\geps$ is uniformly bounded, we have
\Beq
  \frac {\kappa^2} {\gianni 4} \iO |\nabla\ueps(t)|^2
  \leq \frac 1{\gianni 2} \iO |\sigmaeps(t)|^2 + c \,.
  \label{pergradu}
\Eeq
At this point, we sum \accorpa{testprima}{pergradu} to each other.
Then, \gianni{two} terms cancel \ele{and we eventually} obtain
\gianni{%
\Bsist
  && \iO \weps(t)
  - \cz \iO \ln(\thetaeps(t)+\eps)
  + \intQt |\nabla\ln(\thetaeps+\eps)|^2
  \non
  \\
  && \quad {}
  + \kappa \iO |\dt\ueps(t)|^2
  + (1-\lambda) \iO |\ueps(t)|^2
  + \bigl( (\kappa^2/4) - \lambda \bigr) \iO |\nabla\ueps(t)|^2
  \non
  \\
  && \quad {}
  +{} \pier{ \frac 12 \iO |\sigmaeps(t)|^2
  + \frac 14 \intQt |\dt\chieps|^2
  + \eps \intQt |\nabla\dt\chieps|^2
  + \normaV{\chieps(t)}^2}
  \non
  \\
  && \leq c \intQt \bigl(
    |\dt\ueps|^2 + |\ueps|^2 + |\nabla\ueps|^2 + |\chieps|^2
  \bigr)
  + \cdelta \,.
  \label{preprimastima}
\Esist
}%
\gianni{Now, we recall that $\weps\geq\deltastar\thetaeps\geq0$
(see~\pier{\eqref{segnoweps}}), whence 
$\deltastar(\thetaeps+\eps)\leq\weps+\deltastar\eps$,
and observe that $\deltastar r-\cz\ln r\geq(\deltastar/2)(r+|\ln r|)-c$
for some constant $c$ and every $r>0$.
Hence, if we choose $\lambda$ small enough and apply the Gronwall lemma,
we obtain}
\Bsist
  && \norma\weps_{\L\infty\Luno}
  + \norma\thetaeps_{\L\infty\Luno}
  + \norma{\ln(\thetaeps+\eps)}_{\pier{\L\infty\Luno\cap\L2 V}}
  \non
  \\
  && \quad {}
  + \norma\ueps_{\W{1,\infty}H\cap\L\infty V}
  + \norma\sigmaeps_{\L\infty H}
  + \gianni{\norma{\dt\chieps}_{\LQ2}}
  \non
  \\
  && \quad {}
  + \norma\chieps_{\L\infty V}
  + \eps^{1/2} \norma{\nabla\dt\chieps}_{\L\infty H}
  \leq c \,.
  \label{primastima}
\Esist

\step
\gianni{Consequence}

\gianni{A comparison} in \eqref{secondaeps} easily shows that
\Beq
  \norma{\dt^2\ueps}_{\L2\Vp} \leq c \,.
  \label{stimadtseconda}
\Eeq

\step
Second a priori estimate

We rewrite \eqref{secondaeps} as
\Beq
  \pier{ - \Delta\chieps - \eps \dt\Delta\chieps + \thetac \betaeps(\chieps) = \feps}
  \label{newsecondaeps}
\Eeq
\pier{where we have set
\ $ \feps
  := - \dt\chieps
  - \thetac \, \pi(\chieps) 
  - \alphaeps (\thetaeps) \Geps(\chieps)
  + \sigmaeps \cdot \ee \, \geps^{\pier{\,\prime}}(\chieps) .$ \ 
Observe} that $\feps$ is bounded in~$\LQ2$ by~\eqref{primastima}.
By multiplying \eqref{newsecondaeps} by $-\Delta\chieps$
and integrating over~$Q_t$, we thus obtain
\Bsist
  && \intQt |\Delta\chieps|^2
  + \frac \eps 2 \iO |\Delta\chieps(t)|^2
  + \thetac \intQt \betaeps'(\chieps) |\nabla\chieps|^2
  \non
  \\
  && \leq \frac \eps 2 \iO |\Delta\chizeps|^2
  + \frac 12 \intQt |\Delta\chieps|^2
  + c \,.
  \non
\Esist
As $\betaeps'$ is nonnegative 
and $\eps^{1/2}\normaH{\Delta\chizeps}$ is bounded \pier{independently of $\eps$}  by~\eqref{convchizeps}, 
we conclude that
\Beq
  \norma{\Delta\chieps}_{\LQ2} \leq c
  \quad \hbox{whence also} \quad
  \norma\chieps_{\L2{W_0}} \leq c
  \label{secondastima}
\Eeq
by \eqref{primastima} and elliptic regularity.

\step
Consequences

We introduce $\Deltastar:V\to\Vp$ by setting
\Beq
  \< -\Deltastar v , z > := \iO \nabla v \cdot \nabla z
  \quad \hbox{for every $v,z\in V$}.
  \label{defDeltastar}
\Eeq
Then, for $v\in\L2V$, we have
\Beq
  \ioT \< -\dt\Deltastar\chieps(t) , v(t) > \, dt
  = \intQ \nabla\dt\chieps \cdot \nabla v
  \leq \norma{\nabla\dt\chieps}_{\L2H} \norma v_{\L2V}
  \non
\Eeq
and the estimate for $\nabla\dt\chieps$ given by \eqref{primastima} implies that
\Beq
  \eps^{1/2} \norma{\dt\Deltastar\chieps}_{\L2\Vp} \leq c \,.
  \label{stimaDeltaVp}
\Eeq
We deduce an estimate for $\betaeps(\chieps)$ as follows.
We observe that equation \eqref{secondaeps} for $\chieps$ complemented with $\chieps\in\L2\Wz$
can be written as the abstract equation in~$\Vp$
\Beq
  \dt\chieps - \Deltastar\chieps - \eps \dt\Deltastar\chieps 
  + \thetac \, \betaeps(\chieps)
  + \thetac \, \pi(\chieps) 
  + \alphaeps(\thetaeps) \pier{\Gpeps}(\chieps)
  - \sigmaeps \cdot \ee \, \geps^{\pier{\,\prime}}(\chieps)
  = 0
  \non
\Eeq
and that $\Deltastar\chieps=\Delta\chieps$ since $\Deltastar v=\Delta v$ whenever $v\in\Wz$.
Then, \eqref{primastima}, \eqref{secondastima} and \eqref{stimaDeltaVp} yield by comparison 
\Beq
  \norma{\betaeps(\chieps)}_{\L2\Vp} \leq c \,.
  \label{stimabeta}
\Eeq

\step
Third a priori estimate

We adapt the technique of \cite{BoGa} to the present situation
and give the details, for the reader's convenience,
since some modifications of the argument of \cite{BoGa} are spread in the calculation.
\gianni{Here, Remark~\ref{Parabeps} \pier{plays a role}}.
For every nonnegative integer~$k$, we introduce the truncation function
$\Tk:[0,+\infty)\to\erre$ and the set $\Qk$ defined~by 
\Bsist
  && \Tk(r) := \int_0^r \min\{ (s-k)^+,1 \} \ds
  \quad \hbox{for $r\geq0$}
  \non
  \\
  && \Qk := \graffe{(x,t)\in Q:\ k\leq\weps(x,t)<k+1}
  \non
\Esist
and test \eqref{primaepster} written at the time $t$ by~$v=\Tk^{\pier{\, \prime}} (\weps(t))$.
Then, we integrate over~$(0,T)$ with respect to~$t$ and rearrange.
\elena{Once $k$ is fixed, we} easily obtain
\Bsist
  &&\iO \Tk(\weps(T)) 
  + {\intQk} \nabla\thetaeps \cdot \nabla\weps 
  \non
  \\
  && \leq \iO \Tk(\weps(0))
  + \elena{\intQk |\bigl( \alphaeps(\thetaeps) \dt\Geps(\chieps) + \Reps + |\dt\chieps|^2 \bigr)| \, |\Tk^{\pier{\, \prime}} (\weps)|}.
  \non
\Esist
Now, we notice that the first integral on the \lhs\ is nonnegative.
Moreover, the whole \rhs\ is bounded since 
\elena{$|\Tk^{\pier{\, \prime}} |\leq 1$ and} 
\pier{$\iO\Tk(\weps(0))\leq\iO |\weps(0)|\leq c$
by~\eqref{cauchyeps}} and \eqref{convwzeps}.
Hence, we~\pier{infer}
\Beq
  \intQk \nabla\thetaeps \cdot \nabla\weps \leq c \,.
  \label{perBG}
\Eeq
On the other hand, \pier{in view of \eqref{defweps} we have \ $\nabla\weps = a_\eps \, \nabla\thetaeps + b_\eps \, \nabla\chieps$, \ where}
\gianni{%
\Beq
 a_\eps := \cz - (\thetaeps+\eps) \alphaeps''(\thetaeps) \Geps(\chieps)
  \aand
  b_\eps := \bigl( \alphaeps(\thetaeps) - (\thetaeps+\eps) \alphaeps'(\thetaeps) \bigr) \Geps'(\chieps) 
  \non
\Eeq
whence immediately
\Beq
  \nabla\thetaeps  
  = \frac {\nabla\weps - b_\eps \nabla\chieps}{a_\eps} 
  \aand
  \nabla\thetaeps \cdot \nabla\weps 
  = \frac {|\nabla\weps|^2} {a_\eps} 
  - \frac {b_\eps}{a_\eps} \, \nabla\chieps \cdot \nabla\weps \,.
  \qquad
  \label{nablathetaeps}
\Eeq
}%
By accounting for 
\gianni{\eqref{bddalphaeps}, \eqref{disugGeps}, \eqref{parabeps} and estimate~\eqref{primastima}},
we thus obtain
\Beq
  \nabla\thetaeps \cdot \nabla\weps
  \geq \frac {|\nabla\weps|^2} {\gianni\cstar}
  - \frac c \deltastar \, |\nabla\chieps| \, |\nabla\weps| 
  \geq \frac {|\nabla\weps|^2} {2\gianni\cstar}
  - c \, |\nabla\chieps|^2 
  \geq \frac {|\nabla\weps|^2} {2\gianni\cstar}
  - c 
  \non 
\Eeq
and combining with \eqref{perBG}, we conclude~that
\Beq
  \intQk |\nabla\weps|^2 \leq c \,.
  \label{comeBG}
\Eeq
Assume now $q\in[1,5/4)$ and let $|\Qk|$ \gianni{be} the Lebesgue measure of~$\Qk$.
\pier{As $\displaystyle Q= \cup_{k=0}^\infty \Qk$ by \eqref{segnoweps},} we~have
\Beq
  \intQ |\nabla\weps|^q
  = \somma k0\infty \intQk |\nabla\weps|^q
  \leq \somma k0\infty \Bigl( \intQk |\nabla\weps|^2 \Bigr)^{\pier{q/2}} |\Qk|^{(2-q)/2}
  \leq \cq \somma k0\infty |\Qk|^{(2-q)/2} .
  \label{auxBG}
\Eeq
On the other hand, it is clear that for every~$k$
\Beq
  \intQk \weps^{4q/3}
  \geq k^{4q/3} |\Qk|
  \quad \hbox{whence} \quad
  |\Qk| \leq k^{-4q/3} \intQk \weps^{4q/3}
  \non
\Eeq
so that \eqref{auxBG} and the \holder\ inequality for infinite sums yield
\Bsist
  && \intQ |\nabla\weps|^q
  \leq \somma k0\infty k^{-2q(2-q)/3} \Bigl( \intQk \weps^{4q/3} \Bigr)^{(2-q)/2}
  \non
  \\
  && \leq \Bigl( \somma k0\infty k^{-4(2-q)/3} \Bigr)^{q/2}
  \Bigl( \somma k0\infty \intQk \weps^{4q/3} \Bigr)^{(2-q)/2}
  \non
  \\
  && = \cq \Bigl( \intQ \weps^{4q/3} \Bigr)^{(2-q)/2}
  = \cq \norma\weps_{\LQ{4q/3}}^{(2-q)/2}
  \label{quasiBG}
\Esist
where $\cq$ \pier{may denote} the sum of the above numeric series.
Notice that such a series actually converges since $q<5/4$ implies $4(2-q)/3>1$.
Now, we choose $v=\weps(t)$ in the following interpolation and Sobolev-Poincar\'e inequalities
\Beq
  \norma v_{4q/3}
  \leq \norma v_1^{1/4} \norma v_{3q/(3-q)}^{3/4}
  , \quad   \norma v_{3q/(3-q)} \leq \cq \bigl( \norma{\nabla v}_q + \norma v_1 \bigr)
  \quad \hbox{for every $v\in\Wx{1,q}$}
  \non
\Eeq
and integrate over~$(0,T)$. 
\pier{Recalling} the estimate for $\weps$ given by~\eqref{primastima} and combining 
with \eqref{quasiBG}, we obtain 
\Bsist
  && \norma\weps_{\LQ{4q/3}}^{4q/3}
  = \ioT \norma{\weps(t)}_{4q/3}^{4q/3} \, dt
  \leq \cq \ioT \norma{\weps(t)}_{3q/(3-q)}^q \, dt
  \non
  \\
  && \leq \cq \pier{\ioT} \bigl( \norma{\nabla\weps(t)}_q^q + 1 \bigr) \, dt
  = \cq \norma{\nabla\weps}_{\LQ q}^q + \cq 
  \leq \cq \norma\weps_{\LQ{4q/3}}^{(2-q)/2} + \cq \,.
  \non
\Esist
As $(2-q)/2<4q/3$, we infer that $\norma\weps_{\LQ{4q/3}}$ is bounded.
By using~\eqref{quasiBG} again, we conclude that
\Beq
  \norma\weps_{\LQ{4q/3}\cap\L q{\Wx{1,q}}}
  \leq \cq 
  \quad \hbox{for every $q\in[1,5/4)$}.
  \label{stimaBGw}
\Eeq
Due to \eqref{nablathetaeps}, \eqref{parabeps} and estimate~\eqref{primastima},
we derive that
\Beq
  \norma\thetaeps_{\LQ{4q/3}\cap\L q{\Wx{1,q}}}
  \leq \cq
  \quad \hbox{for every $q\in[1,5/4)$}.
  \label{stimaBGtheta}
\Eeq

\step
Consequence

We write \pier{\eqref{primaepster}} \aat\
and take any $v\in\Wx{1,q'}$ as a test function,
by noting that $\Wx{1,q'}\subset\Linfty$ since $q<5/4$ implies $q'>3$.
\gianni{As $\alphaeps$ and $\Geps'$ are uniformly bounded}, we obtain \aat
\Bsist
  && \iO \dt\weps(t) \, v
  = \iO \alphaeps(\thetaeps(t)) \dt \Geps(\chieps(t)) \, v
  - \iO \nabla\thetaeps(t) \cdot \nabla v
  + \iO \bigl( \Reps(t) + |\dt\chieps(t)|^2 \bigr) v
  \non
  \\
  && \leq c \normaH{\dt\chieps(t)} \norma v_\infty
  + \norma{\nabla\thetaeps(t)}_q \norma{\pier{\nabla v}}_{q'}
  + \normaH{R(t)} \norma v_\infty
  + \normaH{\dt\chieps(t)}^2 \norma v_\infty
  \non 
  \\
  \noalign{\smallskip}
  && \leq c \bigl( 1 + \normaH{\dt\chieps(t)}^2 + 
  \pier{\norma{\thetaeps(t)}_{W^{1,q}(\Omega)}}  + \normaH{R(t)} \bigr) 
  \norma v_{\Wx{1,q'}} \,.
  \non
\Esist
This means that
\Beq
  \norma{\dt\weps(t)}_{(\Wx{1,q'})^*} 
  \leq c \bigl( 1 + \normaH{\dt\chieps(t)}^2 + 
  \pier{\norma{\thetaeps(t)}_{W^{1,q}(\Omega)}}+ \normaH{R(t)} \bigr).
  \non
\Eeq
As $R\in\LQ2$ and \eqref{primastima}\pier{, \eqref{stimaBGtheta} hold},
we conclude that
\Beq
  \norma{\dt\weps}_{\L1{(\Wx{1,q'})^*}} \leq \pier{c_q} \,.
  \label{stimadtw}
\Eeq

\step
Convergence and first consequences

From \elena{\accorpa{primastima}{stimadtseconda}}, \eqref{secondastima},
\gianni{\eqref{bddalphaeps}},
\eqref{stimabeta} and \accorpa{stimaBGw}{stimaBGtheta},
we deduce that a sextuplet \pier{$(w, \theta, u,
\sigma, \chi  ,\xi)$} and a pair $(a,\ell)$ exist such~that,
for a subsequence \pier{of $\eps \searrow 0$} and for every $q\in[1,5/4)$,
the following convergence holds~true
\Bsist
  && \weps \to w
  \quad \hbox{weakly in $\LQ{4q/3}\cap\L q{\Wx{1,q}}$}
  \label{convweps}
  \\
  && \thetaeps \to \theta
  \quad \hbox{weakly in $\LQ{4q/3}\cap\L q{\Wx{1,q}}$}
  \label{convthetaeps}
  \\
  && \alphaeps(\thetaeps) \to a
  \quad \hbox{weakly star in \gianni{$\LQ\infty$}}
  \label{convalpha}
  \\
  && \ln(\thetaeps+\eps) \to \ell
  \quad \hbox{weakly in $\L2V$}
  \label{convln}
  \\
  && \ueps \to u
  \quad \hbox{weakly star in $\H2\Vp\cap\W{1,\infty}H\cap\L\infty V$}
  \qquad
  \label{convueps}
  \\
  && \sigmaeps \to \sigma
  \quad \hbox{weakly star in $\L\infty\HH$}
  \label{convsigmaeps}
  \\
  && \chieps \to \chi
  \quad \hbox{weakly star in $\H1H\cap\L2{\pcol{W_0}}$}
  \label{convchieps}
  \\
  && \eps \dt\nabla\chieps \to 0
  \quad \hbox{strongly in $\L2\HH$}
  \label{convperturb}
  \\
  && \betaeps(\chieps) \to \xi
  \quad \hbox{weakly in $\L2\Vp$}.
  \label{convxieps}
\Esist
\elena{Owing} to Lemma~\ref{Energia}, 
we see that the \pier{regularity requirements \eqref{regtheta}, 
\eqref{regu}--\eqref{regchi}
stated in Definition~\ref{Defsol} and 
regarding $w, \, \theta, \, u,\, \sigma, \, \chi$ are fulfilled except for} 
the positivity for~$\theta$.
Moreover, the above convergence implies weak convergence at least in $\C0\Vp$
for $\ueps$, $\dt\ueps$ and $\chieps$
and \pier{\eqref{cauchyeps},} \accorpa{convuzeps}{convchizeps} hold,
so that the Cauchy conditions \eqref{cauchy} 
\gianni{for $u$, $\dt u$ and $\chi$ are satisfied
(while more work is needed for~$w$)}.
Furthermore, \eqref{convchieps}~implies
\Beq
  \chieps \to \chi
  \quad \hbox{weakly in $\C0V$}
  \label{convchiCzV}
\Eeq
due to the continuous embedding $\H1H\cap\L2W\subset\C0V$.
Now, we \pier{recall the compact embeddings
$V\subset H$ and $W\subset V$, as well as} the continuous embeddings
$\Wx{1,q}\subset\Lx q\subset\pier{\Lx{1}}\subset(\Wx{1,q'})^*$, 
\pier{the first one being compact. Then, thanks to 
estimate~\eqref{stimadtw} and accounting for strong} compactness results
(see, e.g., \cite[Sect.~8, Cor.~4]{Simon}),
we derive some strong \gianni{and a.e.}\ convergence (for a subsequence).
Namely, we \pier{deduce that}
\Bsist
  && \weps \to w
  \quad \hbox{strongly in $\pier{\LQ q}$ and \aeQ}
  \label{strongw}
  \\
  && \ueps \to u
  \quad \hbox{strongly in \pier{$\C1{\Vp}\cap\C0H$} \ele{and \aeQ}} 
  \label{strongu}
  \\
  && \chieps \to \chi
  \quad \hbox{strongly in $\C0H\cap\L2V$ and \aeQ}
  \label{strongchi}
\Esist
\pier{for every $q\in[1,5/4)$.} \gianni{As a consequence of~\eqref{strongchi}, 
the limits of all the nonlinear terms involving~$\chieps$, but $\betaeps(\chieps)$,
can be correctly identified.
Namely, we~have \pier{(cf.~\eqref{defgeps}--\eqref{convGeps})}
\Bsist
  && \phieps(\chieps) \to \phi(\chi)
  \quad \hbox{strongly in $\LQ p$ for $p<+\infty$ and \aeQ}
  \non
  \\
  && \quad \hbox{where $\phi=G$, $G'$, $\gamma$, $\gamma^{\pier{\,\prime}}$}.
  \label{convnonlinchieps}
\Esist
Let us comment, e.g., on the limit of~$\Geps(\chieps)$. 
Due to \eqref{strongchi} and assumption~\eqref{convGeps},
we deduce that $\Geps(\chieps)$ a.e.\ converges to~$G(\chi)$. 
Then, \eqref{convnonlinchieps} with $\phi=G$ follows for 
\eqref{disugGeps} implies that $\Geps(\chieps)$ is bounded in~$\LQ\infty$. 
\pier{In addition, as}
\eqref{strongchi} yields $\pi(\chieps)\to\pi(\chi)$, e.g., strongly in $\C0H$
since $\pi$ is \Lip\ continuous,
we} infer that \eqref{defsigma} is satisfied and~that
\Beq
  \dt\chi - \Delta\chi
  + \thetac \bigl( \xi + \pi(\chi) \bigr) + a \, G'(\chi) - \sigma\cdot\ee \, \gamma^{\pier{\,\prime}}(\chi) = 0
  \quad\hbox{in $\Vp$, \aeT}.
  \label{terzaVp}
\Eeq
\pier{Now}, just by comparison in~\eqref{terzaVp}, we deduce~that
\Beq
  \xi \in \L2H
  \label{xiLdueH}
\Eeq
i.e., the first \pier{condition in}~\eqref{regxi}.
Hence, equation \eqref{terza} is satisfied as well
once we prove that $a=\alpha(\theta)$ and $\xi=\beta(\chi)$.
This will be done \pier{in the following}.
As far as \eqref{seconda} is concerned,
we easily recover an integrated version of~it
(in~fact equivalent to \eqref{seconda} itself), namely
\Beq
  \ioT \< \dt^2 u(t), v(t) > \, dt
  + \kappa \intQ \sigma \cdot \nabla v
  = \intQ \BO \, v
  + \ioT \< \BG(t) , v(t) > \, dt
  \label{intseconda}
\Eeq
for every $v\in \L2V$.
Indeed, the analogue of \eqref{intseconda} for the approximating problem
\Beq
  \ioT \< \dt^2 \ueps(t), v(t) > \, dt
  + \kappa \intQ \sigmaeps \cdot \nabla v
  = \intQ \BOeps \, v
  + \ioT \< \BGeps(t) , v(t) > \, dt
  \non
\Eeq
holds true as well for every $v\in \L2V$, so that it suffices to recall~\accorpa{convBOeps}{convBGeps} \pier{and \eqref{convueps}--\eqref{convsigmaeps}}.

\step
\gianni{More identifications and properties}

\gianni{%
We can derive both~\eqref{defw} and positivity for~$\theta$
(we~just have $\theta\geq0$ \pier{for the moment}, as a consequence of \eqref{convthetaeps} and of $\thetaeps\geq0$ for~$\eps>0$), \pier{as well as
we} identify the weak limits $a$ and $\ell$ given by \accorpa{convalpha}{convln} 
as $\alpha(\theta)$ and $\ln\theta$, respectively.
\pier{Regarding the first claim}, we prove that 
\Beq
  \thetaeps \to \theta
  \quad \aeQ .
  \label{convaetheta}
\Eeq
To this aim, we use the analogous convergence for~$\weps$ (cf.~\eqref{strongw}),
the convergence a.e.\ for $\Geps(\chieps)$ just remarked,
and the uniform bounds and convergence properties of the approximating nonlinearities\pier{:
see \eqref{convalphaeps}, \eqref{disugGeps}, \accorpa{parabeps}{defphieps} and
also} \accorpa{hptreD}{hpbdd}. \pier{For two different indices $\eps, \, \eps'$ of the subsequence 
we~have 
\Bsist
  && w_{\eps} - w_{\eps'} 
  \non
  \\
  && {} = \cz \thetaeps + \bigl( \alphaeps(\thetaeps) - (\thetaeps+\eps) \alphaeps'(\thetaeps) \bigr) \Geps(\chieps)
  \non
  \\
  && \quad {}
  - \cz\theta_{\eps'} - \bigl( \alphaeps(\theta_{\eps'}) - (\theta_{\eps'}+\eps) \alphaeps'(\theta_{\eps'}) \bigr) \Geps(\chieps)
  \non
  \\
  && {} + \bigl( 
    \alphaeps(\theta_{\eps'}) - (\theta_{\eps'}+\eps) \alphaeps'(\theta_{\eps'})
    - \alpha_{\eps'}(\theta_{\eps'}) + (\theta_{\eps'} +\eps') \alpha_{\eps'}'(\theta_{\eps'})
  \bigr) \Geps(\chieps)
  \non
  \\
  && {} + \bigl( \alpha_{\eps'}(\theta_{\eps'}) - (\theta_{\eps'} +\eps')
  \alpha_{\eps'}'(\theta_{\eps'}) 
  \bigr)
  \bigl( \Geps(\chieps) - G_{\eps'}(\chi_{\eps'}) \bigr) .
  \non
\Esist
Thus, we deduce that 
\Bsist
  && |\weps - w_{\eps'} |
  \non
  \\
  && {} \geq \deltastar |\thetaeps - \theta_{\eps'}|
  - \sup_{r\geq0} \, \bigl| 
    \alphaeps(r) - (r+\eps) \alphaeps'(r)
    - \alpha_{\eps'}(r) + (r + \eps') \alpha'_{\eps'}(r) 
  \bigr| \, \sup_{s\in\erre} \Geps(s)
  \non
  \\
  && {} - {} \pier{\sup_{r\geq0}} \,\bigl| \alpha_{\eps'}(r) 
  - (r + \eps') \alpha'_{\eps'}(r) \bigr|
  \, \bigl| \Geps(\chieps) - G_{\eps'}(\chi_{\eps'}) \bigr| 
  \non
\Esist
which implies that $\{\thetaeps \} $ is a Cauchy sequence and consequently converges 
almost everywhere in $Q$ to some measurable function $\Theta$. Then, using \eqref{convthetaeps}
and the Egorov theorem, it is not difficult to find out that $\Theta = \theta$ and 
\Beq
  \thetaeps \to \theta
  \quad \hbox{strongly in $\pier{\LQ q}$ for every $q\in[1,5/4)$.}
  \label{strongtheta}
\Eeq
In particular, \eqref{convaetheta} follows. Moreover, owing to~\eqref{bddalphaeps}, for every $p<+\infty$ a strong convergence in $\LQ p$ to the correct limits holds true  
for all the nonlinear terms involving~$\alphaeps$, like
$\alphaeps(\thetaeps)$ and $(\thetaeps+\eps)\alphaeps'(\thetaeps)$.
Therefore, \eqref{defw} comes out as a consequence.
Next,}} we prove that $\theta>0$ \aeQ\ 
and that the weak limit $\ell$ given by \eqref{convln} coincides with~$\ln\theta$.
To this aim, we recall the bound for $\ln(\thetaeps+\eps)$ given by~\eqref{primastima}.
\pier{Thanks to} \eqref{convaetheta} and to the Fatou lemma,
we deduce that $\ln\theta\in\LQ1$, whence $\theta>0$ \aeQ.
More precisely, we have $\iO|\ln(\thetaeps(t)+\eps)|\leq c$ \aat, whence also $\iO|\ln\theta(t)|\leq c$,
i.e., $\ln\theta\in\L\infty\Luno$ \elena{and $\ell=\ln\theta$ as well.} 

\elena{Now, we aim to 
identify $\xi$ in \eqref{terzaVp} 
\pier{as a selection from $\beta(\chi)$} (see \eqref{regxi}). 
We introduce three nonnegative functionals 
on $H$, $V$ and $\L2V$, respectively, by setting
(being understood that the integrals are possibly infinite)
\Bsist
  && \jH(v) := \iO \Funo(v)
  \quad \hbox{for $v\in H$} 
  \aand
  \jV(v) := \iO \Funo(v)
  \quad \hbox{for $v\in V$} 
  \non
  \\
  && \JV(v) := \ioT \jV(v(t)) \, dt 
  = \intQ \Funo(v)
  \quad \hbox{for $v\in\L2V$} 
  \non
\Esist
(thus, $\jV$ is the restriction of $\jH$ to~$V$).
Their subdifferentials are (possibly multi-valued) maps from the above spaces
\pier{to} the corresponding dual spaces (see, e.g., \cite[p.~52]{Barbu}).
The \characteriz ation of $\jH$ and $\JV$ we use can be obtained by applying, e.g., 
\cite[Ex.~3 and Prop.~2.8]{Barbu} and adapting the argument, respectively,
while the property of $\partial\jV$ we are going to mention can be found, e.g., in \cite[Prop.~2.5]{BCGG}.
So, by identifying $H^*$ with $H$ and $(\L2V)^*$ with $\L2\Vp$ as usual, we have
\Bsist
  && \hbox{for $v\in H$ and $\vstar\in H$}
  \non
  \\
  && \quad
  \vstar \in \partial \jH(v) 
  \quad \hbox{if and only if} \quad
  \vstar(x) \in \partial\Funo(v(x)) = \beta(v(x)) \quad \aaO
  \non
  \\
  && \hbox{for $v\in\L2V$ and $\vstar\in\L2\Vp$}
  \non
  \\
  && \quad
  \vstar \in \partial \JV(v) 
  \quad \hbox{if and only if} \quad
  \vstar(t) \in \partial\jV(v(t)) \quad \aat
  \non
  \\
  && \hbox{for $v\in V$ and $\vstar\in H $}
  \non
  \\
  && \quad
  \vstar \in \partial \jV(v) 
  \quad \hbox{if and only if} \quad
  \vstar \in \partial \jH(v) .
  \non
\Esist
By recalling that $\chi\in\L2V$ and $\xi\in\L2H$ (cf.~\eqref{xiLdueH}),
observing that $\xi\in\beta(\chi)$ \aeQ\ if and only if
$\xi(t)\in\beta(\chi(t))$ \aeO, \aat,
and combining the above statements, we deduce that
\Beq
  \xi \in \beta(\chi) \quad \aeQ 
  \quad \hbox{if and only if} \quad
  \xi \in \partial\JV(\chi) .
  \non
\Eeq
Thus, we prove that $\xi\in\partial\JV(\chi)$, i.e.
\Beq
  \intQ \Funo(\chi)
  + \ioT \< \xi(t) , z(t) - \chi(t) > \, dt
  \leq \intQ \Funo(z)
  \quad \hbox{for every $z\in\L2V$}
  \label{xibetastar}
\Eeq
(in fact, the above duality is an integral since $\xi\in\L2H$).
So, we fix $z\in\L2V$ and assume that $\Funo(z)\in\LQ1$, without loss of generality.
By convexity and $\betaeps=\Funoeps'$ (see~\eqref{defFeps}), we have for $\eps>0$
\Beq
  \intQ \Funoeps(\chieps) 
  + \pcol{\ioT  \<\betaeps(\chieps) , z-\chieps > (t) \, dt}{}  \leq \intQ \Funoeps(z) .
  \non
\Eeq
Moreover, the weak convergence \eqref{convxieps} is coupled 
with the strong convergence \eqref{strongchi} in the \pcol{duality paring} 
on the \lhs\ of the above inequality,
and $\Funoeps(s)\leq\Funo(s)$ for every $s\in\erre$ by~\eqref{defFeps}.
Therefore, \eqref{xibetastar} immediately follows once we prove~that
\Beq
  \intQ \Funo(\chi) \leq \liminf_{\eps\seto0} \intQ \Funoeps(\chieps) .
  \label{perxi}
\Eeq
To this end, we fix $\eps'>0$ for a while. 
By accounting for~\eqref{strongchi},
the lower semicontinuity of $\Funoepsp$ and the inequality 
$\Funoepsp(s)\leq\Funoeps(s)$ for every $s\in\erre$ and $\eps\in(0,\eps')$
(trivially from~\eqref{defFeps}),
we~obtain 
\Beq
  \intQ \Funoepsp(\chi)
  \leq \liminf_{\eps\seto0} \intQ \Funoepsp(\chieps)
  \leq \liminf_{\eps\seto0} \intQ \Funoeps(\chieps) .
  \label{perxiA}
\Eeq
Now, we let $\eps'$ vary and recall that
$\Funoepsp(s)\Neto\Funo(s)$ monotonically for every $s\in\erre$ as $\eps'\Seto0$.
Thus, the Beppo Levi monotone convergence theorem yields
\Beq
  \intQ \Funo(\chi)
  = \lim_{\eps'\seto0} \intQ \Funoepsp(\chi).
  \label{perxiB}
\Eeq 
By combining \eqref{perxiB} and \eqref{perxiA}, we obtain~\eqref{perxi}.
Therefore, even \eqref{xibetastar} is established and the proof is complete.
}

\step
Further strong convergence

\pier{In order to} \elena{pass to the limit in \eqref{primaepster} we need to 
prove that $\dt\chieps$ strongly converges to $\dt\chi$  in~$\LQ2$.}
To this end, it suffices to show~that
\Beq
  \limsup_{\eps\seto0} \intQ |\dt\chieps|^2
  \elena{{}\leq{}} \intQ |\dt\chi|^2 
  \label{perdtchi}
\Eeq
where it is understood that $\eps$ tends to zero along the subsequence satisfying
all the convergence properties just proved,
in particular~\eqref{convchieps}.
To achieve~\eqref{perdtchi}, we compute the integral on the \lhs\
by testing \eqref{terzaeps} by~$\dt\chieps$.
We~have
\Bsist
  && \intQ |\dt\chieps|^2
  = - \frac 12 \iO |\nabla\chieps(T)|^2
  +  \frac12 \iO |\nabla\chizeps|^2
  - \eps \intQ |\dt\nabla\chieps|^2
  \non
  \\
  &&\quad \pier{{}- \thetac \iO \Funoeps (\chieps (T) )
  + \thetac \iO \Funoeps (\chizeps )}
    \non
  \\
  && \quad {}
  - \thetac \intQ \pier{\pi}(\chieps) \, \dt\chieps
  - \intQ \alpha(\thetaeps) \Geps'(\chieps) \dt\chieps
  + \intQ \sigmaeps \cdot \ee \, \geps^{\pier{\,\prime}}(\chieps) \dt\chieps \,.
  \label{pier12}
\Esist
\pier{%
As $\chieps(T) \to \chi(T)$ weakly in $V$ and strongly in $H$ 
(cf.~\eqref{convchiCzV} and \eqref{strongchi}), we have 
$$
 \frac 12 \iO |\nabla\chi(T)|^2 \leq \liminf_{\eps\seto0} \, 
 \frac 12 \iO |\nabla\chieps(T)|^2
$$ 
and, \pcgg{arguing as in the proof of \eqref{perxi},} we can show that
\Beq
  \thetac  \iO \Funo(\chi(T)) \leq \liminf_{\eps\seto0}\, \thetac \iO \Funoeps(\chieps(T)) .
  \non
\Eeq
Moreover, in view of \eqref{convchizeps} we also infer that 
$$
   \lim_{\eps\seto0} \tonde{\frac12 \iO |\nabla\chizeps|^2
    + \thetac \iO \Funoeps (\chizeps )} =  \frac12 \iO |\nabla\chiz|^2
    + \thetac \iO \Funo (\chiz ).  
$$
Hence, using the weak convergence of $\dt \chieps $ in $L^2(Q) $ and 
the strong convergences of $\pi(\chieps)$ and 
$\alpha(\thetaeps) \Geps'(\chieps) $ in $L^2(Q)$, 
from \eqref{pier12} it is straightforward to deduce that
\Bsist
  && \limsup_{\eps\seto0} \intQ |\dt\chieps|^2
  \elena{{}\leq{}} - \frac 12 \iO |\nabla\chi(T)|^2  - \thetac  \iO \Funo(\chi(T)) 
  \non
  \\
  && \quad {}
   + \frac12\iO |\nabla\chiz|^2 + \thetac \iO \Funo (\chiz )  
  - \thetac \intQ \pi (\chi) \, \dt\chi
  \non
  \\
  && \quad {}  
  - \intQ \alpha(\theta) G'(\chi) \dt\chi
  + \limsup_{\eps\seto0} \intQ \sigmaeps \cdot \ee \, \geps^{\pier{\,\prime}}(\chieps) \dt\chieps \,.
  \label{perdtchiA}
\Esist
Unfortunately,} the last term of \eqref{perdtchiA} cannot be immediately identified
since it couples two weakly convergent factors.
In order to estimate it,
we compute the integral \pier{with the help of}~\eqref{energiaeps} written with $t=T$,
combine weak convergence for the terms involving the solution
and strong convergence for the data (see~\accorpa{convBOeps}{convchizeps}),
and use \pier{weak semicontinuity as before. In particular, \pcgg{due to
\eqref{convueps} and \eqref{strongu} we note that}
$\dt \ueps (T) \to  \dt u(T)$ weakly in~$H$
and $\ueps (T) \to  u(T)$ weakly in~$V$ 
\pcgg{(whence easily $\sigmaeps (T) \to \sigma (T) $ weakly in~$\HH$ as well)}.}
Finally, we account for identity~\eqref{energia}.
We obtain
\Bsist
  && \limsup_{\eps\seto0} \intQ \sigmaeps \cdot \ee \, \geps^{\pier{\,\prime}}(\chieps) \dt\chieps
  \non
  \\
  && \elena{{}\leq{}}- \frac \kappa 2 \iO |\dt u(T)|^2
  + \frac \kappa 2 \iO |\uz'|^2
  - \frac 12 \iO |\sigma(T)|^2
  + \frac 12 \iO |\sigma(0)|^2
  \non
  \\
  && + \kappa \intQ \BO \dt u
  + \kappa \< \BG(T) , u(T) > 
  - \kappa \< \BG(0) , \uz >
  - \kappa \ioT \< \dt\BG(s) , u(s) > \ds
  \non
  \\
  && = \intQ \sigma \cdot \ee \, \gamma^{\pier{\,\prime}}(\chi) \dt\chi.
  \non
\Esist
Hence, \pier{observing also that 
\Bsist
  &&- \frac 12 \iO |\nabla\chi(T)|^2  - \thetac  \iO \Funo(\chi(T)) 
  \non
   \\
   &&\qquad {}+ \frac12\iO |\nabla\chiz|^2 + \thetac \iO \Funo (\chiz )  
   = {}
     - \int_Q ( -\Delta\chi  + \thetac \, \xi )\,  \dt \chi  
  \non
\Esist thanks to integration by parts, \eqref{regxi}, and the chain rule shown e.g. in \cite[Lemme~3.3, p.~73]{Brezis}, the inequality \eqref{perdtchiA} entails
\Bsist
  && \limsup_{\eps\seto0} \intQ |\dt\chieps|^2 \leq
  - \int_Q ( - \Delta\chi  + \thetac \, \xi )\,  \dt \chi  
  \non
  \\
  && \quad {}
  - \thetac \intQ \pi (\chi) \, \dt\chi
  - \intQ \alpha(\theta) G'(\chi) \dt\chi
  + \intQ \sigma \cdot \ee \, \gamma^{\pier{\,\prime}}(\chi) \dt\chi.
  \label{perdtchiB}
\Esist
On} the other hand, by testing \eqref{terza} by $\dt\chi$
and integrating over~$Q$,
one immediately sees that the \rhs\ of \eqref{perdtchiB}
is precisely $\intQ|\dt\chi|^2$.
Therefore, \eqref{perdtchi}~is proved.

\step
\pier{End of the proof}

\pier{Now, we can take the limit in \eqref{primaepster}. 
Taking $v \in W^{1,q'}(\Omega)$, with $q'$ as in \eqref{regwt}, 
and integrating from $0$ to $t\in (0,T]$, 
thanks to \eqref{cauchyeps} we have that
\Beq
 \langle \weps (t) ,  v \rangle = \langle \wzeps ,  v \rangle  
  + \int_{Q_t} \alphaeps(\thetaeps) \, \Gpeps(\chieps) \, \dt\chi_\eps  \, v
  - \int_{Q_t} \nabla\thetaeps \cdot \nabla v
  + \int_{Q_t}\bigl( R + |\dt\chieps|^2 \bigr) v . 
  \non
\Eeq
\pcgg{We observe that \eqref{strongw} yields $ \weps (t) \to w(t) $ strongly in $\Lx q$ \aat.
On the other hand,} owing to \pcgg{our} convergence properties and \eqref{convwzeps}, 
\pcgg{the above \rhs\ converges to the expected limit for every $t\in[0,T].$
Therefore, it turns out that}
\Beq
 \langle w (t) ,  v \rangle = \langle \wz ,  v \rangle  
  + \int_{Q_t} \alpha(\theta) \, G' (\chi) \, \dt\chi  \, v
  - \int_{Q_t}\nabla\theta \cdot \nabla v
  + \int_0^t \langle R + |\dt\chieps|^2 , v \rangle  
  \label{pier13}
\Eeq
\pcgg{\aat. In particular, $w$  belongs to $\C0{(W^{1,q'}(\Omega))^*}$  and
the initial condition for $w$ in \eqref{cauchy} is satisfied. Furthermore,} by differentiating \eqref{pier13} with respect to $t$, we finally recover \eqref{primater} and the regularity 
\eqref{regwt} for $\dt w$.}

\pier{About the $L^\infty(0,T;L^1(\Omega))$-regularity of $w$ 
(cf.~\eqref{regtheta2}), \eqref{strongw} implies that $ \weps \to w $ in $L^1 (0,T;L^1(\Omega))$, whence
$$
\| \weps (t)\|_{L^1(\Omega)} \ \to \ \| w (t)\|_{L^1(\Omega)} \quad\hbox{ for a.e. } \, t\in(0,T), \, \hbox{ at least for a subsequence.}
$$
Then, recalling \eqref{primastima} we infer that
$$ \| w (t)\|_{L^1(\Omega)}\leq \sup_{\eps \in (0,1)}  \| \weps \|_{L^\infty(0,T;L^1(\Omega))} \leq c,$$
and consequently $w \in L^\infty(0,T;L^1(\Omega)).$ The same property can be deduced for $\theta$,  so that \eqref{regtheta2} holds. The proof of Theorem~\ref{Esistenza} is then complete.}


\Begin{thebibliography}{10}

\elena{
\bibitem{Auricchio} 
F. Auricchio, L. Petrini, \pier{A three-dimensional model describing stress-temperature 
induced solid phase transformations: solution algorithm and boundary value problems, 
{\it  Internat. J. Numer. Methods Engrg.} {\bf 61}} (2004) 807-836.}

\elena{
\bibitem{AuricchioB} 
F. Auricchio, E. Bonetti, A new ``flexible'' 3D macroscopic model 
for shape memory alloys, {\it Discr. Cont. Dynam. Syst Ser. S}  
{\bf 6} (2013) 277-291.
}

\bibitem{Baiocchi}
C.  Baiocchi, 
Soluzioni ordinarie e generalizzate
del problema di Cauchy per equazioni differenziali astratte
lineari del secondo ordine in spazi di Hilbert, 
{\it Ricerche Mat.} {\bf 16} (1967) 27-95.

\bibitem{Barbu}
V. Barbu,
``Nonlinear semigroups and differential equations in Banach spaces'',
Noord\-hoff,
Leyden,
1976.

\bibitem{BCGG}
V. Barbu, P. Colli, G. Gilardi, M. Grasselli,
Existence, uniqueness, and longtime behavior
for a nonlinear Volterra integrodifferential equation,
{\it Differential Integral Equations} {\bf 13} (2000) 1233-1262. 

\elena{
\bibitem{bfg} V. Berti, M. Fabrizio, D. Grandi, Phase transitions 
in shape memory alloys: a non isothermal Ginzburg-Landau model, 
\pier{{\it Phys. D} {\bf 239} (2010) 95-102.}}

\ele{
\bibitem{bfg2} V. Berti, M. Fabrizio, D. Grandi, Hysteresis and phase transitions for 
one-dimensional and three-dimensional models in shape memory alloys, 
{\it J. Math. Phys.} {\bf 51} (2010), 13~pp.}

\bibitem{BoGa}
L. Boccardo, T. Gallou\"et, 
Non-linear elliptic and parabolic equations involving measure data,
{\it J. Funct. Anal.} {\bf 87} (1989) 149-169.

\pier{\bibitem{B1}
E. Bonetti, 
Global solution to a nonlinear phase transition model with dissipation,
{\it  Adv. Math. Sci. Appl.} {\bf 12} (2002) 355-376.}

\pier{\bibitem{b03-1}
E. Bonetti,
\newblock Global solvability of a dissipative {F}r\'emond model for shape
  memory alloys. {I}. {M}athematical formulation and uniqueness,
\newblock {\em Quart. Appl. Math.}  {\bf 61} (2003) 759-781.}

\pier{\bibitem{b03-2}
E. Bonetti,
\newblock Global solvability of a dissipative {F}r\'emond model for shape
  memory alloys. {II}. {E}xistence,
\newblock {\em Quart. Appl. Math.} {\bf 62} (2004) 53-76.}

\pier{\bibitem{bcl}
E. Bonetti, P. Colli, Ph. Lauren\c cot, Global existence 
for a hydrogen storage model with full energy balance, 
{\it Nonlinear Anal.}  {\bf 75} (2012) 3558-3573.}

\pier{\bibitem{bfl06}
E. Bonetti, M. Fr{\'e}mond, Ch. Lexcellent,
\newblock Global existence and uniqueness for a thermomechanical model for
  shape memory alloys with partition of the strain,
\newblock {\em Math. Mech. Solids} {\bf 11} (2006) 251-275.}


\pier{\bibitem{BFL}
G. Bonfanti, M. Fr\'emond, F. Luterotti,
Global solution to a nonlinear system for irreversible phase changes,
{\it  Adv. Math. Sci. Appl.}  {\bf 10} (2000) 1 - 24.}

\bibitem{Brezis}
H. Brezis,
``Op\'erateurs maximaux monotones et semi-groupes de contractions
dans les espaces de Hilbert'',
{\it North-Holland Math. Stud.}
{\bf 5},
North-Holland,
Amsterdam,
1973.

\pier{\bibitem{c95}
P. Colli,
\newblock Global existence for the three-dimensional {F}r\'emond model of shape
  memory alloys,
\newblock {\em Nonlinear Anal.}  {\bf 24} (1995) 1565-1579.}

\pier{\bibitem{cfrs06}
P. Colli, M. Fr{\'e}mond, E. Rocca, K. Shirakawa,
\newblock Attractors for a three-dimensional thermo-mechanical model of shape
  memory alloys,
\newblock {\em Chinese Ann. Math. Ser. B} {\bf 27} (2006) 683-700.}

\bibitem{CGG}
P. Colli, G. Gilardi, M. Grasselli,
Well-posedness of the weak formulation for the phase-field model with memory,
{\it Adv. Differential Equations} {\bf 2} (1997) 487-508. 

\pier{\bibitem{cls00}
P. Colli, Ph. Lauren{\c{c}}ot, U. Stefanelli,
\newblock Long-time behavior for the full one-dimensional {F}r\'emond model for
  shape memory alloys,
\newblock {\em Contin. Mech. Thermodyn.} 12(6):423--433, 2000.}

\pier{\bibitem{CLSS}
P. Colli, F. Luterotti, G. Schimperna, and U. Stefanelli,
Global existence for a class of generalized systems for irreversible phase changes,
{\it NoDEA Nonlinear Differential Equations Appl.}  {\bf 9} (2002) 255-276.}

\pier{\bibitem{cs93}
P.~Colli, J.~Sprekels.
\newblock Positivity of temperature in the general {F}r\'emond model for shape
  memory alloys,
\newblock {\em Contin. Mech. Thermodyn.} {\bf 5} (1993) 255-264.}

\elena{\bibitem{dfg} 
F. Daghia, M. Fabrizio, D. Grandi, A non isothermal Ginzburg-Landau 
model for phase transitions in shape memory alloys,
{\it Meccanica} {\bf 45} (2010) 797-807.}

\bibitem{DauLio}
R. Dautray and J. L. Lions, 
``Analyse math\'ematique et calcul num\'erique pour 
les sciences et les techniques'', Tome~3,
Masson, Paris, 1985.

\ele{\bibitem{dfmz} R.P. Dhote, M. Fabrizio, R.N.V. Melnik, J. Zu, 
Hysteresis phenomena in shape memory alloys
by non-isothermal Ginzburg-Landau models, {\em Commun. Nonlinear 
Sci. Numer. Simul.}  {\bf 18} (2013) 2549-2561.}

\ele{\bibitem{FaPe} M. Fabrizio, 
M. Pecoraro, Phase transitions and thermodynamics for the shape
memory alloy AuZn, {\it Meccanica}, published on line (2013), 
doi 10.1007/s11012-013-9701-3}

\elena{\bibitem{Fremond}
M. Fr\'emond,
``Non-smooth Thermomechanics'',
Springer-Verlag, Berlin, 2002.}

\pier{\bibitem{frem2}
M. Fr\'emond,
``Phase Change in Mechanics'', 
{\it Lect. Notes Unione Mat. Ital.} {\bf 13}, 
Springer, Heidelberg; UMI, Bologna, 
2012.}

\bibitem{GiRo}
G. Gilardi, E. Rocca,
Well-posedness and long-time behaviour for a singular phase
field system of conserved type,
{\it IMA J. Appl. Math.} {\bf 72} (2007) 498-530.

\bibitem{Jerome}
J.W. Jerome, ``Approximation of Nonlinear Evolution Systems'', 
{\it Math. Sci. Engrg.}  {\bf 164}, Academic Press, Orlando,
1983.

\pier{\bibitem{LSS}
Ph. Lauren\c cot, G. Schimperna, and U. Stefanelli,
Global existence of a strong solution to the one-dimensional 
full model for irreversible phase transitions,
{\it  J. Math. Anal. Appl.}  (\bf 271) (2002) 426-442.}

\pier{\bibitem{LS}
F. Luterotti, U. Stefanelli,
Existence result for the one-dimensional full model of phase transitions,
{\it  Z. Anal. Anwendungen}  {\bf 21} (2002) 335-350.}

\pier{\bibitem{Rou}
T. Roub{\'{\i}}{\v{c}}ek,
Modelling of thermodynamics of martensitic transformation in
              shape-memory alloys, 
{\it Discrete Contin. Dyn. Syst.}, 
``Dynamical Systems and Differential Equations. Proceedings of
              the 6th AIMS International Conference, suppl.'',
 (2007) 892-902.}

\bibitem{Simon}
J. Simon,
{Compact sets in the space $L^p(0,T; B)$},
{\it Ann. Mat. Pura Appl.~(4)} {\bf 146} (1987) 65--96.

\End{thebibliography}

\End{document}

\bye